\documentclass{alggeom}

\usepackage{mathrsfs}
\usepackage{amscd}
\usepackage{amsmath}
\usepackage{amssymb}
\usepackage{latexsym}
\usepackage{comment}
\usepackage{amscd}
\usepackage{wasysym}  
\usepackage{tikz}
\usetikzlibrary{matrix,arrows}
\usepackage{tikz-cd}
\usepackage{appendix}
\usepackage[nice]{nicefrac}
\usepackage{fix-cm} 
\usepackage{comment}
\usepackage{amscd}
\usepackage{wasysym}  
\usepackage{tikz}
\usetikzlibrary{matrix,arrows}
\usepackage{tikz-cd}
\usepackage{appendix}
\usepackage[nice]{nicefrac}

\usepackage{marvosym}
\usepackage{tikz}
\usepackage[all,cmtip]{xy}
\usetikzlibrary{matrix,arrows}

\usepackage[pagebackref,hyperindex,linktocpage=true]{hyperref}
\hypersetup{
    colorlinks,
    linkcolor={blue!35!black},
    citecolor={blue!35!black},
    urlcolor={blue!35!black}
}

\usepackage{color}

\usetikzlibrary{decorations.markings}

\makeatletter
\tikzcdset{
  open/.code     = {\tikzcdset{hook, circled};},
  closed/.code   = {\tikzcdset{hook, slashed};},
  open'/.code    = {\tikzcdset{hook', circled};},
  closed'/.code  = {\tikzcdset{hook', slashed};},
  circled/.code  = {\tikzcdset{markwith = {\draw (0,0) circle (.375ex);}};},
  slashed/.code  = {\tikzcdset{markwith = {\draw[-] (-.4ex,-.4ex) -- (.4ex,.4ex);}};},
  markwith/.code ={
    \pgfutil@ifundefined%
    {tikz@library@decorations.markings@loaded}%
    {\pgfutil@packageerror{tikz-cd}{You need to say %
      \string\usetikzlibrary{decorations.markings} to use arrows with markings}{}}{}%
    \pgfkeysalso{/tikz/postaction = {
      /tikz/decorate,
      /tikz/decoration={markings, mark = at position 0.5 with {#1}}}
    }
  },
}
\makeatother

\renewcommand{\geq}{\geqslant}
\renewcommand{\leq}{\leqslant}
\renewcommand{\ge}{\geqslant}

\newtheorem{thm}{Theorem}[section]

\newtheorem{propo}[thm]{Proposition}
\newtheorem{lem}[thm]{Lemma}
\newtheorem{sublem}[thm]{Sublemma}
\newtheorem{lem-def}[thm]{Lemma-Definition}
\newtheorem{cor}[thm]{Corollary}
\newtheorem{conject}[thm]{Conjecture}
\newtheorem{propert}[thm]{Properties}
\newtheorem{observ}[thm]{Observation}

\newtheorem{fac}[thm]{Fact}
\newtheorem{ex}[thm]{Example}

\theoremstyle{definition}
\newtheorem*{ack}{Acknowledgements}

\newtheorem{rmk}[thm]{Remark}
\newtheorem{dfn}[thm]{Definition}
\newtheorem{quest}[thm]{Question}
\newtheorem{expec}[thm]{Expectation}
\newtheorem*{abs}{Abstract}

\numberwithin{equation}{section}

\newcommand{\nc}{\newcommand}

\nc{\theo}{\begin{thm}} \nc{\xtheo}{\end{thm}}
\nc{\prop}{\begin{propo}} \nc{\xprop}{\end{propo}}
\nc{\lemm}{\begin{lem}} \nc{\xlemm}{\end{lem}}
\nc{\sublemm}{\begin{sublem}} \nc{\xsublemm}{\end{sublem}}
\nc{\lemmdefi}{\begin{lem-def}} \nc{\xlemmdefi}{\end{lem-def}}
\nc{\coro}{\begin{cor}} \nc{\xcoro}{\end{cor}}
\nc{\conj}{\begin{conject}} \nc{\xconj}{\end{conject}}
\nc{\proper}{\begin{propert}} \nc{\xproper}{\end{propert}}
\nc{\obse}{\begin{observ}} \nc{\xobse}{\end{observ}}
\nc{\ques}{\begin{quest}} \nc{\xques}{\end{quest}}

\nc{\fact}{\begin{fac}} \nc{\xfact}{\end{fac}}
\nc{\expe}{\begin{expec}} \nc{\xexpe}{\end{expec}}

\nc{\ackn}{\begin{ack}} \nc{\xackn}{\end{ack}}
\nc{\exam}{\begin{ex}} \nc{\xexam}{\end{ex}}
\nc{\rema}{\begin{rmk}} \nc{\xrema}{\end{rmk}}
\nc{\defi}{\begin{dfn}} \nc{\xdefi}{\end{dfn}}
\nc{\abst}{\begin{abs}} \nc{\xabst}{\end{abs}}

\nc{\pf}{\begin{proof}} \nc{\xpf}{\end{proof}}

\nc{\on}{\operatorname}
\nc{\fraka}{{\mathfrak a}} \nc{\bba}{{\mathbf a}}
\nc{\frakb}{{\mathfrak b}}
\nc{\frakc}{{\mathfrak c}}
\nc{\frakd}{{\mathfrak d}}
\nc{\frake}{{\mathfrak e}}
\nc{\frakf}{{\mathfrak f}}
\nc{\frakg}{{\mathfrak g}}
\nc{\frakh}{{\mathfrak h}}
\nc{\fraki}{{\mathfrak i}}
\nc{\frakj}{{\mathfrak j}}
\nc{\frakk}{{\mathfrak k}}
\nc{\frakl}{{\mathfrak l}}
\nc{\frakm}{{\mathfrak m}}
\nc{\frakn}{{\mathfrak n}}
\nc{\frako}{{\mathfrak o}}
\nc{\frakp}{{\mathfrak p}}
\nc{\frakq}{{\mathfrak q}}
\nc{\frakr}{{\mathfrak r}}
\nc{\fraks}{{\mathfrak s}}
\nc{\frakt}{{\mathfrak t}}
\nc{\fraku}{{\mathfrak u}}
\nc{\frakv}{{\mathfrak v}}
\nc{\frakw}{{\mathfrak w}}
\nc{\frakx}{{\mathfrak x}}
\nc{\fraky}{{\mathfrak y}}
\nc{\frakz}{{\mathfrak z}}
\nc{\frakA}{{\mathfrak A}}
\nc{\frakB}{{\mathfrak B}}
\nc{\frakC}{{\mathfrak C}}
\nc{\frakD}{{\mathfrak D}}
\nc{\frakE}{{\mathfrak E}}
\nc{\frakF}{{\mathfrak F}}
\nc{\frakG}{{\mathfrak G}}
\nc{\frakH}{{\mathfrak H}}
\nc{\frakI}{{\mathfrak I}}
\nc{\frakJ}{{\mathfrak J}}
\nc{\frakK}{{\mathfrak K}}
\nc{\frakL}{{\mathfrak L}}
\nc{\frakM}{{\mathfrak M}}
\nc{\frakN}{{\mathfrak N}}
\nc{\frakO}{{\mathfrak O}}
\nc{\frakP}{{\mathfrak P}}
\nc{\frakQ}{{\mathfrak Q}}
\nc{\frakR}{{\mathfrak R}}
\nc{\frakS}{{\mathfrak S}}
\nc{\frakT}{{\mathfrak T}}
\nc{\frakU}{{\mathfrak U}}
\nc{\frakV}{{\mathfrak V}}
\nc{\frakW}{{\mathfrak W}}
\nc{\frakX}{{\mathfrak X}}
\nc{\frakY}{{\mathfrak Y}}
\nc{\frakZ}{{\mathfrak Z}}
\nc{\bbA}{{\mathbb A}}
\nc{\bbB}{{\mathbb B}}
\nc{\bbC}{{\mathbb C}}
\nc{\bbD}{{\mathbb D}}
\nc{\bbE}{{\mathbb E}}
\nc{\bbF}{{\mathbb F}} \nc{\bbf}{{\mathbf f}}
\nc{\bbG}{{\mathbb G}}
\nc{\bbH}{{\mathbb H}}
\nc{\bbI}{{\mathbb I}}
\nc{\bbJ}{{\mathbb J}}
\nc{\bbK}{{\mathbb K}}
\nc{\bbL}{{\mathbb L}}
\nc{\bbM}{{\mathbb M}}
\nc{\bbN}{{\mathbb N}}
\nc{\bbO}{{\mathbb O}}
\nc{\bbP}{{\mathbb P}}
\nc{\bbQ}{{\mathbb Q}}
\nc{\bbR}{{\mathbb R}}
\nc{\bbS}{{\mathbb S}}
\nc{\bbT}{{\mathbb T}}
\nc{\bbU}{{\mathbb U}}
\nc{\bbV}{{\mathbb V}}
\nc{\bbW}{{\mathbb W}}
\nc{\bbX}{{\mathbb X}}
\nc{\bbY}{{\mathbb Y}}
\nc{\bbZ}{{\mathbb Z}}
\nc{\calA}{{\mathcal A}}
\nc{\calB}{{\mathcal B}}
\nc{\calC}{{\mathcal C}}
\nc{\calD}{{\mathcal D}}
\nc{\calE}{{\mathcal E}}
\nc{\calF}{{\mathcal F}}
\nc{\calG}{{\mathcal G}}
\nc{\calH}{{\mathcal H}}
\nc{\calI}{{\mathcal I}}
\nc{\calJ}{{\mathcal J}}
\nc{\calK}{{\mathcal K}}
\nc{\calL}{{\mathcal L}}
\nc{\calM}{{\mathcal M}}
\nc{\calN}{{\mathcal N}}
\nc{\calO}{{\mathcal O}}
\nc{\calP}{{\mathcal P}}
\nc{\calQ}{{\mathcal Q}}
\nc{\calR}{{\mathcal R}}
\nc{\calS}{{\mathcal S}}
\nc{\calT}{{\mathcal T}}
\nc{\calU}{{\mathcal U}}
\nc{\calV}{{\mathcal V}}
\nc{\calW}{{\mathcal W}}
\nc{\calX}{{\mathcal X}}
\nc{\calY}{{\mathcal Y}}
\nc{\calZ}{{\mathcal Z}}

\nc{\scrA}{{\mathscr A}}
\nc{\scrE}{{\mathscr E}}
\nc{\scrR}{{\mathscr R}}

\nc{\Bmu}{\mbox{$\raisebox{-0.59ex}{$l$}\hspace{-0.18em}\mu\hspace{-0.88em}\raisebox{-0.98ex}{\scalebox{2}{$\color{white}.$}}\hspace{-0.416em}\raisebox{+0.88ex}{$\color{white}.$}\hspace{0.46em}$}{}}

\nc{\bnu}{{\bar{ \nu}}}

\nc{\olO}{\bar{\calO}}

\nc{\al}{{\alpha}} 
\nc{\be}{{\beta}}
\nc{\ga}{{\gamma}} \nc{\Ga}{{\Gamma}}
 \nc{\hGa}{\hat{\Gamma}}
\nc{\ve}{{\varepsilon}} 
\nc{\la}{{\lambda}} \nc{\La}{{\Lambda}}
\nc{\om}{\omega} \nc{\Om}{\Omega} 
\nc{\sig}{{\sigma}} \nc{\Sig}{{\Sigma}}

\nc{\tnb}{\psi_{\rm tame}}
\nc{\oM}{\overline{{M}}}
\nc{\op}{{\on{op}}}
\nc{\ad}{{\on{ad}}}
\nc{\alg}{{\on{alg}}}
\nc{\Ad}{{\on{Ad}}}
\nc{\Adm}{{\on{Adm}}} \nc{\aff}{{\on{aff}}}
\nc{\Aut}{{\on{Aut}}}
\nc{\Bun}{{\on{Bun}}}
\nc{\cha}{{\on{char}}}
\nc{\der}{{\on{der}}}
\nc{\Der}{{\on{Der}}}
\nc{\diag}{{\on{diag}}}
\nc{\End}{{\on{End}}}
\nc{\Fl}{{\calF\!\ell}}
\nc{\Tr}{{\on{Transp}}}
\nc{\TR}{{\calT\!\calR}}
\nc{\Gal}{{\on{Gal}}}
\nc{\Gr}{{\on{Gr}}}
\nc{\rH}{{\on{H}}}
\nc{\Hom}{{\on{Hom}}}
\nc{\IC}{{\on{IC}}}
\nc{\id}{{\on{id}}}
\nc{\Id}{{\on{Id}}}
\nc{\ind}{{\on{ind}}}
\nc{\Ind}{{\on{Ind}}}
\nc{\Lie}{{\on{Lie}}}
\nc{\Pic}{{\on{Pic}}}
\nc{\pr}{{\on{pr}}}
\nc{\Res}{{\on{Res}}}
\nc{\res}{{\on{res}}} \nc{\Sat}{{\on{Sat}}}
\nc{\s}{{\on{sc}}}
\nc{\drv}{{\on{der}}}
\nc{\sgn}{{\on{sgn}}}
\nc{\Spec}{{\on{Spec}}}\nc{\Spf}{\on{Spf}} 
\nc{\Sph}{\on{Sph}}
\nc{\St}{{\on{St}}}
\nc{\tr}{{\on{tr}}}
\nc{\Mod}{{\mathrm{-Mod}}}
\nc{\Hilb}{{\on{Hilb}}} 
\nc{\Ext}{{\on{Ext}}} 
\nc{\vs}{{\on{Vec}}}
\nc{\ev}{{\on{ev}}}
\nc{\nO}{{\breve{\calO}}}
\nc{\tS}{{\tilde{S}}}
\nc{\spe}{{\on{sp}}}
\nc{\loc}{{\on{loc}}}
\nc{\Sym}{{\on{Sym}}}
\nc{\Cone}{{\on{C}}}
\nc{\syn}{{\on{syn}}}
\nc{\reg}{{\on{reg}}}
\nc{\colim}{{\on{colim}}}
\nc{\Norm}{{\on{N}}}

\nc{\nscrR}{{\mathscr{R}^{\on{nr}}}}

\nc{\GL}{{\on{GL}}}
\nc{\U}{{\on{U}}}
\nc{\Gl}{\on{Gl}} 
\nc{\GSp}{{\on{GSp}}}
\nc{\gl}{{\frakg\frakl}}
\nc{\SL}{{\on{SL}}} 
\nc{\SU}{{\on{SU}}} 
\nc{\SO}{{\on{SO}}}
\nc{\PGL}{{\on{PGL}}}
\newcommand*{\calhom}{\mathit{\mathcal{H}{om}}}
\nc{\Conv}{{\on{Conv}}}
\nc{\Rep}{{\on{Rep}}}
\nc{\Dom}{{\on{Dom}}}
\nc{\red}{{\on{red}}}
\nc{\act}{{\on{act}}}
\nc{\nr}{{\on{nr}}}
\nc{\ctf}{{\on{ctf}}}

\nc{\str}{{\on{-}}} 
\nc{\os}{{\bar{s}}}
\nc{\oeta}{{\bar{\eta}}}

\nc{\hookto}{\hookrightarrow}
\nc{\longto}{\longrightarrow}
\nc{\leftto}{\leftarrow}
\nc{\onto}{\twoheadrightarrow}
\nc{\lonto}{\twoheadleftarrow}

\nc{\uG}{{\underline{G}}}
\nc{\uA}{{\underline{A}}}
\nc{\uS}{{\underline{S}}}
\nc{\uT}{{\underline{T}}}
\nc{\uM}{{\underline{M}}}
\nc{\uP}{{\underline{P}}}
\nc{\uB}{{\underline{B}}}
\nc{\uN}{{\underline{N}}}

\nc{\ucG}{{\underline{\calG}}}
\nc{\ucA}{{\underline{\calA}}}
\nc{\ucS}{{\underline{\calS}}}
\nc{\ucT}{{\underline{\calT}}}
\nc{\ucalM}{{\underline{\calM}}}
\nc{\ucP}{{\underline{\calP}}}
\nc{\ucalN}{{\underline{\calN}}}

\nc{\bF}{{\breve{F}}}

\nc{\oFl}{{\overline{\Fl}}} 
\nc{\bU}{{\overline{U}}}
\nc{\tGr}{{\tilde{\Gr}}}
\nc{\cGr}{\calG\! r}
\nc{\oGr}{\overline{\on{Gr}}} 
\nc{\ocGr}{\overline{\calG\! r}}
\nc{\co}{{\colon}}
\nc{\sch}[1]{(Sch/{#1})}
\nc{\HypLoc}[1]{HypLoc({#1})}

\nc{\ohtimes}{\stackrel{!}{\otimes}}
\nc{\boxtilde}{\widetilde{\boxtimes}}
\nc{\vstar}{{\varhexstar}}

\nc{\Div}{\on{Div}}
\nc{\Sht}{\on{Sht}}
\nc{\Frob}{\on{Frob}}

\nc{\x}{\times}
\nc{\bsl}{\backslash}
\nc{\algQl}{{\bar{\bbQ}_\ell}}
\nc{\sF}{{\bar{F}}}
\nc{\nF}{{\breve{F}}}
\nc{\nW}{{W^{\on{nr}}}}
\nc{\sk}{{\bar{k}}}
\nc{\cont}{\on{c}}
\nc{\Supp}{\on{Supp}}
\nc{\blt}{\bullet}  
\nc{\dom}{\on{dom}}
\nc{\scon}{{\on{sc}}} 
\nc{\Affine}{\on{Aff}} 
\nc{\nscrA}{\mathscr{A}^{\on{nr}}} 
\nc{\nfraka}{{\bbf^{\on{nr}}}}
\nc{\ran}{{\rangle}}
\nc{\lan}{{\langle}}
\nc{\bk}{{\bar{k}}}
\nc{\tF}{{\tilde{F}}}
\nc{\sS}{{\bar{S}}}
\nc{\LG}{{^\text{L}\hspace{-0.04cm}G}}
\nc{\LL}{{^\text{L}\hspace{-0.07cm}L}}
\nc{\et}{{\text{\rm \'et}}}
\nc{\inv}{{\on{inv}}}
\nc{\Hecke}{{\on{Hecke}}}
\nc{\Isom}{{\on{Isom}}}
\nc{\oSht}{{\overline{\on{Sht}}}}
\nc{\umu}{{\underline \mu}}
\nc{\AIJ}{{\calO_X[{\scriptstyle{\calI\over \calJ}}]}}
\nc{\Proj}{{\on{Proj}}}
\nc{\Bl}{{\on{Bl}}}
\nc{\Stab}{{\on{Stab}}}
\nc{\cl}{{\on{cl}}}

\nc{\Pos}{{\on{Pos}}}
\nc{\Sets}{{\on{Sets}}}
\nc{\AffSch}{{\on{AffSch}}}
\nc{\Groups}{{\on{Groups}}}
\nc{\Gpds}{{\on{Groupoids}}}
\nc{\Sch}{{\on{Sch}}}
\nc{\fl}{{\on{flat}}}

\nc{\pot}[1]{ [\hspace{-0,5mm}[ {#1} ]\hspace{-0,5mm}] }
\nc{\rpot}[1]{ (\hspace{-0,7mm}( {#1} )\hspace{-0,7mm}) }

\nc{\defined}{\hspace{0.1cm}\stackrel{\text{\tiny \rm def}}{=}\hspace{0.1cm}}

\usepackage{moresize}
\usepackage{graphicx}

\begin{document}

\title
[ALGEBRAIC MAGNETISM] 
{\resizebox{\linewidth}{!}{$~~~~~~~~~$Algebraic Magnetism$~~~~~~~~~~$}}

\shortauthors{ARNAUD MAYEUX}

\author{Arnaud Mayeux}

\email{mayeux@wisc.edu}

\address{University of Wisconsin–Madison, Madison, WI, USA}

\classification{14L15, 14L30, 20M32, 14M25, 13F65, 20G07, 20G35}
\keywords{algebraic magnetism,  magnets,  attractors associated to magnets, monoids, faces of monoids, Bialynicki-Birula decomposition, Grothendieck schemes, algebraic spaces, diagonalizable monoid schemes, diagonalizable group schemes, group schemes, reductive group schemes, parabolic groups, Levi subgroups, groups of type R, root groups, actions of group schemes on algebraic spaces, fixed-points, Z-FPR atlases}
\maketitle

\begin{center}
{\it \large With an appendix by Matthieu Romagny}
\end{center}

$~~$

\abst {
We introduce Algebraic Magnetism, a theory of attractors associated to magnets within the framework of actions of diagonalizable monoid schemes on algebraic spaces. For a diagonalizable monoid scheme $A(M)_S$ acting on an algebraic space $X$, we introduce for any submonoid $N$ of $M$ an attractor space $X^N$. We then investigate and study various aspects of attractors associated to monoids. This leads to the notion of pure magnets, which encode the combinatorial structure underlying attractor phenomena in terms of semigroups.
In the affine case, we give an explicit description of attractors in terms of graded algebra, yielding representability as closed subspaces in this case. Using descent methods, technical preliminary results, and fixed-point-reflecting atlases, we prove representability for general algebraic spaces. We establish compatibility of attractors with fiber products, base change, faces of monoids, intersections of submonoids, subgroups, dilatations, and many other operations.
 We study pure magnets and prove, and conjecture, finiteness results under suitable hypotheses. In this way, Algebraic Magnetism provides invariants to study diagonalizable actions in algebraic geometry.
 }  \xabst

$~~~$
\tableofcontents
\newpage
\section{Introduction } \label{sectionintro}

Let $M$ be a finitely generated abelian group and let $S$ be a base scheme. Let $D(M)_S $ be the associated diagonalizable group scheme (i.e. if $M = \bbZ ^r \times \prod _{i=1}^{n} \bbZ /n_i \bbZ $, then $D(M)_S = \mathbb{G} _{m,S}^r \times \prod _{i = 1} ^n \mu _{n_i , S}$ is the product of a split torus with group schemes of roots of unity.) Algebraic actions of diagonalizable group schemes appear systematically in many areas of mathematics. This article is devoted to introduce a new tool  in the general context of an arbitrary algebraic action of a diagonalizable group scheme on a scheme or on an algebraic space. The style is foundational and we in fact develop the theory for any diagonalizable monoid scheme. 

\subsection{Algebraic attractors} \label{para1} Let $S$ be a scheme. 
 Let $M$ be an abelian monoid. Let $\bbZ [M]$ be the ring $\bigoplus _{m \in M} \bbZ X^m$ where the multiplication is induced by the structure of monoid on $M$. Let $A(M)$ be $\Spec (\bbZ[M])$, it is a monoid scheme over $\Spec (\bbZ ) $. We consider the base change $A(M)_S = A(M) \times _{\Spec (\bbZ)} S$, it is called a diagonalizable monoid scheme over $S$ (if $M$ is moreover a group, then $A(M)_S$ is also denoted $D(M)_S$ and is a diagonalizable group scheme). For any submonoid $N$ of $M$, $A(M)_S$ acts canonically on $A(N)_S$.
  Let $X$ be an algebraic space over $S$ with an action $a$ of $A(M)_S$. The main idea in this paper is to introduce the following definition, for any submonoid $N$ of $M$:
\defi \label{defiintropara1} Let $X^N$ be the contravariant functor from schemes over $S$ to $\Sets$ given by 
\[ (T \to S) \mapsto \Hom _T ^{A(M)_T}(A(N)_T, X _T) \]
where $\Hom _T ^{A(M)_T}(A(N)_T, X _T)$ is the set of $A(M)_T$-equivariant $T$-morphisms from $A(N)_T$ to $ X_T =X \times _S T$. The functor $X^N$ is called the attractor associated to the submonoid $N$ under the action of $A(M)_S$ on $X$.
\xdefi

\subsection{Magnets and attractors} \label{firstheuristic}

We proceed with the notation from §\ref{para1} and assume that $X \to S$ is separated.
If $N \subset L $ are submonoids of $M$, then we have a canonical monomorphism $X^N \subset X^L$. Moreover $X^0 $ identifies with the fixed-points functor $X^{A(M)_S}$ and $X^M $ identifies with $X$. So for any submonoid $N \subset M$, we have monomorphisms $X^0 \subset X^N \subset X$. We now use the following terminology: a magnet for the action $a$ is a submonoid $N$ of $M$. A magnet $N\subset M$ is thought as something which algebraically attracts the subspace $X^N$ of $X$. In this paper, the word magnetism refers to a form of attraction. It is not electromagnetism, nor social-science magnetism.
For each attractor space $X^N$ there is a minimal magnet $E$ such that $X^N = X^E$, $E$ is called a pure magnet. We have a bijection between attractor spaces and the set $\mho (a) $ of pure magnets  (cf. Theorem \ref{efficienttheo}).
Assume now moreover that $X$ is finitely presented over $S$, Theorem \ref{strumagnet} says that the set of pure magnets $\mho (a)$ is finite (cf. Theorem \ref{strumagnet} for our precise assumptions).

\subsection{Faces of monoids}
We proceed with the notation from §\ref{firstheuristic} and assume for simplicity that $M$ is a finitely generated abelian groups. We fix a submonoid $N \subset M$. A face of $N$ is by definition a submonoid $F$ of $N$ such that the canonical $M$-graded projection $\bbZ [N] \to \bbZ [F]$ is a morphism of rings. Equivalently, $F \subset N$ is a face if for any two elements $x,y \in N$ the following equivalence holds: \[x+y \in F \Leftrightarrow x \in F \text{ and }y \in F.\] Note that the only face of the group $M$ is $M$ itself, so faces matter only in the world of monoids. Each face of $N$ contains $N^{*}$, the face of invertible elements in $N$.  Now if $F \subset N$ is a face, we get a canonical transformation of functor $\mathfrak{p}_{N,F}: X^N \to X^F$. Recall that on the other hand we have a monomorphism $X^F \to X^N$. The map $\mathfrak{p}_{N,F}: X^N \to X^F$ could be thought as a directional limit. Now if $Z \subset X^F$ is a monomorphism, we put $X^N_{F,Z} = X^N \times _{X^F} Z$ and we call it the attractor associated to $N$ with prescribed limit $Z$ relatively to the face $F$, if $F=N^*$ we also use the notation $X^N_Z$. We have a canonical monomorphism $X^N_{F,Z} \to X^N$. The concept of attractors with prescribed limits allows to reduce the fixed-points parts of attractors.

\subsection{Summary of results}
This article studies intrinsically algebraic magnetism, i.e. the formalism of attractors associated to magnets. Consequently, we prove in this paper a large number of results. We list here the most significant ones as informal slogans, with references to precise statements. \begin{enumerate}
\item 
Attractors are compatible with fiber products and base changes (cf. Proposition \ref{produitNNN}, Proposition \ref{fiberprodu} and Proposition \ref{basechangeee}).
\item Attractors preserve monoid and group structures (cf. Proposition \ref{groupmonoidstru}). 
\item Attractors are compatible with equivariant actions in a natural sense (cf. Proposition \ref{propequiv}). 
\item
Attractors are representable in many cases (cf. Theorems \ref{representableaffine} and \ref{representable}, Remark \ref{remamonooorep} and Proposition \ref{propequiv} (iii)).
\item In the affine case, attractors are representable by explicit closed subspaces and intersections of attractors correspond to intersections of magnets (cf. Theorem \ref{representableaffine} and Proposition \ref{inter}).

\item Inclusions of monoids give (mono)morphisms on attractors  and face inclusions provide retractions (cf. Remark \ref{equal0}, Fact \ref{equalM}, Fact \ref{morphismsgeneral}, Fact \ref{faceattractors}, Proposition \ref{monosep}, Corollary \ref{monoX}, Proposition \ref{retract}).
\item Equivariant morphisms of spaces induce morphisms on attractors; moreover closed immersions give closed immersions, open immersions give open immersions, smooth morphisms give smooth morphisms, étale morphisms give étale morphisms, unramified morphisms give unramified morphisms, monomorphisms give monomorphisms, locally finitely presented morphisms give locally finitely presented morphisms (cf. Fact \ref{xyxyxyxyx}, Corollary \ref{corosmoothxy}, Corollary \ref{coroétt}, Corollary \ref{corounramifiedd}, Fact \ref{monomono}, Lemma \ref{equiclosed}, Proposition \ref{opclosmoo}, Fact \ref{xypres}). However, flat morphisms do not give flat morphisms on attractors in general (cf. Remark \ref{flatvistoli}).
\item Attractors of attractors make sense and correspond to intersections of magnets (cf. Remarks \ref{actiongroup} and \ref{actionmono} and Proposition \ref{NLNL}).
\item Attractors associated to subgroups correspond to fixed-points under diagonalizable group schemes (cf. Proposition \ref{Group}).
\item Faces and attractors allow to obtain easily non trivial cartesian diagrams (cf. Remark \ref{cartesian-para} and Proposition \ref{cartesian}). 
\item Attractors preserve smoothness in two different senses (cf. Section \ref{fosmo}, e.g. Corollaries \ref{corosmoothxy} and \ref{123456789}). 
\item The morphism $X^N \to X^{N^*}$ induces a bijection on the sets of connected components (cf. Proposition \ref{proptopo}). Moreover if $X$ is smooth, then $X^N \to X^{N^*}$ is an affine space fibration (cf. Theorem \ref{BB}).
\item Attractors commute with tangent spaces and Lie algebras (cf. Propositions \ref{tangent} and \ref{lieN}).
\item Attractors make sense for ind-spaces and behave as one can expect (cf. Proposition \ref{indprop}).
\item Attractors commute with dilatations (cf. Proposition \ref{dilattra}).
\end{enumerate}

\subsection{Examples}
Let us discuss some examples. 
We refer to §\ref{examplesse} and Example \ref{exammh} for some other examples and applications.
\exam
Assume that $X = \mathbb{A}^n_S=: \mathbb{V}$ and $A(M)_S$ acts linearly so that the action induces a direct sum decomposition in weight spaces $\mathbb{V} = \bigoplus _{m \in M } \mathbb{V} _m $, then $\mathbb{V}^N = \bigoplus _{n \in N } \mathbb{V} _n  \subset \mathbb{V}$.
\xexam

\exam \label{exintro1} (Magnetic point of view on reductive groups) Let $G$ be a reductive group scheme over $S$. Let $T = D (M)_S$ be a maximal split torus of $G$. Let $a$ be the adjoint action of $T$ on $G$. The set of pure magnets is given by additively stable sets of roots. Algebraic attractors associated to the action $a$ give all the well-known classical objects of the theory of reductive groups (Levi subgroups containing $T$, parabolic subgroups containing $T$, ``groupes de type R à fibres résolubles"  \cite[Exp. XXII §5.6]{SGA3}). Moreover attractors with prescribed limits give unipotent radicals of such objects, namely root groups and unipotent radicals of parabolic groups. We refer to §\ref{sectroots}, Proposition \ref{conradrootSGA} and \cite[§6.3]{ALRR22} for some precise statements. 
\xexam

\subsection{Related literature}  \label{relation}
Let us list the main sources of inspiration for our work.\begin{enumerate}
\item Our work is of a completely different nature than logarithmic algebraic geometry.  However it has in common with logarithmic algebraic geometry to use schemes associated to monoids as background and \cite[Part I]{Og} was useful at some stages of the realization of our work (cf. e.g. Proposition \ref{flatflatflat} and Fact \ref{fifififif}).
 \item 
If $M= \bbZ$ and $N = \bbN$, then our attractors $X^N$ equal the well-known attractors $X^+$ associated to $\bbG _m$-actions. As a consequence, our work was partly inspired by the following works, which include studies of $\bbG_m$-actions \[[\text{Ref}_{\bbG _m}]:= \Big\{\cite{Bi73}, \cite{He80}, \cite{Ju85}, \cite{CGP10},\cite{Dr15}, \cite{DG14}, \cite{Mar15}, \cite{Ri16}, \cite{Mi17}, \cite{HR21}\Big\}.\] We now provide some examples of relations between statements in our work and in $[\text{Ref}_{\bbG _m}]$: \begin{enumerate} \item statement of Definition \ref{definitio} was partly inspired by \cite[Definition 1.3]{Ri16}, \cite[Definition 1.3.2]{Dr15}, \cite[II.4.1]{He80}, etc, 
\item statement of Theorem \ref{representableaffine} was partly inspired by \cite[Lemma 2.1.4]{CGP10}, \cite[§1.3.4]{Dr15}, \cite[Lemma 1.9]{Ri16},  etc,
\item statement of Proposition \ref{fixxx} was partly inspired by \cite[Lemma 1.11]{Ri16},
\item statements of Propositions \ref{lemmpf} and \ref{repgroup}  and Theorem \ref{G_connected_implies_fixed_pts_closed} were partly inspired by \cite[Lemma 1.10]{Ri16}, \cite[Theorem 1.8 (i)]{Ri16}, \cite[Proposition 1.2.2.]{Dr15}, etc,
\item statement of Theorem \ref{representable} was partly inspired by \cite[Theorem 1.8]{Ri16}, \cite[Theorem 1.4.2]{Dr15}, etc
\item statement of Proposition \ref{proptopo} was partly inspired by \cite[Corollary 1.12]{Ri16},
\item statement of Corollary \ref{corosmoothxy} was partly inspired by \cite[Theorem 1.1]{Mar15}, \cite[Theorem 1.8 (iii)]{Ri16}, \cite[Lemma 2.2 (ii)]{HR21}, 
\item statement of Corollary \ref{123456789} was partly inspired by \cite[Proposition 1.4.20]{Dr15},
\item statement of Proposition \ref{opclosmoo} was partly inspired by \cite[Corollary 2.3]{HR21},
\item statement of Theorem \ref{BB} was partly inspired by \cite{Bi73} and \cite[Theorem 13.47]{Mi17}.
\end{enumerate}

Similarly, some of our proofs are also partly inspired by the proofs in $[\text{Ref}_{\bbG _m}]$. Our proof of smoothness results (Corollaries \ref{corosmoothxy} and \ref{123456789}) use formal smoothness (Propositions \ref{formsmooth} and \ref{formsmoothface}) and was partly inspired by \cite[Exp. XII Théorème 9.7 (unpublished)]{SGA3}.
Other statements on attractors associated with monoids in our paper do not have analogues in $[\text{Ref}_{\bbG _m}]$ cf. e.g. Fact \ref{faceattractors}, Proposition \ref{inter}, Proposition \ref{cartesian}, Proposition \ref{dilattra}, Corollary \ref{dilamagnecor}, Theorem \ref{efficienttheo}, Theorem \ref{strumagnet}, etc. We note that Proposition \ref{NLNL} makes sense for $\bbG _m$-actions (it says that $(X^+)^-=(X^-)^+=X^0$) but we did not come across it in $[\text{Ref}_{\bbG _m}].$ Even in the case where $M=\bbZ$ (i.e. $D(M)_S = \bbG_{m,S}$), our work does not specify exactly to the formalism of $\bbG _m$-attractors as in [Ref$_{\bbG _m}$] as we require a choice of magnet before considering an attractor. Our work refines strictly the formalisms in [Ref$_{\bbG _m}$], even for $\mathbb{G}_m$. Again, the concept of magnets and the use of monoids was absent in [Ref$_{\bbG _m}$].
Let us inform that \cite{JS18} and \cite{JS20} generalize attractors associated to $\bbG _m$-actions for general reductive groups over fields. In the special case where $M= \bbZ^r$ (and $S$ is the spectrum of a field as in \cite{JS18} and \cite{JS20}), all the $N \subset M$ such that $A(N)_S$ fulfill the conditions of \cite{JS18, JS20} have to satisfy several conditions. Indeed, we quote from \cite{JS18} where $\overline{G}$ plays the role of our $A(N)_S$ (and $\bbA^1$ in [$\text{Ref}_{\mathbb{G}_m}$]):
{``The functorial description readily generalizes
once we understand what should be put in place of $\bbA^1$. It turns out that a suitable replacement
is a linearly reductive monoid $\overline{G}$, i.e., an affine variety with multiplication $\overline{G} \times \overline{G} \to \overline{G}$ and a unit,
such that $G \subset \overline{G}$ is the dense submonoid consisting of invertible elements."} Note that \cite{JS18} and \cite{JS20} do not work in the framework of abelian monoids or actions of diagonalizable monoid schemes, whereas our theory is based crucially on their rigid structure. In our theory, magnets, i.e. submonoids of $M$, are the main objects of study, and all submonoids of $M$ are allowed, which is essential for many of our results and is possible only due to the rigidity of diagonalizable monoid actions. In Algebraic Magnetism, the morphism $A(M)_S \to A(N)_S$ is not required to be an open immersion (compare with $ \mathbb{G}_m \subset \mathbb{A}^1$ in [$\text{Ref}_{\mathbb{G}_m}$] and $G \subset \overline{G}$ in \cite{JS18, JS20}), $M$ is allowed to have torsion, and all submonoids $N \subset M$ are considered (we repeat this point for emphasis).
 
\item  This work is written in the language of Grothendieck schemes \cite{EGA} and Artin's algebraic spaces \cite{Ar71}, we mainly use \cite{stacks-project} as treatment for this theory in the present text. We use algebraic spaces instead of schemes for the same reasons than \cite{Dr15} and \cite{Ri16}: some actions of group schemes that we are interested in are not Zariski locally linearizable but are étale locally linearizable by deep results of Alpher-Hall-Rydh \cite{AHR21}. This leads to use étaleness as local notion instead of openess and so to use algebraic spaces instead of schemes.
\item  Our work was inspired by \cite{SGA3} for some technical aspects around group schemes.
 \end{enumerate}

\subsection{Organization of the paper}
\subsubsection{Sections  \ref{sectsemigroup}-\ref{examplesse}.} 
Section \ref{sectsemigroup} focuses on preliminaries on diagonalizable monoid schemes. Section \ref{sectdefinition} introduces algebraic attractors associated to monoids and prove several results. Section \ref{prescribed} introduces attractors with prescribed limits. Sections  \ref{basechange}-\ref{repandprop} take care of several results used to prove the representability of attractors in non-affine cases, in particular the notion of $Z$-FPR and strongly-FPR atlases are introduced. Sections \ref{hoschsection}-\ref{fosmo} deal with formal smoothness and formal étaleness results.   Section \ref{topo} is about topology and geometry of attractors, in particular we generalize the Bialynicki-Birula decomposition. Section \ref{dilatationse} shows that attractors are compatible with dilatations. Section \ref{indspaces} takes care of ind-spaces. Section \ref{sectionefficient} studies pure magnets. Section \ref{examplesse} provides some examples and applications to the structure of algebraic groups, in particular we study some general results about actions of diagonalizable group schemes on group schemes.

\subsubsection{Appendix \ref{appendixx}.} Appendix \ref{appendixx}, written by M. Romagny, is devoted to prove the existence of $Z$-FPR atlases (cf. Theorem \ref{adapted-atlas}) using deep results of Alper-Hall-Rydh \cite{AHR21}.  Some results of the paper are stated under the existence of $Z$-FPR atlases, cf. e.g Theorem \ref{representable}. One knows that such atlases exist for Sumihiro's actions, but some actions are not Sumihiro (cf. Section \ref{Zariskicase}) and in this case one needs to know the existence of $Z$-FPR atlases using other tools. This is why Appendix \ref{appendixx} is useful in this paper.
Appendix \ref{appendixx} also contains generalizations of some results stated in other sections (cf. Proposition \ref{lemmpf}, Proposition \ref{repgroup}, and Theorem \ref{G_connected_implies_fixed_pts_closed}).

\subsection{Guide to the reader}

We recommend that the inexperienced reader start with Section~\ref{sectdefinition} and refer to Section~\ref{sectsemigroup} when indicated. Then we suggest admitting that attractors are representable, i.e.\ skipping Sections~\ref{basechange}--\ref{FPR}--\ref{Zariskicase}--\ref{repandprop} and Appendix \ref{appendixx}, and directly reading Sections~\ref{topo} (Topology and Geometry) and \ref{sectionefficient} (Pure Magnets). In this way, the reader will have a good vision of the phenomena at play. Reading examples in Section~\ref{examplesse} will give a further idea of the subject.

Sections~\ref{hoschsection}--\ref{deformation}--\ref{fosmo} (smoothness), Section~\ref{dilatationse} (compatibility with dilatations), and Section~\ref{indspaces} (Ind-spaces) are more advanced and can be skipped on a first read.

\newpage 
\section{Rings associated to monoids and their spectra} \label{sectsemigroup}
 
This section focuses on diagonalizable monoid schemes associated with commutative monoids. Subsection \ref{ss:og} recalls some basic terminology and indicates a reference for abelian monoids. Subsection \ref{ss:dms} contains several results on diagonalizable monoid algebraic spaces and diagonalizable monoid schemes. Section \ref{sectsemigroup} will be used throughout the paper, as diagonalizable monoid schemes are at the heart of the theory. Definitions \ref{ringN} and \ref{defiAA} are the core of this section. The section is otherwise organized as a collection of various results, and the reader may refer to it as needed and return to it later.

\subsection{Commutative monoids} \label{ss:og}
We refer to \cite{Og} for a detailed and beautiful introduction to monoids and related structures. We recall in this section some basic definitions and facts that we frequently use in our work. 
In this article, monoids and rings are always commutative.
 Readers unfamiliar with monoids should read \cite[I.1]{Og}. Let us recall some very basic notations, again we refer to \cite[I]{Og} for a more general and conceptual presentation.
Let $P$ be a monoid and $N\subset P$ and $L \subset P$ be two submonoids. Let $N+ L $ be $ \{ n+ l \in P  | n  \in N , l \in L \} $. Then $N+L$ is a submonoid of $P$. An arbitrary intersection of submonoids is a submonoid.
Let $P$ be a monoid and let $E$ be a subset of $P$, then we denote by $[E\rangle$ the smallest submonoid of $P$ containing $E$, this is the intersection of all monoids of $P$ containing $E$.
Similarly if $M$ is a group and $E$ is a subset of $M$, we denote by $( E ) $ the subgroup of $M$ generated by $E$.
 Let $M$ be an abelian group and let $N$ be a submonoid of $M$. The subgroup generated by $N $ in $M$ is denoted $N^{\mathrm{gp}} $, in fact $N^{\mathrm{gp}}= \{x-y |x,y \in N\}$. We have an obvious notion of finitely generated monoids. The monoid $(\bbN _{\geq1} \times \bbN ) \cup (0,0)$ is not finitely generated and it is a submonoid of the finitely generated abelian group $\bbZ \times \bbZ $. 

\subsection{Diagonalizable monoid schemes} \label{ss:dms}

We fix in this section an arbitrary commutative monoid $M$. We will often work with submonoids of $M$ and will often use the symbols $N,N',L,Q$ or $F$ to denote them. 

\defi \label{defiZn}
Let $N$ be a commutative monoid. Let $\bbZ[N]$ be the ring
whose underlying abelian group is $\bigoplus _{n \in N} \bbZ X^n$, where for any $n $ the abelian group $\bbZ X^n$ is a formal copy of $\bbZ$, and multiplication is induced by the operation of the monoid: $X^n \times X^{n'}=X^{n+n'}$. We write $X^0=1$. The ring $\bbZ[N]$ is called the ring associated to the monoid $N$.
\xdefi

\fact \label{bialgebra} Let $N$ be a commutative monoid. The ring $\bbZ[N]$ is a bialgebra over $\bbZ$.
The augmentation is the map $ \bbZ [{N}] \to \bbZ$ sending $X^n$ to $1$ for every $n \in N $.
The comultiplication is the map $\bbZ[N] \to \bbZ[N] \otimes \bbZ [N] $ sending $X^n$ to $X^n \otimes X^n$ for  $n \in N$.
Moreover, if $N$ is a group,  $\bbZ[N]$ is a Hopf algebra over $\bbZ$, the antipode being the map $\bbZ [{N}] \to \bbZ [N]$ sending $X^n$ to $X^{-n}$ for every $n \in N $.

\xfact

\defi \label{ringN} Let $N$ be a commutative monoid. Let $R$ be a ring. Let $B $ be an algebraic space over a scheme $S$. Let $\calO_S$ be the structure sheaf of $S$ (cf. \cite[\href{https://stacks.math.columbia.edu/tag/01IJ}{Tag 01IJ}]{stacks-project} and  \cite[\href{https://stacks.math.columbia.edu/tag/0091}{Tag 0091}]{stacks-project}). Let $\calO _B$ be the structure sheaf of $B$ in the sense of \cite[\href{https://stacks.math.columbia.edu/tag/04KD}{Tag 04KD}]{stacks-project}.
 \begin{enumerate} \item We put $R[N]= \bbZ [N] \otimes _{\bbZ} R$. Obviously, $R[N] = \bigoplus _{n \in N } R X^n$. The ring $R[N]$ is canonically a $R$-bialgebra.  
 \item Let $\calO _S[N]$ be the $\calO _S$-algebra obtained by sheafification of the presheaf of algebras given by $ \calO _S (U) [N]$  for any open subset $U \subset S$, this is a sheaf of algebras (cf. \cite[\href{https://stacks.math.columbia.edu/tag/00YR}{Tag 00YR}]{stacks-project}).  Note that the underlying $\calO _S$-module is called the free $\calO_S$-module with basis $N$. In particular, as $\calO _S$-module,  $\calO _S [N]$ is isomorphic to $\bigoplus _{n \in N} \calO _S$.  If $U$ is a quasi-compact open subset of $S$, then $\calO _S[N] (U) = \calO _S (U) [N]$ (cf. \cite[\href{https://stacks.math.columbia.edu/tag/01AI}{Tag 01AI}]{stacks-project}). The $\calO _S$-algebra $\calO _S[N]$ is canonically an $\calO _S$-bialgebra.
 
 \item Let $\calO _B[N]$ be the $\calO _B$-algebra obtained by sheafification of the presheaf of algebras given by $ \calO _B (T) [N]$  for any  $T \in B_{spaces, \acute{e}tale} $ (cf. \cite[\href{https://stacks.math.columbia.edu/tag/03G0}{Tag 03G0}]{stacks-project} for the notation $B_{spaces, \acute{e}tale}$), this is a sheaf of algebras (cf. \cite[\href{https://stacks.math.columbia.edu/tag/00YR}{Tag 00YR}]{stacks-project}).  Note that the underlying $\calO _B$-module is called the free $\calO_B$-module with basis $N$ (cf. \cite[\href{https://stacks.math.columbia.edu/tag/03DD}{Tag 03DD}]{stacks-project}).  In particular, as $\calO _B$-module,  $\calO _B [N]$ is isomorphic to $\bigoplus _{n \in N} \calO _B$. If $T \in B_{spaces , \acute{e}tale}$ is quasi-compact, then $\calO _B [N] (T) = \calO _B (T) [N]$ by  \cite[\href{https://stacks.math.columbia.edu/tag/0935}{Tag 0935}]{stacks-project}. The $\calO _B$-algebra $\calO _B[N]$ is canonically an $\calO _B$-bialgebra.

\item  Assume that moreover $N$ is a group, then $R[N]$ (resp. $ \calO _S [N]$, resp. $\calO _B[N]$) is canonically a Hopf algebra over $R$ (resp. $ \calO_S $,  resp. $\calO _B $).
\end{enumerate}
\xdefi

\defi \label{defiAA} (Diagonalizable monoid schemes and algebraic spaces) Let $N$ be a commutative monoid. Let $S$ be a scheme and let $B$ be an algebraic space over $S$. \begin{enumerate} \item Let $A(N)$  be the scheme $\Spec (\bbZ [N])$. This is a monoid scheme over $\Spec (\bbZ )$.

\item We put $A(N)_S= A(N) \times _{\Spec (\bbZ)} S$, this is canonically a monoid scheme over $S$.
\item  We put $A(N)_B=A(N)_S\times _S B$, this is canonically a monoid algebraic space over $B$.
\item The objects $A(N), A(N)_S$ and $A(N)_B$ are called diagonalizable monoid schemes and algebraic spaces. If $N=M$ is an abelian group, then $A(M)$ is denoted $D(M)$ and is a group scheme over $\Spec (\bbZ)$ called the diagonalizable group scheme associated to $M$. We define similarly $D(M)_S$ and $D(M)_B$.\end{enumerate}
\xdefi

\rema Let $M $ be an abelian group and $N$ be a monoid. The reference \cite{SGA3} uses the notation $D_S (M)$ where we use $D(M)_S$. Similarly, the notation $A_B (N)$ can be used to denote $A(N)_B$. We like to use standard notation and see $A(N)_B$ and $D(M)_B$ as base change from $\Spec (\bbZ)$ to $B$ of $A(N)$ and $D(M)$. 
\xrema

\fact \label{fact-RNAN}
Let $N$ be a commutative monoid. Let $B$ be an algebraic space over a scheme $S$. \begin{enumerate}
 \item The monoid scheme $A(N)_S$ is affine and flat over $S$, moreover its quasi-coherent $\calO _S$-bialgebra is $\calO_S [N]$ (cf. Definition \ref{ringN}).
 \item The monoid algebraic space $A(N)_B$ is affine and flat over $B$, moreover its quasi-coherent $\calO _B$-bialgebra is $\calO _B [N]$ (cf. Definition \ref{ringN}).
 \end{enumerate} 
\xfact 
\pf
Being affine and flat is stable by base change, so it is enough to prove that $A(N)$ is affine and flat over $\Spec (\bbZ)$. It is obvious that $A(N)$ is affine and flat over $\Spec (\bbZ)$ and this proves the first parts of both assertions. 
 We claim that the  $\calO _{\Spec (\bbZ)}$-bialgebra of $A(N)$ is $\calO_{\Spec (\bbZ)} [N]$. Let $p = A(N) \to \bbZ$. We know that  $p_*\calO _{A(N)} = \widetilde{\bbZ [N]}$.
 So by \cite[\href{https://stacks.math.columbia.edu/tag/01ID}{Tag 01ID}]{stacks-project}, we have $p_*\calO _{A(N)} = \widetilde{\bigoplus _{n \in N}  \bbZ } = \bigoplus _{n \in N } \widetilde{\bbZ} = \bigoplus _{n \in N } \calO _{\Spec (\bbZ)}= \calO _{\Spec (\bbZ)} [N]$ and this finishes to prove the claim.
  Now by \cite[\href{https://stacks.math.columbia.edu/tag/01SA}{Tag 01SA}]{stacks-project}, we have that $q_*\calO _{A(N)_S} = r^* (p_* \calO _{A(N)})$ where $q: A(N)_S \to S$ and $r: S\to \bbZ$. So by  \cite[\href{https://stacks.math.columbia.edu/tag/01AJ}{Tag 01AJ}]{stacks-project}, we have $q_* \calO _{A(N)_S} = \calO _S [N]$. The last assertion is proved similarly using  \cite[\href{https://stacks.math.columbia.edu/tag/03M1}{Tag 03M1}]{stacks-project}, \cite[\href{https://stacks.math.columbia.edu/tag/081V}{Tag 081V}]{stacks-project} and \cite[\href{https://stacks.math.columbia.edu/tag/03DC}{Tag 03DC}]{stacks-project}. 
\xpf

\fact \label{NNNN} We have a canonical identification $A(N \times L) = A(N) \times _{\Spec (\bbZ) }A(L)$ for any pair of monoids $N,L$. If $B$ is an algebraic space over a scheme $S$, then $A(N \times L)_S = A(N)_S \times _{S }A(L)_S$ and $A(N \times L)_B = A(N)_B \times _{B }A(L)_B$.
\xfact

\pf
It is enough to prove the first assertion. This follows from the identity 
\[ \bbZ[N \times L] = \bigoplus _{(n,l) \in N \times L} \bbZ X^{(n,l)} = \big( \bigoplus _{n \in N } \bbZ X^n \big) \otimes _{\bbZ}  \big( \bigoplus _{l \in L } \bbZ X^l \big) = \bbZ [N] \otimes _{\bbZ} \bbZ [L].\]
\xpf 

Recall that a monoid $N$ is cancellative if for all $x,y,z \in N, $ $x+y=x+z$ implies $y=z$. In \cite{Og}, the word ``integral" is used instead of ``cancellative", cf. \cite[Definition 1.3.1]{Og}. A monoid is cancellative if and only if it identifies with a submonoid of a group.

\prop  \label{flatflatflat} Let $N$ be a cancellative commutative monoid. Assume that the action of $N$ on $N \times N$ given by $n \cdot (a,b)=(a+n,b+n)$ is free (in the sense of \cite[§1.2]{Og}). Let $B$ be an algebraic space over a scheme $S$.\begin{enumerate} \item  The morphism of rings $\bbZ[N] \to \bbZ[N \times N ] $, $X^n \mapsto X^{(n,n)}$ is flat. 
\item The multiplication $m: A(N) \times _{\Spec (\bbZ)} A(N) \to A(N)$ of the monoid scheme $A(N)$ is flat.
\item The multiplication morphism $m:A(N)_S \times _{S} A(N)_S \to A(N)_S$ is flat.
\item  The multiplication morphism $m:A(N)_B \times _{B} A(N)_B \to A(N)_B$ is flat. \end{enumerate}
\xprop

\pf The monoid $N$ acts freely on $N \times N$ via $n \cdot (m,l)=(m+n,l+n)$.
So the $\bbZ[N]$-module $\bbZ[N \times N]$ is flat by \cite[Prop. 4.5.12 p.134]{Og}, and so (i) holds. Now (ii) follows using Facts \ref{bialgebra} and \ref{NNNN}. Assertions (iii) and (iv) follow since flatness is preserved by base change.
\xpf

\rema
If $M$ is a commutative group, the action of $M$ on $M\times M $ given by $n \cdot (a,b)=(a+n,b+n)$ is free. Let $N $ be $\bbN \setminus \{1\}$. Then the action of $N$ on $N\times N $ given by $n \cdot (a,b)=(a+n,b+n)$ is not free because the orbits of $(2,3)$ and $(3,4)$ have non-trivial intersection (e.g. the intersection contains $(4,5)$) but are different (e.g. $(3,4)$ is not in the orbit of $(2,3)$).
\xrema

\fact \label{fifififif} Let $N$ be a finitely generated monoid. Let $B$ be an algebraic space over a scheme $S$. \begin{enumerate}
\item The monoid $N$ is finitely presented.
\item The $\bbZ$-algebra $\bbZ [N]$ is finitely presented.
\item The structural morphism $A(N) \to \Spec (\bbZ )$ is finitely presented.
\item The multiplication morphism $A(N) \times _{\Spec (\bbZ)} A(N) \to A(N)$ is finitely presented.
\item The morphisms $A(N)_S \to S $ and $A(N)_S \times _S A(N)_S \to A(N)_S$ are finitely presented.
\item The morphisms $A(N)_B \to B $ and $A(N)_B \times _B A(N)_B \to A(N)_B$ are finitely presented. \end{enumerate}
\xfact 

\pf \begin{enumerate}
\item This is \cite[Theorem 2.1.7]{Og}.
\item This follows from (i).
\item This follows from (ii).
\item Since $N$ is finitely generated, $N \times N $ is finitely generated and so $A(N) \times _{\Spec (\bbZ)} A(N)$ is finitely presented over $\Spec (\bbZ)$. Now the assertion follows using  \cite[\href{https://stacks.math.columbia.edu/tag/02FV}{Tag 02FV}]{stacks-project}.
\item This follows using previous assertions and \cite[\href{https://stacks.math.columbia.edu/tag/01TS}{Tag 01TS}]{stacks-project}.
\item This follows using previous assertions and \cite[\href{https://stacks.math.columbia.edu/tag/049M}{Tag 049M}]{stacks-project}.
\end{enumerate}
\xpf

\defi \begin{enumerate}${}$ \item 
Let $e_\bbZ$ be the trivial monoid scheme over $\Spec (\bbZ)$, it is a group scheme over $\Spec (\bbZ)$ equal to $\Spec (\bbZ)$ as scheme. Note that $e_\bbZ \simeq  A(0)$ where $0$ is the trivial monoid. 
\item Let $S $ be a scheme. We put $e_S = e_\bbZ \times _{\Spec (\bbZ )} S$, this is the trivial monoid scheme over $S$.
\item Let $B$ be an algebraic space over a scheme $S$. We put $e_B = e_S \times _S B $, this is the trivial monoid algebraic space over $B$. \end{enumerate}
\xdefi

\fact \label{actionelemen}
Let  $N$ be a submonoid of $ M$. Let $B$ be an algebraic space over $S$. Then we have an algebraic action of $A(M)$ on $A(N)$ over $\Spec (\mathbb{Z} )$ given by:
\begin{align*}
 \bbZ [{N}] &\xrightarrow{} \bbZ [{N}] \otimes \bbZ [{M}]\\
 X^n & \mapsto X^n \otimes X^n  .
\end{align*}
Equivalently, this action comes from the general Proposition \ref{4.7.3.1mono} that we discuss below in this section.
Remark also that the action of $A(M)$ on $A(N)$ comes from the morphism of monoid schemes $A(M)\to A(N)$. 
By base change, we obtain actions of $A(M)_S $ on $A(N)_S$ and of $A(M)_B$ on $A(N)_B$. 
\xfact 

\fact \label{incluA}
Let  $N' \subset N$ be submonoids of $M$, let $B$ be an algebraic space over $S$. Then we have a canonical morphism of bialgebras over $\bbZ$, $\bbZ[N'] \to \bbZ [N]$ (sending $X^{n'} $ to $ X^{n'}$ for any $n' \in N'$). This induces a canonical morphism $A(N) \to A(N')$ of monoid schemes over $\bbZ$, this morphism is $A(M)$-equivariant. We obtain canonical equivariant morphisms of monoid objects $A(N)_S \to A(N')_{S}$ and $A(N)_B \to A(N')_{B}$. 
\xfact
\pf
 The inclusion morphism preserves $M$-gradings on $\bbZ[N']$ and $\bbZ[N]$.
 \xpf

 \defi Let $F \subset N $ be submonoids of $M$.
 We say that $F $ is a face of $N$ if the projection map
 \[\bbZ [N] \to \bbZ [F], X^n \mapsto 0\text{ if } n \in {N} \setminus{{F}} \text{ and } X^n \mapsto X^n \text{ if } n \in {F}\] is a morphism of rings. 
 \xdefi
 
 \fact \label{faceANmor}
 If $F$ is a face of a submonoid $N$ of $M$, then the associated morphism of schemes $A(F) \to A(N)$ is $A(M)$-equivariant. For any algebraic space $B$ over any scheme $S$,  $A(F)_S \to A(N)_S$ is $A(M)_S$-equivariant and $A(F)_B \to A(N)_B$ is $A(M)_B$-equivariant.
 \xfact
 
 \pf
 The projection morphism preserves $M$-gradings on $\bbZ[F]$ and $\bbZ[N]$.
 \xpf
 
 \prop \label{nicepropprop} Let $F \subset N$ be submonoids of $M$. Then $F $ is a face of $N$ if and only if for all $x,y \in N$ \[ x+y \in F \Leftrightarrow x \in F \text{ and } y \in F.\]  
 \xprop
 
 \pf Let $\phi$ denote the projection and assume it is a morphism of rings. Let $x,y \in N$. Then  $x+y \in F$ $ \Leftrightarrow$ $\phi (X^{x+y})=X^{x+y}= \phi(X^x)\phi(X^y)$ is not zero $ \Leftrightarrow $ both $x$ and $y$ are in $F$.
Reciprocally assume that for all $x,y \in N$, $x+y \in F\Leftrightarrow x \in {F}$ and $ y \in {F}$, then we have $\phi (X^{x}X^y) = \phi (X^x ) \phi (X^y)$.
 \xpf
 
 \prop \label{abelianimplynice}
  Let $N $ be a monoid. Let $N^{*} =\{ x \in N | ~\exists y \in N \text{ such that } x+y =0 \} $, then $N^{*}$ is a submonoid of $N$ and a group, moreover $N^{*} $ is a face of $ N$. The group $N^*$ is the largest subgroup of $N$, called the face of invertible elements.
 \xprop
 
 \pf  We apply Proposition \ref{nicepropprop} as follows. Take $x,y \in N$. Assume $x+y \in N^{*}$, then there exists $z \in N$ such that  $x+y+z =0$, this shows that $x$ and $y$ are in $N ^{*}$.
 \xpf

 \fact Let $F$ be a face of a monoid $N$, then $N^* \subset F$.
 \xfact
\pf
Let $n \in N^*$. Then $n+(-n)=0 $ belongs to $F$ and so $n$ belongs to $F$.
\xpf

  \prop \label{lemmpush} Let $ N \subset L \subset L' $ be submonoids of $M$. Assume that $L$ is a face of $L'$.  Then $N' := L' \setminus (L \setminus N )$ is a submonoid of $M$  and $N$ is a face of $N'$; moreover for any scheme $S$, $A(N') _ S$ is the push-out, in the category of schemes, of  \[ \begin{tikzcd}
   & A(L)_ S  \ar[rd] \ar[dl]& \\
   A(L')_S & & A(N)_S
  \end{tikzcd}.\] 
  \xprop
  \pf We have $L' \setminus (L \setminus N )= (L' \setminus L) \sqcup N$ and it is clearly a submonoid of $M$. It is also clear that $N$ is a face of $N'$. Let us assume firstly that $S= \Spec (R)$ is affine.
      We have $R[L'] \times _{R[L]} R[N] \cong \{ (x,y) \in R[L'] \times R[N] | f(x)=g(y) \} \cong R[N']  $, indeed an element $x$ in $R[L']$ maps to an element in $R[N]$ under the projection $R[L']\to R[L]$ if and only if $x \in R[(L' \setminus L )\sqcup N]$.
      The map $R[N'] \to R[N]$  is the projection morphism associated to the face inclusion $N \subset N'$. The map $R[N'] \to R[L']$ is the morphism associated to the inclusion $N' \subset L'$.
      So by \cite[0ET0]{stacks-project} the scheme $A(N')_R$ is the push-out, in the category of schemes, of the diagram 
 \[ \begin{tikzcd}
   & A(L)_R \ar[rd] \ar[dl]& \\
   A(L')_R & & A(N)_R
  \end{tikzcd}.\] 
  Now let us prove the general case.
  Let $S = \cup _{i \in I} U_i$ be an affine open covering and write $U_i=\Spec (R_i)$. 
  Let $Y$ be a scheme and let $A(L')_S \to Y$ and $A(N)_S \to Y$ be two morphisms such that the following diagram commutes  \[ \begin{tikzcd}
   & A(L)_ S  \ar[rd] \ar[dl]& \\
   A(L')_S \ar[dr] & & A(N)_S \ar[dl]\\ 
   & Y & 
  \end{tikzcd}.\] 
  We then obtain, for any $i \in I$, a commutative diagram 
  \[ \begin{tikzcd}
   & A(L)_ {U_i}  \ar[rd] \ar[dl]& \\
   A(L')_{U_i} \ar[dr] & & A(N)_{U_i} \ar[dl]\\ 
   & Y & 
  \end{tikzcd}.\] 
  Now since $U_i$ is affine, we obtain a unique morphism $f_i : A(N')_{U_i} \to Y$ such that the following diagram commutes
   \[ \begin{tikzcd}
   & A(L)_ {U_i}  \ar[rd] \ar[dl]& \\
   A(L')_{U_i} \ar[ddr] \ar[dr] & & A(N)_{U_i} \ar[ddl] \ar[dl]\\ 
   &A(N')_{U_i} \ar[d,"f_i"] &\\
   & Y & 
  \end{tikzcd}.\]
 For $i ,j \in I$, we have $U_i \times _S U_j = U_i \cap U_j $. Let $U_i \cap U_j = \cup _{q \in Q} V_q$ be an affine open covering. We have  $f_i |_{ A(N')_{V_q} }= f_j  | _{A(N')_{V_q}}$ for all $ q \in Q$ by the affine case done before. So we have $f_i |_{ A(N')_{U_i \cap U_j} }= f_j  | _{A(N')_{U_i \cap U_j}}$ by \cite[Prop. 3.5]{GW}. Thus using \cite[Prop. 3.5]{GW} again, we obtain a unique morphism $f: A(N')_S \to Y $ such that the following diagram commutes
 \[ \begin{tikzcd}
   & A(L)_ {S}  \ar[rd] \ar[dl]& \\
   A(L')_{S} \ar[ddr] \ar[dr] & & A(N)_{S} \ar[ddl] \ar[dl]\\ 
   &A(N')_{S} \ar[d,"f"] &\\
   & Y & 
  \end{tikzcd}.\] This finishes the proof.
  \xpf 
The category of so-called $G$-$\calO_S$-module introduced in \cite[§4.7]{SGA3} for a group scheme $G$ over $S$ naturally extends to monoid schemes over $S.$ In particular $A(M)_S$-$\calO_S$-modules are well-defined. For the convenience of the reader we recall the definition here. Recall that if $\calF$ is an $\calO_S$-module, then $W(\calF)$ is defined as a functor on $Sch/S$ by the formula $W(\calF ) (T) = ( p^* \calF) (T)$ where $p:T \to S$ is the associated morphism. Then $W(\calF) $ is an $O_S$-module where $O_S$ is the ring functor on $Sch/S$ given by the formula $O_S (T)= \calO _T (T)$.

\defi \label{GOmodZarsch}An $A(M)_S$-$\calO_S$-module is an $\calO_S$-module $\calF$ such that \begin{enumerate} \item  for any $S$-scheme $T$, the monoid $A(M)_S (T) $ acts on $W(\calF)(T)$  and $h \cdot ( x + \lambda y ) = h \cdot x  + \lambda (h \cdot y )$ for all $h \in A(M)_S(T)$, $x,y \in W(\calF ) (T)$ and $\lambda \in O_S (T)$, \item  for any $S$-morphism $T \to T'$ the actions of  $A(M) _S(T) $ on $W(\calF)(T)$ and of  $A(M) _S(T') $ on $W(\calF)(T')$ are compatible with the morphisms $A(M)_S(T') \to A(M)_S(T)$ and $W(\calF) (T') \to W(\calF) (T)$.
\end{enumerate}
In other words, it is an $\calO_S$-module $\calF$ with an $O_S$-linear action of the $S$-monoid $A(M)_S$ on $W( \calF)$.
\xdefi

\prop \label{4.7.3mono} Let $S$ be a scheme.  The category of quasi-coherent $A(M)_S$-$\calO_S$-modules is equivalent to the category of quasi-coherent $M$-graded $\calO _S$-modules. 
\xprop

\pf
This is proved as \cite[Exp. I Proposition 4.7.3]{SGA3}.
\xpf 


\prop \label{4.7.3.1mono} Let $S$ be a scheme and let $X$ be a scheme over $S$ such that $X \to S$ is affine. Actions of $A(M)_S$ on $X$ correspond to $M$-gradings of the quasi-coherent $\calO _S$-algebra of $X$.
\xprop

\pf
This is proved as \cite[Exp. I Corollaire 4.7.3.1]{SGA3}.
\xpf

\section{Algebraic attractors associated to magnets and properties} \label{sectdefinition}

In this section we introduce attractors associated to monoids and state many properties. We recommend a complete reading of this section, since it informs a lot about the behavior of attractors.
We also recommend a linear reading, because of the logical structure.

Let us fix a base scheme $S$. In this paper, an $S$-space is an algebraic space over $S$.
Let us fix an $S$-space $X$. Let us fix a commutative monoid $M$ and let $A(M)_S$ be the associated diagonalizable monoid scheme over $S$ (cf. Definition \ref{defiAA}).
We assume that $ A(M)_S$ acts on $X$, this means that we have a morphism of $S$-spaces \[A(M)_S  \times _S   X \xrightarrow{action} X\] such that the following diagrams of $S$-morphisms commute 
\[
\begin{tikzcd} A(M)_S \times _S A(M)_S \times _S X  \ar[d, "\Id \times action"] \ar[rr, "{m \times \Id} "]  & & A(M)_S \times _S X \ar[d, "{action}"] \\
 A(M)_S \times _S X \ar[rr,"action"]& & X \end{tikzcd}
 \]
 
 \[ \begin{tikzcd} 
  e_S \times _S X  \ar[d,"{\Id \times \Id} " ] \ar[rr,"\varepsilon \times action" ]  & & A(M)_S \times _S X \ar[d, "{action}"] \\
 S \times _S X \ar[rr,"\Id"]& & X \end{tikzcd}
 \]where $\varepsilon$ is the unit, i.e. the unique morphism of $S$-monoids from $e_S$ to $ A(M)_S$. In this section, a magnet $N$ is by definition a submonoid $N \subset M$, this terminology is specific to the purpose of algebraic magnetism.
For any magnet $N\subset M$, we consider $A(N)_S$ with the canonical action of $A(M)_S$ as in Fact \ref{actionelemen}.
For any scheme $T \to S$,   $\Hom _T ^{A(M)_T}(A(N)_T, X _T)$ denotes the set of $A(M)_T$-equivariant $T$-morphisms from $A(N)_T$ to $ X_T = X \times _S T$. 
We now introduce the attractor $X^N$ associated to a magnet $N \subset M$ under the action of $A(M)_S$ on $X$.

\defi \label{definitio} Let $X^N:(Sch/S) \to Set$ be the contravariant functor that associates to an object $(T \to S)$ the set \[ \Hom _T ^{A(M)_T}(A(N)_T, X _T)\]
 and that associates to the $S$-morphism $T' \to T$ the map \begin{align*}
 \Hom _T ^{A(M)_T}(A(N)_T, X _T)&\to \Hom _{T'} ^{A(M)_{T'}}(A(N)_{T'}, X _{T'}) \\
 ( f: A(N)_T \to X_T ) & \mapsto (f_{T'} : A(N)_T \times _T T' \xrightarrow{f \times \Id} X_T \times _T T' ).
\end{align*}
\xdefi

\prop \label{fppfsheaftopos} The attractor $X^N$ is an object in the topos $Sh \big( (Sch/S)_{fppf}\big)$, i.e. $X^N$ is a sheaf on the site $(Sch/S)_{fppf}$.
\xprop

\pf \begin{sloppypar}Let $PSh((Sch/S)_{fppf})$ be the category of all contravariant functors from $(Sch/S)_{fppf}$ to $ Set$.  The functor $X^N$ is the equalizer, in the category  ${PSh}((Sch/S)_{fppf})$, of 
\[ 
\begin{tikzcd}  \underline{\Hom} _S ( A(N)_S , X) \ar[r,shift left=.75ex,"\Phi"]
  \ar[r,shift right=.75ex,swap,"\Psi"] & \underline{\Hom} _S ( A(M)_S \times _S A(N)_S , X)  \end{tikzcd}  \] where $\Phi, \Psi $ are defined by: for any $T /S$ and any $f \in \Hom _T ( A(N)_T , X)$,  $ \Phi (f) (h,a)= f(h \cdot a)$ and  $\Psi (f) (h,a) = h \cdot f(a)$ for any $T'/T$ and any $h \in A(M)_T (T') , a \in A(N)_T (T')$.  
The topos $Sh \big( (Sch/S)_{fppf}\big)$ is a full subcategory of ${PSh}((Sch/S)_{fppf})$.
The functor $\underline{\Hom}_S (X' , X) : T/S \mapsto \Hom _{T} ( X'_T , X_T)$ is an fppf sheaf for any pair of $S$-schemes $X,X'$, so it belongs to the topos  $Sh \big( (Sch/S)_{fppf}\big)$.  Equalizers in a category are unique when they exist. Equalizers exist in any topos. The forgetful functor $Sh \big( (Sch/S)_{fppf}\big)\to PSh \big( (Sch/S)_{fppf}\big)$ is a right adjoint (cf. \cite[\href{https://stacks.math.columbia.edu/tag/00WH}{Tag 00WH}]{stacks-project}) so it preserves limits and in particular equalizers.
   We conclude that $X^N$ is equal to the equalizer of $ (\Phi , \Psi)$ in the category $Sh \big( (Sch/S)_{fppf}\big)$, in particular $X^N$ belongs to $Sh \big( (Sch/S)_{fppf}\big)$.
\end{sloppypar}
\xpf

\rema Definition \ref{definitio} makes sense for an arbitrary $S$-functor endowed with an action of $A(M)_S$. We sometimes use this implicitly.
\xrema

\rema \label{Notationavecindicegroup}
In the context of Definition \ref{definitio}, when we want to indicate that we take the attractor $X^N$ relatively to the action of the monoid scheme $A(M)_S$, we sometimes use the notation $X^{N \subset M}$.
\xrema

\rema \label{equal0} We have an identification $X^0 = X^{A(M)_S}$ where $X^{A(M)_S}$ denotes the functor of fixed-points (cf. Proposition \ref{Group} for an extension of this remark in the case where $M$ is a finitely generated group). This is tautological because $A(0)_S $ identifies with $S$ endowed with the trivial action. So the concept of attractors refines the concept of fixed-points. 
\xrema

\fact \label{equalM} We have an identification $X =X^M$ given as follows. Let   $T \to S$ be a scheme. Let $\phi $ be in $ X^M (T)$, we associate the morphism \[  T= e_T \xrightarrow{\varepsilon} A(M)_T  \xrightarrow{\phi} X_T  \in X(T)\] where $\varepsilon $ is the unit. Let $\varphi$ be in $X (T)$, we associate the morphism \[A(M)_T = A(M)_T \times _T T \xrightarrow{ \Id \times \varphi} A(M)_T \times X_T \xrightarrow{action} X_T .\]
\xfact

\prop \label{produitNNN} Let $X_1 , X_2 , X_3$ be $S$-spaces endowed with actions of $A(M)_S$. Let $X_1 \to X_2 $ and $X_3 \to X_2$ be two $A(M)_S$-equivariant morphisms of $S$-spaces. Let $N \subset M$ be a magnet.

\begin{enumerate} \item We have a canonical action of $A(M)_S$ on $X_1 \times _{X_2} X_3$, the canonical maps $X_1 \times _{X_2} X_3 \to X_1 $ and $X_1 \times _{X_2} X_3 \to X_3 $ are $A(M)_S$-equivariant and moreover the following is a cartesian square in the category of $S$-spaces endowed with action of $A(M)_S$ (morphisms in this category are $A(M)_S$-equivariant morphisms of $S$-spaces)
\[
\xymatrix{X_1 \times _{X_2} X_3 \ar[d] \ar[r]& X_1 \ar[d] \\ X_3 \ar[r] & X_2 .}
\]
\item We have an isomorphism of functors \[ (X_1 \times _{X_2} X_3)^N = {X_1}^N \times _{X_2^N} {X_3}^N.\]\end{enumerate}
\xprop

\pf \begin{enumerate} \item Let $T$ be a scheme. Then $A(M)_S(T)$ acts on $X_1(T) \times _{X_2(T)} X_3(T)$ via $g.(x_1,x_3) = (g.x_1,g.x_3)$ (recall that $(X_1 \times _{X_2} X_3) (T) = X_1(T) \times _{X_2(T)} X_3(T) = \{ (x_1, x_3) \in X_1(T ) \times X_3 (T) | x_1\equiv x_3 \text{ in } X_2(T) \}$). This proves the first assertion. The projection maps on $X_1$ and $X_3$ are $A(M)_S$-equivariant by definition. Now let $Y$ be an $S$-space with an action of $A(M)_S$. An $A(M)_S$-equivariant morphism $Y \to X_1 \times _{X_2} X_3$ gives birth by composition with projections to canonical $A(M)_S$-equivariant morphisms $Y \to X_1$ and $Y \to X_2$. Reciprocally let $Y \to X_1$ and $Y \to X_3$ be $A(M)_S$-equivariant. Then we get a morphism in the category of $S$-spaces $Y \to X_1 \times _{X_2} X_3. $ This morphism is $A(M)_S$-equivariant.

\item This is a direct consequence of the previous assertion.
\end{enumerate}
\xpf

\prop \label{fiberprodu} Let $N \subset M$ be a magnet. Let $Y$ and $Z$ be two algebraic spaces over $S$ endowed with actions of $A(M)_S$. Let $A(M)_S$ act componentwise on $Y\times _S Z$.  
We have \[ (Y \times _{S} Z)^N = {Y}^N \times _{S} Z^N.\] 
\xprop
\pf This follows from Proposition \ref{produitNNN}.
\xpf

\prop \label{basechangeee} Let $T \to S$ be an $S$-scheme. The $T$-space $X_T = X \times _S T $ is canonically endowed with an action of $A(M)_T$. Let $N \subset M$ be a magnet. We have an identification of $T$-functors
\[ X^N \times _S T = (X_T )^N . \] 
\xprop
\pf
Indeed, if $T' \to T$ is a $T$-scheme, we have identifications of sets
\begin{align*}
(X^N \times _S T )(T')&= \Hom _T ( T' , X^N \times _S T ) \\ 
&= \Hom _S ( T' , X^N) \\
&= \Hom _{T'}^{A(M)_{T'}} ( A(N)_{T'}, X_{T'})  \text{ and,}\\
 & \\
(X_T)^N (T') &= \Hom _{T'}^{A(M)_{T'}} (A(N)_{T'}, X_T \times _T T')\\
&=\Hom _{T'}^{A(M)_{T'}} ( A(N)_{T'}, X_{T'}) .
\end{align*}
\xpf 
\fact \label{morphismsgeneral} Let $N \subset L$ be magnets of $M$.
For $T \to S$, we have a morphism \[  \Hom _T ^{A(M)_T}(A(N)_T , X_T) \to \Hom _T ^{A(M)_T}(A(L)_T, X_T)\] obtained using the morphism $A(L)_T \to A(N)_T$  (cf. Fact \ref{incluA}). This is functorial, so we get a morphism of functors $\iota _{N,L}	 : X^N \to X^L$. 
\xfact

\rema \label{actiongroup} Let $N \subset M$ be a magnet. We have an action of the monoid scheme $A(M)_S$ on $X^N$ given as follows. For any $S$-scheme $T$, we have actions of $A(M)(T)$ on $X_T$ and $A(N)_T$, in particular for any $t \in A(M)(T)$, we have arrows $X_T \xrightarrow{t \cdot } X_T$ and $A(N)_T \xrightarrow{t \cdot} A(N)_T$. Now let $f \in X^N (T) = \Hom _T ^{A(M)_T} (A(N)_T , X_T) $ and $t \in A(M)(T)$. 
The composition $A(N) _T \xrightarrow{f} X_T \xrightarrow{t \cdot } X_T $ equals the composition $A(N) _T \xrightarrow{t \cdot } A(N)_T \xrightarrow{f} X_T  $ and is denoted $t \cdot f$. This defines an action of $A(M)_S $ on $X^N $.
 \xrema
 
 \rema \label{actionmono} We proceed with the notation from Remark \ref{actiongroup}.  We have an action of the monoid scheme $A(N)_S$ on $X^N$ given as follows. Let $f \in X^N (T) = \Hom _T ^{A(M)_T} (A(N)_T , X_T ) $ and $t \in A(N)(T)$. We define $t \cdot f$ to be  $A(N) _T \xrightarrow{t \cdot } A(N)_T  \xrightarrow{f} X_T $.
 \xrema 
 \rema \label{remaactions}
 The action of Remark \ref{actiongroup} can be obtained from the action of Remark \ref{actionmono} via the canonical morphism of monoid schemes $A(M)_S \to A(N)_S$. This follows from definitions and the fact that the action of $A(M)_S$ on $A(N)_S$ comes from the canonical morphism of monoid schemes $A(M)_S \to A (N)_S$ induced by the inclusion $N \subset M$. Let us mention a sanitary check. If $L \subset M$ is an other magnet, then under some working assumptions $
(X^{N \subset M})^{L \subset M} = (X^{N \subset M})^{(N \cap L  ) \subset N}$ (e.g. use directly \ref{propequiv}, \ref{NLNL}, \ref{representableaffine} or \ref{representable}).
 \xrema
 
\prop Let $N \subset M$ be a magnet. \label{groupmonoidstru} \begin{enumerate}
\item If $X$ is a monoid algebraic space over $S$ and the action of $A(M)_S$ is by monoid endomorphisms, then $X^N$ is a monoid functor.
 \item  If $X$ is a group algebraic space over $S$ and the action of $A(M)_S$ is by group automorphisms, then $X^N$ is a group functor.

\end{enumerate}
\xprop
\pf 
\begin{enumerate} 
\item  Let $T$ be a scheme over $S$. Let $A(N)_T \overset{g}{\to}  X_T$ and $A(N)_T \overset{h}{\to} X_T$ be two elements in $X^N(T)$. Then we define $gh$ as the composition
\[ A(N)_T \to A(N)_T \times _T A(N)_T \to X_T \times _T X_T  \to X_T \] where the first morphism is the diagonal morphism, the second is $g \times h$, and the third is the multiplication morphism coming from the group structure on $X$. The two firsts are equivariant by definitions and the third is equivariant because $A(M)_S$ acts on $X_S$ by monoid endomorphisms. This defines a monoid law on $X^N$. 
\item By (i) we have a monoid law on $X^N$. This law is a group law.
\end{enumerate}
\xpf
 \fact \label{faceattractors} Let $F$ and $N$ be magnets of $M$ and assume that $F $ is a face of $N$. Then for all $T \to S$ the $A(M)_T$-equivariant morphism $A(F)_T \to A(N)_T$ (cf. Fact \ref{faceANmor}) induces a morphism 
 \[ \Hom _T ^{A(M)_T}(A(N)_T, X_T) \to \Hom _T ^{A(M)_T}(A(F)_T , X_T). \]
 So we obtain a morphism of functors $X^N \to X^F$, that we denote $\mathfrak{p}_{N,F}$. The morphism $\mathfrak{p}_{N,F}$ (sometimes also denoted $r$) satisfies $\mathfrak{p}_{N,F} \circ \iota _{F,N} = \Id_{X^F}$.
\xfact
\pf Clear. \xpf

\fact We proceed with the notation from Fact \ref{faceattractors}. Assume that $X^{N} $ and $X^F$ are representable and $X^N \to X^F$ is separated over $S$ (e.g. if $X^N /S$ is separated by \cite[\href{https://stacks.math.columbia.edu/tag/03KR}{Tag 03KR}]{stacks-project}, cf. also Theorem \ref{representable}). Then $\iota _{F,N}: X^F \to X^N$ is a closed immersion.
\xfact
\pf
The composition $X^F \to X^N \to X^F$ being the identity, it is a closed immersion. Now $X^N \to X^F$ is separated and we apply \cite[\href{https://stacks.math.columbia.edu/tag/0AGC}{Tag 0AGC}]{stacks-project}. 
\xpf 

\fact \label{xyxyxyxyx} If $f:X \to Y$ is an $A(M)_S$-equivariant morphism of algebraic spaces, then for any magnet $N$ of $M$, we have a morphism of functors $f^N:X^N \to Y^N$.
\xfact

\pf
For any $T \to S$, send an equivariant arrow $ A(N)_T \to X_T $ to $A(N)_T \to  X_T \to Y_T $.
\xpf

\fact \label{monomono}  Let $N \subset M$ be a magnet. Assume $X \to Y$ is an $A(M)_S$-equivariant monomorphism of $S$-spaces, then $X^N \to Y^N$ is a monomorphism.
\xfact

\pf Let $T \to S$ be a scheme.
Let \begin{tikzcd}  A(N)_T \ar[r,shift left=.75ex,"f"]
  \ar[r,shift right=.75ex,swap,"g"] & X_T \end{tikzcd} be two equivariant $S$-morphisms and assume the compositions \begin{tikzcd}  A(N)_T \ar[r,shift left=.75ex,"f"]
  \ar[r,shift right=.75ex,swap,"g"] & X_T \ar[r]& Y_T \end{tikzcd} are equal. Since $X_T \to Y_T$ is a monomorphism, we get $f=g$.
\xpf 

\fact \label{xypres}  Let $N \subset M$ be a magnet. Assume $f:X \to Y$ is an $A(M)_S$-equivariant morphism of $S$-spaces. If $f$ is locally of finite presentation, then $f^N : X^N \to Y^N$ is locally of finite presentation.
\xfact
\pf
This follows from the definitions and \cite[\href{https://stacks.math.columbia.edu/tag/04AK}{Tag 04AK}]{stacks-project}.
\xpf

\theo \label{representableaffine} Assume that $X$ is affine over $S$. Let $p : X \to S$ be the structural morphism and $\calA = p_* \calO _X$ so that $X = \Spec _S (\calA)$. Let $N$ be a magnet of $M$.
The functor $X^{N}$ is representable by a closed subscheme of $X$ whose quasi-coherent ideal sheaf $\calJ _N$ is the ideal sheaf generated $\{\calA _m |m \in M \setminus N \}$ where $\calA_m$ is the component appearing in the direct sum decomposition $\calA = \bigoplus _{m \in M} \calA _m$ coming from the action of $A(M)_S $ on $X$ (cf. Proposition \ref{4.7.3.1mono}).
\xtheo
\pf Since $X^N$ is a fppf sheaf by Proposition \ref{fppfsheaftopos} and in particular a Zariski sheaf, the statement is local on $S$ using e.g. \cite[\href{https://stacks.math.columbia.edu/tag/01JJ}{Tag 01JJ}]{stacks-project}. We assume $S = \Spec(B)$ and $X = \Spec (A)$, let $A= \bigoplus _{m \in M} A_m$ be the decomposition coming from the action. Let $J$ be the ideal of $A$ generated by $\{ A_m | m \in M \setminus N \}$. Let $B'$ be a $B$-algebra. It is enough to define functorial maps $\Theta $ and $\Psi$ 
\[ \Hom _B ( A/J, B') \overset{\overset{\Theta}{\to}}{\underset{\underset{\Psi}{\leftarrow}}{{\longleftrightarrow }}} \Hom _{B} ^{M\text{-graded}}( A, B' [N]) \]
such that $\Theta \circ \Psi = \Id $ and $\Psi \circ \Theta = \Id$.
Take $A/J \overset{F}{\to} B'$ on the left-hand-side and define a map $f=\Theta (F)$ on the right-hand-side as 
\begin{align*}
 A= &\oplus _{m \in M } A_m  \overset{f}{\to} B'[N] \\
 &a_m \in A_m \mapsto F([a_m]) X^m
\end{align*}
Let us check that this map $\Theta$ is well-defined. Since $F([a_m]) = 0$ if $m\in {M} \setminus N$, the element $F([a_m]) X^m$ belongs to $B'[N] $. We have to explain that $f$ is a morphism of $B$-algebras.
This is a consequence of the identity \[f(a_m a_{m'})= F([a_m a_{m'}]) X ^{m+m'} = F([a_m][a_{m'}]) X^m X^{m'}= f(a_m) f(a_{m'}).\] 
So $\Theta$ is well-defined.
Now take $A\overset{f} \to B'[N]$ on the right-hand-side. Then $f(a_m)=0$ for all $a_m \in A_m$ for all $m  \in {M} \setminus N$, so $f$ vanishes on $J$, 
i.e $f$ factors through $A \to A /J \overset{f}{\to} B'[N]$.
Now we define $F= \Psi (f)$ as the composition $A/J\overset{f}{\to} B'[{N}] \overset{X^n \mapsto 1}{\to} B'$, this is a morphism of $B$-algebras.
Let us prove that $\Theta \circ \Psi = \Id$. Let $f$ be a morphism on the right-hand-side. Let $a_n  \in A_n$ for $n \in N$, we have $f(a_n) = \lambda _n X^n$. We have
\[((\Theta \circ \Psi )(f))(a_n)=(\Theta (\Psi (f)))(a_n)=(\Psi(f))([a_n])\cdot X^n=(f(a_n) ) |_{X^n =1} \cdot X^n = \lambda _n X^n = f(a_n).\] Now let $a_m \in {M} \setminus N$, then 
\[((\Theta \circ \Psi )(f))(a_m)=(\Theta (\Psi (f)))(a_m)=(\Psi(f))([a_m])\cdot X^m=0 =  f(a_m).\] This proves that $\Theta \circ \Psi = \Id$.
Let us prove that $\Psi \circ \Theta = \Id$. Let $F$ be a morphism on the left-hand-side, and let us look at the image of $[a_n]$ for some $a_n \in A_n $ with $n \in N$ under $(\Psi \circ \Theta) (F) = \Psi (\Theta (F)):$
\begin{align*}
 A/J& \xrightarrow{\Theta (F)}B'[N] \to B' \\
 [a_n]&\mapsto F([a_n])X^n \mapsto F([a_n]).
\end{align*}
This finishes the proof of Theorem \ref{representableaffine}.
\xpf

\exam \label{examALSF} \cite[Proposition 2.1]{Ma25} Let $L$ be a submonoid of $M$. If $X= A(L)$ is endowed with the action of $D(M)_S$, then for any submonoid $N$ of $M$ the scheme $X^{N}$ is a closed subscheme of $X$ and equals $  \Spec (\mathbb{Z} [L] / J_N )$ where $J_N $ is the ideal of $\mathbb{Z} [L]$ given by $ \langle \mathbb{Z} X^m | m \in L \setminus (N\cap L) \rangle $.
\xexam

\coro  \label{inclusionclosed} Assume $X$ is affine over $S$. If $N \subset L$ are magnets of $ M$, then $ X^N \to X^L$ is a closed immersion.
\xcoro

\pf With the notation of Theorem \ref{representableaffine}, we have $\calJ _L \subset \calJ _N$.
\xpf 

\rema
Corollary \ref{inclusionclosed} does not generalize outside the $S$-affine case.
\xrema 
\prop  \label{inter} Assume $X$ is affine over $S$. Let $\{N_i\}_{i \in I}$ be magnets of $M$. Then \[\bigcap _{i \in I } X^{N_i} = X^{\cap _{i \in I} N_i} = {\underset{i \in I}{\prod } }{}_X X^{N_i}.\] Here $\cap $ means the scheme theoretic intersection as in  \cite[\href{https://stacks.math.columbia.edu/tag/0C4H}{Tag 0C4H}]{stacks-project} and infinite products of affine morphisms make sense by  \cite[\href{https://stacks.math.columbia.edu/tag/0CNI}{Tag 0CNI}]{stacks-project}.
\xprop

\pf This follows from the description of affine attractors given in Theorem \ref{representableaffine} and the identity $\sum _{i \in I } \calJ _{N_i} = \calJ _{\cap _{i \in I} N_i}$ that follows from the identity $\cup _{i \in I} (M \setminus N_i )= M \setminus (\cap_{i \in I} N_i)$. \xpf

\rema The first equality in Proposition \ref{inter} does not make sense outside the $S$-affine case in general. This is because $X^{N} $ is not a closed subspace of $X$ in general if $X$ is not $S$-affine. Moreover in many non-affine cases $X^{N \cap L} \ne X^N \times _X X^L$ (e.g. cf. \cite[Remark 1.6.3]{Dr15}).
\xrema

        \lemm  \label{equiclosed}  Let $f:Z \to Y$ be an $A(M)_S$-equivariant closed immersion of $S$-affine schemes. Let $N$ be a magnet of $M$. Then the morphism $Z^N \to Y^N$ is a closed immersion, and more precisely $Z^N = Z\times _Y Y^{N}$.      
       \xlemm
       \pf This is local on $S$ so we assume $S=\Spec(R)$ is affine, moreover we identify $Z$ with a closed subscheme of $Y$. Let $A$ be the $R$-algebra of $Y$ and $I$ be the ideal of $A$ defining $Z$.
       Let $J$ be the ideal of $A$ defining $Y^N$, cf. Proof of Theorem \ref{representableaffine} . Then  the ideal of $A/I$ defining $Z^N$ is $I+J/I$ using Theorem \ref{representableaffine}.
       So the ideal of $A$ defining $Z^N$ is $I+J$. Now the isomorphism $A/I \otimes _A A/J \cong A/(I+J)$ finishes the proof.    
       \xpf 
     \rema Proposition \ref{opclosmoo} extends  Lemma \ref{equiclosed} outside the $S$-affine case. Note that the proof of Proposition \ref{opclosmoo} uses Lemma \ref{equiclosed}. \xrema

\prop \label{diagonalizableattra}  Let $\calE$ be a quasi-coherent $A(M)_S$-$\calO _S$-module and let $\mathbb{V} ( \calE ) = \Spec _S ( \Sym ~\calE ) $ be the associated quasi-coherent bundle defined by $\calE$. Then $A(M)_S $ acts linearly on $\mathbb{V} ( \calE )$.  Let $N$ be a magnet of $M$. Then the attractor $\big( \mathbb{V} (\calE) \big) ^N$ associated to $N$ is canonically isomorphic to $\mathbb{V} (\calE ^N)$ where $\calE ^N$ is the $N$-graded component of $\calE$ relatively to the $A(M)_S$-action on $\calE$ (cf. Proposition \ref{4.7.3mono}).
\xprop
\pf Let $p:T \to S$ be a scheme over $S$. We have $\mathbb{V} (\calE )\times _S T = \mathbb{V} (p^* \calE)$. The quasi-coherent $\calO _T$-module inherits a $M$-grading and we have $(p^*\calE )^N = p^* (\calE ^N)$. The following identifications finish the proof
\begin{align*}
\Hom _T^{A(M)_T} (A(N)_T , \mathbb{V} (p^* \calE) ) & = \Hom _{\calO _T\text{-alg}}^{M\text{-graded}} \big( \Sym ~ p^*\calE  , \calO _T[N] \big) \\
& = \Hom _{\calO _T \text{-mod}}^{M\text{-graded}} \big( p^*\calE , \calO _T [N] \big) \\
&= \Hom _{\calO _T \text{-mod}}^{M\text{-graded}} \big( p^*\calE ^N , \calO _T [N] \big) \\
&= \Hom _{\calO _T\text{-mod}} \big( p^* \calE ^N , \calO _T \big) \\
&= \Hom _{\calO _T\text{-alg}} \big(  \Sym ~p^*\calE ^N , \calO _T \big) \\
& = \mathbb{V} (\calE ^N) (T).
\end{align*}
\xpf

\prop \label{monosep} Assume that $X$ is separated and let $N \subset N'$ be two magnets of $M$. Then the map of functors $X^{N} \to X^{N'}$ is a monomorphism, i.e. for any scheme $T$ over $S$ we have a canonical inclusion $X^{N} (T ) \subset X^{N'} (T)$.
\xprop

\begin{proof} Since $N \subset N'$, by \cite[\href{https://stacks.math.columbia.edu/tag/01R8}{Tag 01R8}]{stacks-project}, the scheme theoretic image of $A(N')_T \to A(N)_T $ is $A(N)_T $ for any scheme $T$ over $S$.  Now let $f_1, f_2$ be two elements in $X^{N} (T)$ such that their images in $X^{N'}(T)$ coincide. Consider the schematic kernel of the maps $(f_1,f_2): A(N)_T \to X_T$. Since $X_T $ is separated, $\ker(f_1,f_2)$ is a closed subscheme of $A(N)_T$ (the proof of \cite[Def./Prop. 9.7]{GW} works in this context). So since the scheme theoretic image of $A(N')_T \to A(N)_T $ is $A(N)_T$, we have $\ker (f_1 , f_2) = A(N)_T$ and so $f_1=f_2$.
\end{proof} 

\coro \label{monoX} Assume that $X$ is separated and let $N $ be a magnet of $M$. Then the natural map of functors $X^{N} \to X$ is a monomorphism.
\xcoro 

\pf Combine Proposition \ref{monosep} and Remark \ref{equalM}.
\xpf

\prop \label{propequiv}
Let $Z$ be a monoid and let $f:M \to Z $ be a morphism of monoids. Let $Y $ be a magnet of $Z$ and let $N $ be $f^{-1} (Y)$. Then $N $ is a magnet of $M$. Assume that one of the following conditions holds
\begin{enumerate}
\item $X$ is affine over $S$,
\item $M$ and $Z$ are finitely generated groups,
\item $M$ is cancellative, $Z=M^{\mathrm{gp}}$, $N=Y$ and $X$ is separated,
\end{enumerate} 
then in each case we have an identification of functors $X^{N}= X^{Y}$, where on the left-hand-side $X$ is seen as an $A(M)_S$-space and on the right-hand-side as an $A(Z)_S$-space (via $A(Z)_S \to A(M)_S$ dual to $f: M \to Z$). In other words, with the notation of Remark \ref{Notationavecindicegroup}, we have $X^{N \subset M}=X^{Y \subset Z}$.
\xprop

\pf \begin{enumerate} \item  We reduce to the case where $S$ and $X= \Spec( A)$ are affine. We have two compatible gradings on $A$, one given by $Z$ and one given by $M$. For any $y \in Y$, we have $A_y = \oplus _{n \in f^{-1} (y)} A_n$.
So $\oplus _{z \in Z \setminus {Y}} A_z= \oplus _{n \in M \setminus N} A_n$.
Then the ideal defining $X^N $ equals the ideal defining $X^Y$, cf. Theorem  \ref{representableaffine}.
\item Let $T \to S$ be a scheme. We have a canonical map \[\theta: \Hom _T ^{D(M)_T} ( A(N)_T , X_T) \to \Hom ^{D(Z)_T} (A(Y)_T,X_T )\] obtained by precomposition with $A(Y)_T \to A(N)_T$. Let us show that $\theta$ is an isomorphism. We construct the reciprocal map. Let $f : A(Y)_T \to X_T$ be a $D(Z)_T$-equivariant map. We get a map $f':D(M)_T \times _T A(Y)_T \to X_T $ given on points by $(g,x) \mapsto g \cdot f(x)$. Let $K$ be $f(M)+ Y^{\mathrm{gp}}$, this is a subgroup of $Z$. We have a morphism of groups $\phi : f(M) \times Y^{\mathrm{gp}} \to  f(M) + Y^{\mathrm{gp}} $ given by $(f(m),y) \mapsto f(m) +y$. The morphism $\phi$ induces a closed immersion of group schemes $D( f(M) + Y^{\mathrm{gp}} )_T \to D(f(M) \times Y^{\mathrm{gp}})_T = D(f(M))_T \times _T D( Y^{\mathrm{gp}})_T$. Let us consider the action of $D(f(M))_T \times _T D(Y^{\mathrm{gp}})_T $ on $D(M)_T \times _T A(Y)_T $ given by $(g,h) \cdot (g',x)= (g \cdot g' , h^{-1} \cdot x) $, remark that this action is free. So we obtain by composition a free action $\star$ of $D(K)_T$ on $D(M)_T \times _T A(Y)_T $.
Let us consider the morphisms
\[\xi_{\star}, \xi : D(K)_T \times _T D(M)_T \times _T A(Y)_T \to X_T \] where $\xi_{\star}$ is given by $(k,g,y) \mapsto f'(k \star (g,y) )$ and $\xi$ is given by $(k,g,y) \mapsto f'(g,y )$. Let us prove that $\xi_{\star}= \xi$. Let us consider the canonical morphism  \[p: D(Z)_T \times _T D(M)_T \times _T A(Y)_T \to D(K)_T \times _T D(M)_T \times _T A(Y)_T \] induced by the inclusion $K \subset Z$.
Then $p \circ \xi = p \circ \xi _{\star}$, and so $\ker (\xi , \xi _{\star} )=D(K)_T \times _T D(M)_T \times _T A(Y)_T $ because the schematic image of $p$ is  $D(K)_T \times _T D(M)_T \times _T A(Y)_T$. This finishes the proof of the claim $\xi = \xi_{\star}$.
So we have $f'(k \star (g,x))=f'((g,x))$ for all $T'\to T$, $k \in D(K)_T(T')$, $g \in D(M)_T(T')$ and $x \in A(Y)_T (T')$.
 So the map $f'$ induces by factorization a map $(D(M)_T \times _T A(Y)_T ) / D(K)_T \to X_T$ where $(D(M)_T \times _T A(Y)_T ) / D(K)_T $ is the fpqc quotient as in \cite[VIII Th. 5.1]{SGA3}. Now let $\mathcal{A}$ be the quasi-coherent algebra of the $T$-affine scheme $(D(M)_T \times _T A(Y)_T ) $. Consider the $K$-grading on $\mathcal{A}$ associated to the action of $D(K)_T $ on  $(D(M)_T \times _T A(Y)_T )$. By \cite[VIII Th. 5.1]{SGA3}, $(D(M)_T \times _T A(Y)_T ) / D(K)_T$ is affine with quasi-coherent algebra the degree zero part $\mathcal{A}_0$ in $\mathcal{A}$.
The $K$-grading on $\mathcal{A}$ is given locally by $\mathrm{deg}(X^{(m,y)})=f(m)-y  ~\in K$. This implies that $\mathcal{A}_0 \simeq  \calO _T [f^{-1} (Y)]$. This identifies $(D(M)_T \times _T A(Y)_T ) / D(K)_T$ with $A(N)_T$. So the map $f$  induces a map $A(N)_T \to X_T$. This map is $D(M)_T$-equivariant. The obtained map $ \Hom ^{D(Z)_T} (A(Y)_T,X_T ) \to \Hom ^{D(M)_T} ( A(N)_T , X_T)$ is the reciprocal map of $\theta$.
\item We have to show that for any $T/S$, $\Hom _T ^{D(M^{\mathrm{gp}})_T} (A(N)_T , X_T) = \Hom _T^{A(M)_T} (A(N)_T , X_T)$. The inclusion $\supset$ is clear. Reciprocally let $\phi : A(N) _T \to X_T$ be $D(M^{\mathrm{gp}})_T$-equivariant. We use that $X$ is separated and that the schematic image of \[A(N)_T \times _T {D(M^{\mathrm{gp}})_T} \to A(N)_T \times _T A(M)_T\] is $A(N)_T \times _T A(M)_T$ to prove that $\phi $ is $A(M)_T$-equivariant (cf. the proof of Proposition \ref{monosep} for similar arguments). \end{enumerate}
\xpf

\prop \label{NLNL}  Let $L$ and $N$ be arbitrary submonoids of $M$. Assume that one of the following conditions holds
\begin{enumerate}
\item $X$ is affine over $S$,
\item $M$ is a finitely generated abelian group,
\item $M$ is cancellative, $M^{\mathrm{gp}} $ is finitely generated as abelian group, and $X$ and $X^N$ are separated algebraic spaces,
\end{enumerate} 
then we have canonical identifications \[(X^N)^L= X^{N \cap L},\] cf. \ref{actiongroup} for the actions of $A(M)_S$ on $X^N$ and $X^L$ that we used implicitely on the left-hand sides.
\xprop

\pf \begin{enumerate} \item We use the explicit description given in Theorem \ref{representableaffine}.
\item 
 Let us use \cite[VIII Th. 5.1]{SGA3} as follows. We remark that \[(X^N)^L (T)= \Hom _T ^{D(M)_T \times _T D(M)_T } (A(N)_T \times _T A(L)_T , X_T)\] where $D(M)_T \times _T D(M)_T $ acts on $X_T$ via $(g,h) \cdot x = g \cdot (h \cdot x)$ and on $ A(N)_T \times _T A(L)_T $ via $(g,h)\cdot ( a,b) = (g \cdot a , h \cdot b)$.
Let us consider the action $\star$ of $D(N^{\mathrm{gp}} + L^{\mathrm{gp}})_T $ on $A(N)_T \times _T A(L)_T $ given by $g \star (a,b)=(g^{-1} \cdot a,g \cdot b)$. The action $\star$ is free because firstly the action of $D(N^{\mathrm{gp}} \times L^{\mathrm{gp}} )_T$ on $A(N \times L)_T$ given by $(\lambda , \beta ) \cdot (a,b) = (\lambda \cdot a, {\beta }^{-1} \cdot b ) $ is free and secondly because the morphism of groups $N^{\mathrm{gp}} \times L^\mathrm{gp} \to N^\mathrm{gp} + L ^{\mathrm{gp}} $ given by $(n,l) \mapsto n+l$ is surjective and so $D(N^{\mathrm{gp}}+L^{\mathrm{gp}})_T $ is a closed subgroup  scheme of $D(N^{\mathrm{gp}} \times L^{\mathrm{gp}})_T$. Now let $F \in  \Hom _T ^{D(M)_T \times _T D(M)_T } (A(N)_T \times _T A(L)_T , X_T)$, $F : A(N)_T \times _T A(L)_T \to  X_T$. Then ${F(g \star (a,b) ) = F((a,b))}$ for any $T' \to T$, $g \in D(N^{\mathrm{gp}}+L^{\mathrm{gp}})_T(T')$ and $(a,b) \in (A(N)_T \times _T A(L)_T)(T')$. So $F $ induces a morphism $f$ from the fpqc quotient, $f:\big( A(N)_T \times _T A(L)_T \big) / D(N^{\mathrm{gp}}+L^{\mathrm{gp}})_T \to X_T$. We have $A(N)_T\times _T A(L) _T  = \Spec ( \bbZ [ N \times L] ) \times _{\bbZ}T$
 and the degree zero part of $\bbZ [ N \times L] $ (relatively to the $N^{\mathrm{gp}}+L^{\mathrm{gp}}$-grading induced by the action $\star$ of $D(N^{\mathrm{gp}}+L^{\mathrm{gp}})$) is $ \bbZ [ N \cap L]$. So we have an identification $\big( A(N)_T \times _T A(L)_T \big) / D(N^{\mathrm{gp}}+L^{\mathrm{gp}})_T = A(N \cap L ) _T$.  So $F$ induces a natural morphism $f:A(N \cap L )_T \to X_T$, that is $D(M)_T$-equivariant. Now let $f \in X^{N \cap L}(T)$, $f : A(N \cap L)_T \to X _T$ and consider the composition 
 \[ F : A(N)_T \times _T A(L)_T \to A (N \cap L) _T \times _T A(N \cap L)_T \to A(N \cap L)_T \xrightarrow{f} X_T \] where the second morphism is the multiplication of $A (N \cap L) _T$. Then $F$ is $D(M)_T \times _T D(M)_T$-equivariant. The previous maps $F \mapsto f$ and $f \mapsto F$ induce a canonical bijection between \[{ \Hom _T ^{D(M)_T \times _T D(M)_T } (A(N)_T \times _T A(L)_T , X_T)} \] and $\Hom _T ^{D(M)_T}( A(N \cap L )_T , X_T )$.
 \item Using three times Proposition \ref{propequiv} and one time (ii), we have \[(X^{N \subset M})^{(L \subset M)}=
  (X^{N \subset M^{\mathrm{gp}}})^{L\subset M^{\mathrm{gp}}}= X^{(L \cap N )\subset M^{\mathrm{gp}}}= X^{(N \cap L ) \subset M}.\]\end{enumerate}
\xpf

\exam
Assume $M = \mathbb{Z}^2$. Let $a$ be the action of $D(M)_S$ on $X=A(\mathbb{N}^2)_S$ corresponding to the inclusion $\mathbb{N}^2 \subset \mathbb{Z}^2$. In other words, $a$ is the action of $\mathbb{G}_{m,S}^2$ on $X=\mathbb{A}_S^2$ given by $(\lambda , \beta ) \cdot ( x,y) = ( \lambda x , \beta y )$.
Then $X^{\mathbb{N} \times 0} = \mathbb{A}_S^1 \times _S 0_{\mathbb{A}_S^1} \subset \mathbb{A}_S^2$, and $(X^{\mathbb{N} \times 0} )^{0 \times \mathbb{N}} = 0_{\mathbb{A}_S^2} = X^{(\mathbb{N} \times 0) \cap (0 \times \mathbb{N})}$ (e.g. by Example \ref{examALSF}).
\xexam

\fact \label{intersectioninfini}
Let $\{N_i\}_{i\in I} \subset M$ be magnets such that $X^{N_i}=X^{N_j}$ for all $i,j \in I$. Assume that $X$ is separated and that one of the following conditions holds
\begin{enumerate}
\item $X$ is affine over $S$,
\item $M$ is a finitely generated abelian group,
\item $M$ is cancellative, $M^{\mathrm{gp}} $ is finitely generated as abelian group, and $X^{\cap_{l \in J}N_l}$ is a separated algebraic spaces for any finite subset $J \subset I$,
\end{enumerate}then, for all $i \in I$,  \[X^{N_i}= X^{\cap_{l\in I} N_l}.\]
\xfact
\pf\begin{sloppypar} If $X$ is affine, we use Proposition \ref{inter}. Now we prove (ii) and (iii). The inclusion $\cap_{l \in I } N_l \subset N_i$ induces a canonical morphism $X^{\cap_{l \in I } N_l} \to X^{ N_i}$. Since $X$ is separated, this canonical morphism is a monomorphism by Proposition \ref{monosep}. So for any scheme $T \to S$, $X^{\cap_{l \in I } N_l} (T) \subset X^{ N_i} (T)$. We now prove the reverse inclusion.
 We write $I $ as a directed colimit of finite sets $J$. Then, in the category of sets,  $\cap _{l \in I} N_l = \lim _{J \subset I} \cap _{l \in J} N_l$. Furthermore, in the category of rings, $\bbZ [ \cap _{l \in I} N_l ] = \lim _{J \subset I} \bbZ [ \cap _{l \in J} N_l ]$. So $A(\cap _{l \in I} N_l )_S = \colim _{J \subset I } A(\cap _{l \in J} N_l )_S $. By Proposition \ref{NLNL}, $X^{\cap _{l\in J} N_l} =X^{N_i}$ for any finite subset $J$ of $I$. We use these observations to obtain
 \begin{align*}
 X^{ \cap _{l \in I} N_l} (T) &= \Hom _T^{A(M)_S}( A(\cap _{l \in I} N_l )_T , X_T) \\
 &= \Hom _T^{A(M)_S}(  \colim _{J \subset I } A(\cap _{l \in J} N_l )_T , X_T) \\
\text{ \cite[\href{https://stacks.math.columbia.edu/tag/002H}{Tag 002H}]{stacks-project} and equivariance is clear}&\supset  \lim _{ J \subset I} \Hom _T^{A(M)_S}(  A(\cap _{l \in J} N_l )_T , X_T) \\
 &= \lim _{J \subset I} X^{\cap _{l \in J} N_l} (T)\\
 &= \lim _{J \subset I} X ^{N_i}(T)\\
 &= X^{N_i}(T).
 \end{align*} This finishes the proof.\end{sloppypar}
\xpf
\prop \label{Group} Assume $M$ is a finitely generated group.
Let $Z \subset M$ be a subgroup. Then the attractor space $X^Z$ for the action of $D(M)_S$ on $X$ is identified with the fixed-points space $X^{D(M/Z)_S}$ of $X$ under the action of $D(M/Z)_S$.
\xprop
\pf We have an exact sequence of abelian groups 
\[ 0 \to Z \to M \to M/Z \to 0 . \] 
By \cite[Exp. VIII]{SGA3} we obtain an exact sequence of diagonalizable group schemes 
\[ 1 \to D(M/Z)_S \to D (M) _S \to D(Z)_S \to 1 .\]
Let $T \to S$ be a scheme and let us prove that $X^Z(T) = X^{D(M/Z)_S} ( T)$.
Note that we have a $D(M)_S$-equivariant identification of $S$-scheme $A(Z)_S =D(Z)_S$.
We have 
\begin{align*} 
X^Z (T) &= \Hom ^{D(M)_T} _{T} ( A(Z)_T , X_T ) \\
&= \Hom ^{D(M)_T}_{T} (D(Z)_T , X_T ) \\
&= \Hom ^{D(M) _T }_T ( D(M)_T / D(M/Z) _T , X_T ) \\
&= \Hom _T ^{D(M)_T} ( D(M)_T ,  (X_T)^{{D(M/Z)_T}} ) .
\end{align*}
Now Remark \ref{equalM} finishes the proof.
\xpf

 \prop \label{cartesian}
  Assume that $M$ is an abelian group and that $X$ is a scheme. Let $N,N',L$ and $L'$ be magnets in $M$ such that $L \subset L'$, $N' \subset L', N = L \cap N'$ and $L'=L+N'$. Assume that $X^E$ is representable by a scheme for all $E \in \{ N,N',L,L'\}$. Assume that  $L$ is a face of $ L'$. Then  $N $ is a face of $N'$. Assume that one of the following conditions holds
 \begin{enumerate}
 \item we have an equality $N' = L' \setminus ( L \setminus N) $,
 \item $S=\Spec (R)$ and $X= \Spec (A)$ are affine, $A_{l}A_{n'}=A_{l+n'}$ for all $l \in L\setminus N $ and $n' \in N'$ (as usual ${A}_m$ denotes the $m$-graded part of $A$),
 \end{enumerate} 
  then the following diagram is a cartesian square in the category of schemes  \[ \begin{tikzcd}                                                                                         & X^{N'} \ar[dl,    "\iota _{N',L'} "']  \ar[rd, "\mathfrak{p}_{N,N'}"']& \\
                                                                                         X^{L'} \ar[rd, "\mathfrak{p}_{L,L'} "] & & X^N\ar[ld,  "\iota _{N,L}"]\\ 
                                                                                         & X^L   & 
                                                                                        \end{tikzcd} . \]
 \xprop
 \pf The monoid $ N $ is clearly a face of $ N'$.
\begin{enumerate}  \item 
  Let $T$ be a scheme and let $T \to X^{L'}$, $T \to X^N$ be two morphisms of schemes such that the following diagram commutes  \[ \begin{tikzcd}
                                                                                         & T \ar[dl]  \ar[rd]& \\
                                                                                         X^{L'} \ar[rd, "\mathfrak{p}_{L,L'} "'] & & X^N\ar[ld,  "\iota _{N,L}"]\\ 
                                                                                         & X^L   & 
                                                                                        \end{tikzcd}  .\]                                                                                        This corresponds to a diagram                                                                                      
                                                                                \[\begin{tikzcd}  & A (L) _T \ar[rd] \ar[dl] & \\ A (L')_T\ar[rd] & & A (N)_T \ar[dl] \\ & X_T & \end{tikzcd} \]
                                                                                        where all arrows are $D(M)_T$-equivariant. By Lemma \ref{lemmpush}, we obtain a unique arrow $A(N')_T \to X_T$ such that the following diagram commutes                                                                      
\[\begin{tikzcd}  & A (L)_T \ar[rd] \ar[dl] & \\ A (L')_T \ar[rd]\ar[rdd] & & A (N)_T \ar[dl]\ar[ddl] \\ 
             & A(N')_T\ar[d] &\\
             & X_T & \end{tikzcd} .\]       
It is enough to show that the arrow $A(N')_T \to X_T$ is $D(M)_T$-equivariant.
Consider the diagram    \begin{tiny}
    \[\begin{tikzcd}  & A (L)_T\times _T D(M)_T \ar[rd] \ar[rrr] \ar[dl] & && A (L) _T\ar[rd] \ar[dl] & \\ A (L')_T\times _T D(M)_T  \ar[rrr, bend left] \ar[rd]\ar[rdd] & & A (N)_T\times _T D(M)_T  \ar[rrr, bend left=50] \ar[dl]\ar[ddl] & A (L')_T\ar[rd]\ar[rdd] & & A (N)_T\ar[dl]\ar[ddl] \\ 
             & A(N')_T\times _T D(M)_T \ar[rrr]\ar[d] && & A(N')_T\ar[d] &\\
             & X_T \times _T D(M)_T \ar[rrr]& & & X_T & \end{tikzcd} \]        
                                                \end{tiny}obtained by fiber product with $D(M)_T$. We have  \[A(E)_T\times _T D(M)_T = (A(E) \times _{\Spec(\bbZ)} T ) \times _T D(M)_T = A(E)_{D(M)_T}\] for $E \in \{ N,N',L,L'\}$, so the left diamond is a push-out by Proposition \ref{lemmpush}. Now we want to show that the lower rectangle is commutative. Consider the upper right composition in this rectangle and precompose it with the right part of the left diamond, denote this arrow by $a_1$. Consider the lower left composition in the rectangle and precompose it with the left part of the left diamond, denote this arrow by $a_2$. Now using the commutative diagrams coming from the $D(M)_T$-equivariant morphisms on the right, we see that $a_1$ and $a_2$ are both equal to the composition \[A(L)_T \times _T D(M)_T  \to A (L)_T \to A(N')_T \to X_T .\] Using the left push-out diamond, this now implies that the lower rectangle is commutative. So the arrow $A(N')_T \to X_T$ is $D(M)_T$-equivariant. This finishes the proof.                                                                                                                                                                                                                                                                                                                                                                                                                                                                                                                                                                                                                                                                                                                                                                               \item Let $x \in X^{N'}(R)$.
 Then $x$ is a morphism $A \to R[N']$. Now we have an equality of compositions 
 \[ (A \to R[N'] \to R[L'] \to R[L]) = (A \to R[N'] \to R[N' \cap L] \to R[L] ).\] This shows that the diagram is commutative. 
 Now let us prove that it is cartesian. Let $Y=\Spec(B)$ be an affine $R$-scheme and let $f:Y \to X^{L'}$ and $g:Y \to X^{N}$ be two morphisms such that $\mathfrak{p}_{L,L'} \circ f= \iota _{L,N} \circ g $.
 So $f$ is a morphism of graded algebras $A \to B[L']$ and $g$ is a morphism of graded algebras $A \to B [N]$.
 Let $m \in L'=L+N'$ and let $A_m$ be the $m$-graded part of $A$.
 Let $x \in A_m$ and let $\lambda _m$ such that $f(x)=\lambda _m X^m$. Then since $\mathfrak{p}_{L,L'} \circ f= \iota _{L,N} \circ g $, we obtain that $\lambda _m=0$ for all $l \in L \setminus ( N' \cap L)$.
 So we get $f(x)=0$ for all $x \in A_m$ for all  $m \in L' \setminus N'$ (we use that $A_{l}A_{n'}=A_{l+n'}$ for all $l \in L\setminus (N' \cap L) $ and $n' \in N'$). So we obtain a unique morphism $h$ from $Y $ to $X^{N'}$ with the cartesian property. 
 \end{enumerate}
 \xpf

\prop \label{colimcolim} Assume that $X$ is separated and locally of finite presentation over $S$. Let $N $ be a magnet of $M$. Write $N = \colim _{i\in I} N_i$ as a directed colimit of submonoids of $N$. Let $T \to S$ be a scheme, then 
\[X^{N}(T) = \colim _{i \in I} X^{N_i} (T). \]
\xprop

\pf  For any scheme $T$, we have $ A(N) _T = \lim _{i \in I} A(N_i) _T $ by \cite[\href{https://stacks.math.columbia.edu/tag/01YW}{Tag 01YW}]{stacks-project}. So $X( A(N)_T ) = \colim X (A(N_i ) _T )$ by  \cite[\href{https://stacks.math.columbia.edu/tag/049J}{Tag 049J}]{stacks-project}. Let $A(N)_T \to X_T$ be an $A(M)_T$-equivariant morphism. By the previous assertion this factorizes through a morphism $A(N_i)_T \xrightarrow{f} X_T $. We want to show that $f$ is $A(M)_T$-equivariant. Consider the diagram \[\xymatrix{A(M)_T \times A(N_i)_T \ar@<-.5ex>[r] \ar@<.5ex>[r] & X_T \\ A(M)_T \times A(N)_T \ar[u] \ar[ur] }\] where the horizontal arrows $(f_1,f_2)$ correspond to $(g,x) \mapsto f(gx)$ and $(g,x) \mapsto gf(x)$. Since $X$ is separated, the kernel of $(f_1,f_2)$ is a closed subscheme of $A(M)_T \times A(N_i)_T$ (the proof of \cite[Def./Prop. 9.7]{GW} works in this context). Since the schematic image of the vertical morphism $\phi$ is $A(M)_T \times _T A(N_i)_T$ and because $f_1 \circ \phi = f_2 \circ \phi$, we have an equality of schemes $\ker (f_1 , f_2)=A(M)_T \times _T A(N_i)_T$. So we have an equality of morphisms of schemes $f_1 =f_2$, and so $f$ is equivariant.
We deduce that 
$X^N (T) = \colim_{i \in I}  X^{N_i} (T) . $
\xpf

 \section{Attractors with prescribed limits} \label{prescribed}
 
We introduce in this section another functor. 
Let $X$ be an algebraic space over a base scheme $S$. Let $M$ be a commutative monoid and let $A(M)_S$ be the associated diagonalizable monoid scheme over $S$.
Assume that $ A(M)_S$ acts on $X$. Let $N \subset M$ be a magnet. Let $F$ be a face of $N$. Let $Z $ be an other $S$-functor with a monomorphism $Z \to X^{F}$.
We now introduce the attractor $X_{F,Z}^N$ associated to the magnet $N$ under the action of $A(M)_S$ on $X$ with prescribed limit in $Z$ relatively to the face $F$. Recall that we have a canonical morphism $X^N \to X^F $ (cf. Fact \ref{faceattractors}). Since $Z(T) \subset X^{F} (T)$ for any $S$-scheme $T$, the following definition makes sense.

\defi \label{definitionfonct} Let $X_{F,Z}^N$ be the contravariant functor \[(Sch/S) \to Set , (T \to S) \mapsto \{ f \in X^N (T) \mid \text{ the image of } f \text{ in } X^{F}(T) \text{ belongs to } Z (T)  \} .\] If $F=N^*$, we omit $F$ in the notation, i.e. we put $X^N_{N^*,Z}=: X^N_Z$.
\xdefi

\prop We have a canonical isomorphism $X^N_{F,Z} \simeq X^N \times _{X^{F}} Z $.
\xprop 

\pf Clear since
\begin{align*} X_{F,Z}^N (T) 
& =\{ f \in X^N (T) \mid \text{ the image of } f \text{ in } X^{F}(T) \text{ belongs in } Z (T)  \}\\
&= \{ (f , g )\in X^N (T) \times Z(T) \mid ~ f=g \text{ in } X^{F} (T) \} \\
&=   X^N (T) \times _{X^{F}(T)} Z(T)  .
\end{align*}
\xpf

\fact Assume that $X$ is separated. Then we have a canonical monomorphism $X^N _{F,Z} \to X^N$ of $S$-functors. In particular we have a canonical monomorphism $X^N_{F,Z} \to X$ of $S$-functors.
\xfact
\pf
Clear by Corollary \ref{monoX}.
\xpf

 \section{Base change of affine étale morphisms along faces} \label{basechange}
In this section, we prove a single specific result that will be applied throughout the remaining parts of the paper.
 Let $S$ be a scheme and let $X$ be an
$S$-algebraic space with an action of $D(M)_S$ for some finitely generated abelian group $M$. Let $N$ be a magnet of $M$.

\prop \label{fixxx} Let $Q$ be a face of $N$ and let $X^N \to X^Q$ be the associated morphism on attractors. Let $U$ be an $S$-affine scheme.
Let $U\to X$ be an \'etale, $D(M)_S$-equivariant morphism of
algebraic spaces. Assume that one of the following conditions hold
\begin{enumerate}
\item $N$ is finitely generated as monoid,
\item $X$ is separated,
\end{enumerate} then the natural map $U^N\to X^N\times_{X^Q} U^Q$ is an isomorphism.
\xprop

\pf We need some terminologies on monoids. Recall that a monoid is fine if it is cancellative and finitely generated. Recall also that a monoid $L$ is sharp if $L^*=0$. See \cite[§I]{Og} for more details on monoids. \begin{enumerate} \item  Firstly, we remark that it is enough to treat the case $Q=N^*.$ Indeed assume that Proposition \ref{fixxx} is true for the face of invertible elements. Since $Q^*= N^*$, we have 
\[X^N \times _{X^Q} U^{Q} = X^N \times _{X^Q} (X^Q \times _{X^{N^*}} U^{N^*} )= X^N \times _{X^{N^*}} U ^{N^*} = U^N .\]
So we assume $Q=N^*$. We now remark that we can assume $N^* =0$ and that $N$ is fine and sharp using the map $M \to M/N^{*}$ and Proposition \ref{propequiv}. So we now assume that $Q=0$ and $N$ is fine and sharp and we adapt \cite[Lemma 1.11]{Ri16}.
It is enough to construct the inverse morphism
$X^N\times_{X^0} U^0\to U^N$.
For this let $p:T \to S$ be an $S$-scheme and let
\[
\begin{tikzcd}
T=A(0)_T \ar[r,"{f_0}"] \ar[d] & U_T \ar[d] \\
A(N)_T \ar[r] & X_T
\end{tikzcd}
\]
be a diagram corresponding to a $T$-point of $X^N\times_{X^0} U^0$. We
want to find a diagonal filling $A(N)_T\to U_T$. Let $I$ be the kernel of $\calO _T [N] \to \calO _T[0]$, i.e the ideal associated to $N \setminus 0 $ in $\calO _T [N]$. Let $V(I^k)= \Spec _T (\calO _T [N] /I^k ) $ be the
infinitesimal neighbourhoods of $A(0)_T$ inside $A(N)_T$. Since $U_T\to X_T$
is \'etale, it is smooth by  \cite[\href{https://stacks.math.columbia.edu/tag/04XX}{Tag 04XX}]{stacks-project}, and formally smooth by \cite[\href{https://stacks.math.columbia.edu/tag/02H6}{Tag 02H6}]{stacks-project}, so by the infinitesimal lifting property the morphism
$f_0:T\to U_T$ lifts uniquely
to a compatible family of morphisms $f_k:V(I^k)\to U_T$. These liftings
are equivariant because the two maps $D(M)_T\times V(I^k)\to U_T$,
$(g,x)\mapsto f_k(gx)$ and $(g,x)\mapsto gf_k(x)$ are common liftings
of the map $D(M)_T\times A(0)_T\to U_T$, $(g,x)\mapsto f_0(gx)=gf_0(x)$
hence by uniqueness they are equal.
Writing $U=\Spec_S(\calA)$ for a $M$-graded quasi-coherent $\calO_S$-algebra~$\calA$, we have and $U\times _S T = \Spec _T (\calA _T)$ where $\calA _T = p^* \calA$. 
Moreover we have a family of $M$-graded morphisms
$\calA_T\to \calO_T[N]/I^k$. This means that ${(\calA _T)}_m$ goes to $0$ when
$m\not\in N$ and to $(\calO_T[N]/I^k)_m$ when $m\in N$.
The induced morphism to the completion 
$\calA _T \to \calO_T[[N]]$ has image in $\calO_T[N]$, yielding the desired lifting $A(N)_T\to U_T$ (cf. \cite[Chap. I,§3 Prop. 3.6.1]{Og} for the local description of the completion $\calO_T[[N]]$).
\item Let $T \to S $ be a scheme. Write $N = \colim _{i \in I} N_i$ as a colimit of finitely generated monoids. Using (i) and Proposition \ref{colimcolim} we have 
\begin{align*}
U^N (T) &= \colim _{i \in I } U^{N_i} (T) \\
& = \colim _{i \in I } \big(X^{N_i} (T) \times _{X^Q (T)} U^Q (T) \big) \\
&= \colim _{ i \in I}  \big(X^{N_i} (T)\big) \times _{X^Q (T)} U^Q (T) \\
&= X^N (T) \times _{X^Q (T)} U^Q (T) \\ 
& = (X^N \times _{X^Q} U^Q ) (T).
\end{align*}
\end{enumerate}
\xpf

\section{Fixed-point-reflecting atlases} \label{FPR}

Let $M$ be an abelian group. Let $S$ be a scheme. Let $X$ be a quasi-separated $S$-space endowed with an action of $D(M)_S$. 
We discuss fixed-point-reflecting atlases in this section. Roughly, fixed-point-reflecting étale atlases are atlases such that we can learn from the charts the fixed points of our space. This concept makes sense for an action of an arbitrary group scheme and is well-studied. In this section we introduce several related concepts and definitions, $Z$-FPR morphisms and strongly-FPR morphisms, useful and motivated by algebraic magnetism ($Z$ is a subgroup of $M$ here). A strongly-FPR atlas is $Z$-FPR for all subgroups $Z \subset M$. We formulate a theorem (Theorem \ref{adapted-atlas}) that provides the existence of $Z$-FPR atlases and a conjecture (Conjecture \ref{conjecture}) about the existence of strongly-FPR atlases. In this work, $Z$-FPR atlases are used to prove representability results (Proposition \ref{repgroup} and Theorem \ref{representable}) and strongly-FPR atlases are needed to study magnets (Theorem \ref{strumagnet}).
 In the next section (Section \ref{Zariskicase}), we will notice that Conjecture \ref{conjecture} is true in many cases (namely for Zariski locally linearizable actions, e.g. Sumihiro's actions \ref{remasum}) in particular one has $Z$-FPR atlases in this case. Romagny's appendix offers a proof of Theorem \ref{adapted-atlas}.

\defi  Let $f: U \to X$ be a $D(M)_S$-equivariant morphism of $S$-spaces. \label{deffpratl}\begin{enumerate}\item We say that $f$ is fixed-point reflecting (FPR) if the canonical morphism of functors $U^{D(M)_S} \to U \times _X X^{D(M)_S} $ is an isomorphim. 
\item We say that $f$ is $Z$-FPR if it is fixed-point reflecting for the induced action of $D(M/Z)_S$, i.e. if the canonical morphism of functors $U^{D(M/Z)_S} \to U \times _X X^{D(M/Z)_S} $ is an isomorphim. 
\item We say that $f$ is strongly-FPR if it is $Z$-FPR for all subgroups $Z \subset M$. 

\item  We say that $f$ is an equivariant atlas if $f$ is étale and surjective and $U= \coprod_{\tau \in  \mathscr{A}} U_{\tau}$ is the disjoint union of $D(M)_S$-stable $S$-affine schemes; moreover in this case we say that the equivariant atlas $f$ is $S$-affine if  $\mathscr{A}$ may be chosen finite.
\item  We say that $f$ is a $Z$-FPR atlas if $f$ is a $Z$-FPR equivariant atlas.

\item We say that $f$ is a strongly-FPR atlas if $f$ is a strongly-FPR equivariant atlas.
\end{enumerate}
\xdefi

\rema Since a finite disjoint union of $S$-affine schemes is $S$-affine, if an equivariant atlas is $S$-affine then we can assume that $\mathscr{A}$ is a singleton. In other words an $S$-affine equivariant atlas for $X$ is a $D(M)_S$-equivariant étale surjective morphism $U \to S$ where $U $ is an $S$-affine scheme.
\xrema

The following is an immediate generalization of \cite[Lemma 1.10]{Ri16}.
\prop \label{lemmpf}Assume that $M \cong \bbZ ^r$ for some integer $r$ and let $U \to X$ be an étale $D(M)_S$-equivariant morphism, then $U^{D(M)_S} \to X^{D(M)_S} \times _X U $ is an isomorphism. In other words, if $M$ is torsion-free, then every étale morphism is FPR.
\xprop 

\pf The proof of \cite[Lemma 1.10]{Ri16} works replacing ``$\mathbb{G} _m$" by ``$D(\bbZ^r)$". 
\xpf

\rema \label{rem1}In his appendix, Romagny provides a generalization of Proposition \ref{lemmpf} to arbitrary flat group schemes with connected fibers and in particular provides a detailled proof of Proposition \ref{lemmpf}, cf. Theorem \ref{G_connected_implies_fixed_pts_closed} .
\xrema

\rema If $M$ is not torsion-free, and $U \to X$ is an étale $D(M)_S$-equivariant morphism, then $U^{D(M)_S} \to X^{D(M)_S} \times _X U $ is obviously not an isomorphism in general. For example, take $M= \bbZ/2 \bbZ$, $S= \Spec (\bbZ)$, $X=S$ endowed with the trivial action of $D(M)$ and $U = D (M ) $ endowed with the non-trivial action of itself by multiplication. Then the $D(M)$-equivariant morphism $U \to S$ is smooth of relative dimension zero and so it is étale; but $U^{D(M)}= \emptyset $, $X^{D(M)}= X $ and therefore  $U^{D(M)} \not\cong X^{D(M)} \times _X U $. 
\xrema

\rema Let us provide an example of a space $X$ endowed with an action of a diagonalizable group $D(M)$ and an FPR equivariant atlas $U \to X$ that is not strongly-FPR. Choose an exact sequence of abelian groups $0 \to Z \to M \to M/Z \to 0$ with $Z \ne M$ and $M/Z$ finite. Then $f: D(M) \to D(Z)$ is étale, surjective and $D(M)$-equivariant. Put $U = D(M) $ and $X= D(Z)$. Since $U^{D(M)} = X^{D(M)} = \emptyset$, $f$ is FPR. Now $f$ is not $Z$-FPR because $U^{D(M/Z)} = \emptyset $ and $X^{D(M/Z)}  =X$.
\xrema

\lemm \label{clofpr} Let $X,X'$ be $S$-algebraic spaces endowed with actions of $D(M)_S$. Let $Z$ be a subgroup of $M$.  Let $X' \to X$ be a $D(M)_S$-equivariant affine morphism of $S$-algebraic spaces. Let $U \to X$ be a $Z$-FPR atlas of $X$. Then 
\begin{enumerate}
 \item The canonical morphism $ U \times _X X' \to X' $ is a $Z$-FPR atlas of $X'$.
\item If  $U \to X$ is an $S$-affine $Z$-FPR atlas of $X$, then  $U \times _X X' \to X' $ is an $S$-affine $Z$-FPR atlas of $X'$. \end{enumerate}
\xlemm 

\pf (i) Since $X' \to X$ is an affine morphism, $U \times _X X' \to U$ is affine and so the composition $U \times _X X' \to S$ is a disjoint union of $S$-affine schemes. The morphism $U \times _X X' \to X'$ is étale and surjective because $U \to X$ is so. By Proposition \ref{produitNNN}, we have canonical identifications
\[ (X' \times _X U )^{Z} = {X'}^{Z} \times _{X^{Z} } U^{Z} = {X'}^{Z} \times_{ X^{Z} }(U \times _X X^{Z}) = {X'}^{Z} \times _X U= {X'}^{Z} \times _{X'} ( X' \times _X U). \]Assertion (ii) is now immediate.
\xpf

 Theorem \ref{adapted-atlas} was motivated by finding a generalization of Proposition \ref{lemmpf} for general diagonalizable group schemes. Theorem \ref{adapted-atlas} is essentially a theorem
of Alper, Hall and Rydh \cite{AHR21}. The present version of Theorem \ref{adapted-atlas} was formulated by Mayeux. The proof of Theorem \ref{adapted-atlas} could be seen as a corollary of \cite{AHR21}. Establishing Theorem \ref{adapted-atlas} as a corollary of \cite{AHR21} is due to Romagny and is the topic of the appendix. Note that Theorem \ref{adapted-atlas} is numbered as Theorem \ref{adapted-atlas-app} in the appendix.

\theo \label{adapted-atlas} Assume that $M$ is finitely generated as abelian group.
Let $X$ be a quasi-separated $S$-algebraic space locally of
finite presentation endowed with an action of $D(M)_S$.
Let $Z \subset M$ be a subgroup and assume that one of the
following assertions holds:
\begin{enumerate}
\item $X$ is separated over $S$,
\item $M/Z$ is torsion-free.
\end{enumerate}
Then there exists a $Z$-FPR atlas $U \to X$, which may be chosen
quasi-compact (in particular $S$-affine here) if $X\to S$ is quasi-compact.
\xtheo

\pf The proof is given in the Appendix of this article, cf. Theorem \ref{adapted-atlas-app}.
\xpf

In fact the author expects that the following stronger assertions should be true. We will study the correct generality in which this holds. We think that this is an interesting problem.

\conj \label{conjecture} Assume that $M$ is finitely generated as abelian group. Let $S$ be an arbitrary scheme.
Let $X$ be a separated $S$-algebraic space locally of
finite presentation endowed with an action of $G=D(M)_S$.  The following assertions hold.
\begin{enumerate}
\item 
There exists a strongly-FPR atlas $U \to X$.
\item 
If $X/S$ is quasi-compact, there exists an $S$-affine strongly-FPR atlas $U \to X$.
\end{enumerate}
\xconj

Note that §\ref{Zariskicase} already ensures the existence of strongly-FPR atlases for Sumihiro actions.

\section{Zariski locally linearizable actions} \label{Zariskicase} 
In this section, we focus on a specific class of actions.

We proceed with the notation from §\ref{FPR}. We recall the definition of Zariski locally linearizable actions and explain that they have strongly-FPR atlases. We recall that some actions are far from being locally linearizable;  this gives a reason to use étale atlases (cf. § \ref{FPR}) and the theory of algebraic spaces.
\fact\label{fpropen}
Let $U \to X$ be a $D(M)_S$-equivariant open immersion, then for any subgroup $Z \subset M$, \[U^{D(M/Z)_S} \to X^{D(M/Z)_S} \times _X U \] is an isomorphism. In other words,  every $D(M)_S$-equivariant open immersion is $Z$-FPR for all subgroups.
\xfact 

\pf This works for any equivariant monomorphism.
For any scheme $T$ over $S$,  $U(T) \times _{X(T)} X^{D(M/Z)_S}(T)= U^{D(M/Z)_S}(T)$.
\xpf 

\defi The action of $D(M)_S$ on $X$ is Zariski locally linearizable if there are $S$-affine open $D(M)_S$-stable subspaces of $X$ covering $X$. We also refer to these actions as Sumihiro actions or Sumihiro spaces.
\xdefi

\fact
Let $X$ be an $S$-space endowed with a Zariski locally linearizable action of $D(M)_S$, then there exists a strongly-FPR equivariant atlas.
\xfact
\pf Let $U$ be the disjoint union of $S$-affine open $D(M)_S$-stable subspaces of $X$ covering $X$, then $U \to X$ is as required by Fact \ref{fpropen}.
\xpf 
\rema\label{remasum} Assume that $M$ is finitely generated. The result \cite[Corollary 3.11]{Su75} shows that if $D(M)_S$ is smooth over $S$ and if $X$ is a scheme and satisfies the condition (N) (e.g. $X$ is normal, cf. \cite[3.4, 3.5]{Su75}), then the action of $D(M)_S$ on $X$ is Zariski locally linearizable. 
\xrema

\rema
There exists a simple example of an action of $D(\bbZ)_S$ on a quasi-separated
 locally finitely presented scheme $X/S$ that is not Zariski locally linearizable (e.g. cf. \cite[§0.2]{Ri16}).
\xrema

 \section{Representability and properties} \label{repandprop}
 
 This section focuses on the representability of attractors. 
 
Let $M$ be an abelian group. Let $S$ be an arbitrary scheme.
Let $X$ be an $S$-algebraic space endowed with an action of $D(M)_S$. We prove that attractors are representable once we know the existence of FPR atlases. Recall that Theorem \ref{adapted-atlas} or Section \ref{Zariskicase} provide FPR atlases in great generality, so that our representability results hold in great generality.

\prop \label{repgroup} 
Let $Z \subset M$ be a
subgroup. Assume that there exists a $Z$-FPR atlas for $X$ under the action of $D(M)_S$ (cf. Definition \ref{deffpratl}). Then the attractor space $X^Z$ (which identifies with the fixed space $X^{D(M/Z)_S}$ by Proposition \ref{Group}) is representable by a closed subspace of $X$.
\xprop

\pf Let $U $ be a $Z$-FPR atlas. Since $U$ is $Z$-FPR, by Definition \ref{deffpratl}, we have a cartesian square \[\xymatrix{U^Z \ar[r] \ar[d] & X^Z \ar[d] \\ U \ar[r] & X.}\] By Proposition \ref{inclusionclosed} the left vertical arrow is a closed immersion of schemes. The lower arrow is surjective étale. So by \cite[\href{https://stacks.math.columbia.edu/tag/03I2}{Tag 03I2}]{stacks-project}, $X^Z$ is an algebraic space and the right vertical arrow is a closed immersion. We used that being a closed immersion
is stable under base change, fppf-local on the base and closed immersions satisfy fppf-descent (cf. \cite[\href{https://stacks.math.columbia.edu/tag/03M4}{Tag 03M4}]{stacks-project} and its proof).
\xpf

\rema \label{rem2} Theorem \ref{G_connected_implies_fixed_pts_closed} of the appendix provides a result similar to Proposition \ref{repgroup} about closedness.
\xrema

\theo \label{representable} Let $N\subset M$ be a magnet, i.e. an arbitrary submonoid. Let $N^* \subset N$ be the face of invertible elements of $N$. Assume that there exists a $N^*$-FPR atlas for $X$ under the action of $D(M)_S$ (cf. Definition \ref{deffpratl}). Assume that one of the following conditions hold
\begin{enumerate}
\item[(a)] $X$ is separated over $S$, or
\item[(b)] $N$ is finitely generated as monoid,
\end{enumerate} then the attractor $X^N$ is representable by an
algebraic space over $S$. Moreover 
 \begin{enumerate}

\item  The morphism of algebraic spaces $X^N \to X^{N^*}$ is affine.

\item If $X/S$ is quasi-separated, then $X^N/S$ is quasi-separated.
\item If $X/S$ is separated, then $X^N/S$ is separated.

\item If $X/S$ is  locally of finite presentation, then $X^N/S$ is  locally of finite presentation.
\item If $X/S$ is quasi-compact then $X^N/S$ is quasi-compact.

\item If $X/S$ if of finite presentation, then $X^N/S$ is of finite presentation.

 \item If $X$ is a scheme then $X^N$ is a scheme.
\end{enumerate}
\xtheo

\pf
Let $U\to X$ be a $N^*$-FPR atlas. Using Definition \ref{deffpratl} and Proposition
\ref{fixxx}, we obtain a diagram with
cartesian squares:
\[
\begin{tikzcd}
U^N \ar[r] \ar[d] & U^{N^*} \ar[r] \ar[d] & U \ar[d] \\
X^N \ar[r] & X^{N^*} \ar[r] & X.
\end{tikzcd}
\]
 By \cite[\href{https://stacks.math.columbia.edu/tag/0466}{Tag 0466}]{stacks-project} and \cite[\href{https://stacks.math.columbia.edu/tag/03MH}{Tag 03MH}]{stacks-project} the vertical maps are \'etale and surjective,
and we know from Theorem \ref{representableaffine} that $U^N$ is
representable by a disjoint sum of $S$-affine schemes. By Proposition \ref{repgroup}, $X^{N^*}$ is an $S$-algebraic space. It follows
from \cite[\href{https://stacks.math.columbia.edu/tag/03I2}{Tag 03I2}]{stacks-project} that $X^N$ is representable.
Now let us prove the listed properties. 
 We have a cartesian square \[ \xymatrix{U^{N} \ar[r] \ar[d] & X^{N} \ar[d] \\ U^{N^{*}} \ar[r]& X^{N^{*}}.}\] As explained before, the lower horizontal arrow is étale and surjective. The left vertical arrow is affine because $U$ is $S$-affine  and $U^N $ and $ U^{N^*}$ are $S$-affine schemes by Theorem \ref{representableaffine}. Now, we apply \cite[\href{https://stacks.math.columbia.edu/tag/03I2}{Tag 03I2}]{stacks-project} to conclude that the right vertical arrow is affine. We used that being affine is preserved under base change, fppf local on the base and satisfies descent for fppf coverings  (\cite[\href{https://stacks.math.columbia.edu/tag/041Z}{Tag 041Z}]{stacks-project} and  \cite[\href{https://stacks.math.columbia.edu/tag/0DFQ}{Tag 0DFQ}]{stacks-project}). Assume that $X$ is locally of finite presentation, i.e. commutes with colimits of affine schemes in $({Sch}/S)_{\mathrm{fppf}}$ (cf.  \cite[\href{https://stacks.math.columbia.edu/tag/049I}{Tag 049I}]{stacks-project}), then the functor $X^N$ commutes with colimits of affine schemes in $({Sch}/S)_{\mathrm{fppf}}$, i.e. $X^N$ is locally of finite presentation.
Assume that $X\to S$ is separated. By Proposition \ref{repgroup}, $X^{N^*} \to X$ is a closed immersion. So by \cite[\href{https://stacks.math.columbia.edu/tag/04CH}{Tag 04CH}]{stacks-project}, $X^{N^*} \to X$ is separated. So by \cite[\href{https://stacks.math.columbia.edu/tag/03KQ}{Tag 03KQ}]{stacks-project} $X^{N^*} \to S$ is separated. Now $X^N \to X^{N^*}$ is separated because it is affine. So $X^N/S$ is separated. 
Assume that $X/S$ is quasi-compact. Note that, by \cite[\href{https://stacks.math.columbia.edu/tag/04CG}{Tag 04CG}]{stacks-project}, $X^{N^*} \to X$ is quasi-compact because it is closed in $X$ (cf. Proposition \ref{repgroup}). Since $X^N \to X^{N^*}$ is affine, it is also quasi-compact. It follows that $X^N/S$ is quasi-compact.
Assume that $X$ is a scheme. Then $X^{N^*}$ is a scheme because it is closed in $X$ (cf. Proposition \ref{repgroup}). Since $X^N \to X^{N^*}$ is affine and $X^{N^*} $ is a scheme, we obtain than $X^N$ is a scheme (\cite[\href{https://stacks.math.columbia.edu/tag/03WG}{Tag 03WG}]{stacks-project}).
 Note that by definition $X$ is of finite presentation if it is locally of finite presentation, quasi-compact and quasi-separated. So if $X/S$ if of finite presentation, then $X^N/S$ is of finite presentation.
\xpf 

\rema \label{remamonooorep}Note that
Proposition \ref{propequiv} (iii) and Theorem \ref{representable} together show representability of many attractors under actions of diagonalizable monoid schemes outside affine cases.
\xrema
\prop \label{existfini} Assume that $M$ is finitely generated as abelian group.
Assume that $X\to S$ is quasi-compact and separated. Then there exists a finitely generated submonoid $N_c $ of $ N$ such that $X^{N}=X^{N_c}$ and $N^* = N_{c}^*$. 
\xprop 

\pf
 Let $\kappa :  \bbN \to N $ be a bijection. For any $i \in \bbN$, let $N_i$ be the submonoid of $N$ generated by $\{ N^*  \cup \kappa ( \bbN _{\leq i} )\}$. Then $N_i ^* = N^* $ and $N_i$ is finitely generated for any $i \in \bbN$, moreover $\bigcup _{i \geq 0} N_i = N$. By Proposition \ref{colimcolim} $
X^N (T) = \colim X^{N_i} (T) . $
Let $U \to X$ be an $S$-affine $N^*$-FPR atlas as in Theorem \ref{adapted-atlas}. Using that $U$ is an $S$-affine scheme of finite presentation over $S$, there exists an integer $c$ such that $U ^{N_c} = U^{N_i}$ for all $i \geq c$. To see this, write $U$ as the spectrum of a graded $\calO _S$-algebra, consider degrees of a finite set of generators of this $\calO _S$-algebra and use that gradings of $\calO _S$-algebras are preserved under equivariant morphisms of $S$-affine schemes.
Let $i \geq c$, and consider the diagram 
\[
\xymatrix{U^{N_{c}} \ar[r] \ar[d] & X^{N_{c}} \ar[d] \\ U^{N_i} \ar[r]& X^{N_i}}.\]
The horizontal arrows are étale and surjective using Proposition \ref{fixxx}. Moreover, the left vertical arrow is an isomorphism since $i \geq c$. Furthermore, the diagram is a cartesian square because by Proposition \ref{fixxx}
\begin{align*}
U^{N_i } \times _{X^{N_i}} X^{N_c} & \simeq (U^{N^*} \times _{X^{N^*}}X^{N_i} ) \times _{X^{N_i}} X^{N_c} \\
& \simeq U^{N^*} \times _{X^{N^*}} X^{N_c} \\
& \simeq U^{N_c}.
\end{align*}
Using  \cite[\href{https://stacks.math.columbia.edu/tag/03I2}{Tag 03I2}]{stacks-project} for the property ``isomorphism", we deduce that for all $i \geq c$, $X^{N_i} = X^{N_c}.$
  This implies that $X^N = X^{N_c}$. 
\xpf

\section{{Hochschild cohomology for diagonalizable monoids via formulas}} \label{hoschsection}

In this section, we develop Hochschild cohomology for actions of monoid schemes. We then specialize to the case of diagonalizable monoid actions. This material will be used in other sections to prove smoothness properties of attractors.

 Let $S$ be a scheme and $B$ be an algebraic space over $S$. 
 Let $B_{spaces , fppf}$ be the small fppf site of $B$ (defined as in the case of schemes, e.g. cf. \cite{Sc17}). In this section $\calO _B$ is the canonical sheaf of rings on $B_{spaces, fppf}$.

\subsection{Affine morphisms and sheaves} 

\prop   Let $f:X \to B$ be an affine morphism of algebraic spaces. Let $\calA $ be the quasi-coherent $\calO_{B}$-algebra of $X$. Let $\calF $ be a quasi-coherent $\calO _{B}$-module. Then 
\[ f_*f^* \calF = \calF \otimes _{\calO _{B}} \calA . \]
\xprop
\pf Let $f^* f_*\calO _X \to \calO _X$ be the canonical morphism coming from adjunction.  It induces a canonical morphism \[f^* \calF \otimes _{\calO _X} f^{*} f_*{\calO _X} \to f^* \calF.\] Since $\calA = f_*\calO _X$, we get obtain a canonical morphism $f^* \calF \otimes _{\calO _X} f^{*} \calA \to f^* \calF$. Now using \cite[\href{https://stacks.math.columbia.edu/tag/03EL}{Tag 03EL}]{stacks-project}, we obtain a canonical morphism \[f^*( \calF \otimes _{\calO _B}  \calA  )\to  f^* \calF.\]By adjunction, we finally get a canonical morphism $\calF \otimes _{\calO _B}  \calA  \to f_* f^* \calF$. It remains to prove that it is an isomorphism. For this, we reduce to the affine case and apply  \cite[\href{https://stacks.math.columbia.edu/tag/01I8}{Tag 01I8}]{stacks-project}.
\xpf 

\subsection{Fppf Hochschild monoid cohomology over algebraic spaces} \label{Hoschgene}

Let $H$ be a monoid algebraic space over $B$. Let $\calO _B$ be the sheaf of rings of $B$ on $B_{spaces, fppf}$.

\defi \label{GOmodetaledef} A $H$-$\calO _B$-module over $B_{spaces, fppf}$ is an $\calO _B$-module $\calF$ in the sense of  \cite[\href{https://stacks.math.columbia.edu/tag/03CW}{Tag 03CW}]{stacks-project} such that $H$ acts $\calO_B$-linearly on $\calF$, i.e:
\begin{enumerate}
\item For any $T \to B \in B _{spaces, fppf}$, $H(T)$ acts on $\calF (T)$. Moreover for any $g \in H(T), x ,y\in \calF(T)$ and $a \in \calO _B (T)$, we have $g \cdot (x+ay) = g \cdot x + a  (g \cdot y) $.
\item If $T' \to T $ is a morphism in $B_{spaces, fppf}$, the action of $H(T')$ on $\mathcal{F} (T')$ is compatible with the action of $H(T)$ on $\calF (T)$.
\end{enumerate}
\xdefi

\rema
Note that Definition \ref{GOmodetaledef} is similar to Definition \ref{GOmodZarsch}, but is different.
\xrema

\defi \label{partialdef}
 Let  $\calF$ be a  $H$-$\calO_B$-module over $B_{spaces, fppf}$. 
 Let $n \geq 0$ be an integer. We put  \[C^n (H , \calF) = Mor _{B_{spaces , fppf}} (H^n , \calF) =\text{The set of natural transformations from }H^n \text{ to } \calF, \]
 where $H^{n}$ and $\calF $ are seen as functors from $B_{spaces, fppf}$ to $\Sets$. The set $C^n(H,\calF )$ is canonically endowed with a structure of $\calO_B (B)$-module via 
 \[\text { for any }\Theta, \Psi \in C^n(H, \calF), a \in \calO _B (B), T \in B_{spaces, fppf}; (\Theta + a \Psi) (T) = \Theta (T) + {a\mid _T} \Psi (T) .\]
 
   We have a morphism of $\calO_B (B)$-modules \[{\partial : C^n(H,\calF) \to C^{n+1} (H , \calF)}\] sending $\Theta $ to $\partial \Theta $ where $\partial \Theta $ is the transformation from $H^{n+1} $ to $\calF$ such that for any $T \in B_{spaces, fppf}$, $\partial \Theta (T)$ is the map sending $g_1 , \ldots , g_{n+1} \in H^{n+1} (T) $ to 
 \[ g_1 \cdot \Big( \Theta (T) \big( g_2 , \ldots , g_{n+1} \big)\Big) + \sum_{i=1}^n (-1)^i \Theta (T) \big( g_1 , \ldots , g_i  g_{i+1} , \ldots , g_{n+1} )+ (-1)^{n+1} f (g_1 , \ldots , g_n \big) \in \calF (T). \] 
 \xdefi 
 We also have an internal version of $C^n(H, \calF)$ as follows.
\defi \label{partialdef}
 Let  $\calF$ be a  $H$-$\calO_B$-module. 
 Let $n \geq 0$ be an integer. We put for any $V \in B_{spaces, fppf}$ \[\underline{C^n (H , \calF)} (V)  = C^n ( H|_V , \calF |_V) = Mor _{V_{spaces, fppf}} (H^n |_V , \calF |_V ) . \] The functor $\underline{C^n(H,\calF )}$ is canonically endowed with a structure of $\calO_B $-module.
   We have a morphism of $\calO_B $-modules ${\underline{\partial} : \underline{C^n(H,\calF) }\to \underline{C^{n+1} (H , \calF)}}$ induced by $\partial$.
 \xdefi 
 
 \rema For any positive integer $n$ the composition $ \underline{\partial} \circ \underline{\partial} $ is the zero map.
 \xrema 

\defi Let $\calF$ be a $H$-$\calO_B$-module  over $B_{spaces, fppf}$.  We put 
\begin{align*}
Z^{n}_{monoid} (H , \calF ) &= \ker (C^n (H , \calF ) \xrightarrow{\partial} C^{n+1} (H , \calF))  \text{ for any } n \geq0\\
B^0 _{monoid}(H, \calF ) &= 0\\
B^n _{monoid} (H , \cal F) &= \mathrm{im} ( C^{n-1} (H , \calF)  \xrightarrow{\partial} C^n (H , \calF )\text{ for any } n \geq 1\\
H^n _{monoid} (H , \calF) &= \frac{Z^n (H , \calF )}{B^n (H , \calF )} \text{ for any } n \geq 0.
\end{align*}

\begin{align*}
\underline{Z^{n}_{monoid} (H , \calF )} &= \ker (\underline{C^n (H , \calF )} \xrightarrow{\underline{\partial}} \underline{C^{n+1} (H , \calF)})  \text{ for any } n \geq0\\
\underline{B^0 _{monoid}(H, \calF ) }&= 0\\
\underline{B^n _{monoid} (H , \cal F) }&= \underline{\mathrm{im} ( C^{n-1} (H , \calF)}  \xrightarrow{\underline{\partial}} \underline{C^n (H , \calF )}\text{ for any } n \geq 1\\
\underline{H^n _{monoid} (H , \calF)} &= \frac{~\underline{Z^n (H , \calF )}~}{~\underline{B^n (H , \calF )}~} \text{ for any } n \geq 0.
\end{align*}\xdefi
 We finish this subsection with the description of $Z^n_{monoid}$ and $B^n_{monoid}$ for $n=1$.
 \rema \label{remzbun} We have \[Z^1_{monoid}(H , \calF)= \\ \{ c\in C^1(H , \calF ) \mid \forall T \in B_{spaces , fppf}, \forall g,g' \in H(T),~ c(gg') = c(g) + g \cdot c(g') \}.\]
 \[ B^1 _{monoid}(H , \calF ) = \{ b \in C^1 (H , \calF ) \mid \exists  v \in \calF (B); \forall T \in B_{sp., fppf} , \forall g \in G(T), b(g) = g\cdot v_T - v_T \}.\]
 \xrema

 \subsection{Modules under affine monoid algebraic spaces}

  Let $M$ be a finitely generated cancellative monoid. The structural morphism $A(M)_B \to B $ and the multiplication morphism $A(M)_B \times _B A(M)_B \to A(M)_B $ are surjective and  finitely presented by Facts \ref{fact-RNAN} and \ref{fifififif}. The morphism $A(M)_B \to B$ is flat. We assume in this section that the multiplication morphism $A(M)_B \times _B A(M)_B \to A(M)_B $ is flat. We refer to Proposition \ref{flatflatflat} for a criterion. For example if $M$ is a group, the multiplication morphism $D(M)_B \times _B D(M)_B \to D(M)_B $ is flat.

 \prop \label{homprecomod}Let $\calF , \calG$ be two $\calO _B$-modules (over $B_{spaces, fppf}$). Assume $\calG$ is quasi-coherent. Let 
 \[ X '= \Spec _B (\calA ' ) \xrightarrow{m} X = \Spec _B (\calA ) \xrightarrow{f} B .\] be affine morphisms of $S$-spaces. Assume $m$ and $f$ are flat and finitely presented. Then we have a canonical identification 
 \[ \calhom _{\calO _B} ( \calF , \calG ) (X) = \Hom _{\calO _B} (\calF , \calG \otimes _{\calO _B} \calA ). \] Moreover the map 
 \[ \calhom _{\calO _B} (\calF , \calG) (X) \to \calhom _{\calO _B} (\calF , \calG ) (X')  \] associated to $m$ corresponds to the map 
 \begin{align*}
 \Hom _{\calO _B} (\calF , \calG \otimes _{\calO _B} \calA ) &\to \Hom _{\calO _B} ( \calF , \calG \otimes _{\calO _B} \calA' ) \\
 \varphi \mapsto ( \calF \xrightarrow{\varphi} \calG &\otimes _{\calO _B} \calA \xrightarrow{\Id \otimes \phi_m } \calG \otimes _{\calO _B} \calA ' )
 \end{align*} where $\phi _m : \calA \to \calA'$ is dual to $m$ via \cite[\href{https://stacks.math.columbia.edu/tag/081V}{Tag 081V}]{stacks-project}.
 \xprop
 
 \pf Note that $m $ and $f$ are flat and finitely presented. The following identifications prove the first assertion\begin{align*}
  \calhom_  {\calO _B} (\calF , \calG) (X ) 
  & = \Hom _{\calO _{X}}(f^* \calF ,f^* \calG )  \\
  & = \Hom _{\calO _B}( \calF ,f_*f^* \calG )  \\
  & = \Hom _{\calO _B}( \calF ,\calG \otimes _{\calO _B} \calA ) . 
 \end{align*}
 The second assertion also follows from these intermediate identifications.
 \xpf

 \prop \label{comodulefppf} Let $H= \Spec _B (\calA)$ be a monoid algebraic space over $B$ such that $H \to B$ is affine, flat and finitely presented. We assume that the multiplication $H \times _B H \to H$ is flat. An $H$-$\calO _{B}$-module structure (over $B_{spaces, fppf}$) on a quasi-coherent $\calO_{B}$-module $\calF$ corresponds to an $\calO _{B}$-comodule structure on $\calF$, i.e a morphism of $\calO _{B }$-modules $\calF \xrightarrow{\delta} \calF \otimes _{\calO _{B} } \calA$ such that $(\delta \otimes \Id _\calA )\circ \delta = ( \Id_{\calF} \otimes \Delta _{\calA} ) \circ \delta $ and $ (\Id_ {\calF} \otimes \varepsilon  _\calA) \circ \delta = \Id _{\calF} $.
 \xprop
 
 \pf
 Let us first remark that the set of $H$-$\calO _{B}$-module structures on an $\calO_{B}$-module $\calF$ is a subset of the set $Mor_{B_{spaces , fppf } }(H, \calhom_  {\calO _{B,fppf}} (\calF , \calF) )$. Using Proposition \ref{homprecomod} we get 
 \begin{align*}
 Mor_{B_{spaces , fppf} } (H, \calhom_  {\calO _B} (\calF , \calF) )&=  \calhom_  {\calO _B} (\calF , \calF) (H ) \\
  & = \Hom _{\calO _B}( \calF ,\calF \otimes _{\calO _B} \calA ) . 
 \end{align*}
 The condition that $\calF$ is an $H$-module means that the transformations of functors 
 \[ T_1: H  \times H \xrightarrow{m} H \xrightarrow{\Psi} \calhom _{\calO _B} (\calF , \calF) \] and 
 \[ T_2 :H \times H \xrightarrow{\Psi \times \Psi } \calhom _{\calO _B} (\calF , \calF) \times \calhom _{\calO _B} (\calF , \calF) \xrightarrow{\circ} \calhom _{\calO _B} (\calF , \calF)   \]
 are equal and that the composition 
 \[ e_B \xrightarrow{} H \xrightarrow{\Psi} \calhom _{\calO _B} (\calF , \calF) \] is the identity.
 The transformations $T_1 $ and $T_2$ correspond to elements in \[\Hom _{\calO _B} (\calF , \calF \otimes _{{\calO} _B} \calA \otimes _{\calO _B} \calA ),\] namely $T_1$ corresponds to $( \Id_{\calF} \otimes \Delta _{\calA} ) \circ \delta $ (cf. Proposition \ref{homprecomod}) and $T_2 $ corresponds to $(\delta \otimes \Id _{\calA} )\circ \delta $. The composition 
 \[ e_B \xrightarrow{} H \xrightarrow{\Psi} \calhom _{\calO _B} (\calF , \calF) \] corresponds to an element in $\Hom _{\calO _B} ( \calF , \calF \otimes _{\calO _B} \calO _B) $, namely $(\Id_ {\calF} \otimes \varepsilon _\calA) \circ \delta$. This finishes the proof. 
 \xpf 
 
 \prop\label{amstrumod} An $A(M)_B$-$\calO _{B}$-module (over $B_{spaces, fppf}$) structure on a quasi-coherent $\calO_{B}$-module $\calF$ corresponds to a collection $(\mu _m) _{m \in M} \in \Hom _{\calO _B } (\calF , \calF)$ such that only finitely many of the $(\mu _m) _{m \in M}$ are nonzero locally and for all $m,k \in M$,
 \[ \begin{cases} \mu _m \circ \mu _{k} = \mu _k \text{ if }k=m, \\  \mu_m \circ \mu_k =0 \text{ if } k\ne m \end{cases} \] $\text{ and }  \sum _{m \in M } \mu _m = \Id _{\calF}.$
 \xprop 
 \pf Proposition \ref{comodulefppf} shows that an $A(M)_B$-$\calO _{B}$-module corresponds to a morphism of $\calO _B$-comodule 
 \[  \mu: \calF \to \calF \otimes _{\calO _B} \calO _B [M] = \bigoplus _{m \in M } \calF .\]
Now $\mu$ gives a collection $(\mu _m) _{m \in M} $ with $\mu _m \in \Hom _{\calO _B} (\calF , \calF)$, such that only finitely many of the $(\mu _m) _{m \in M}$ are nonzero locally. Recall from Proposition \ref{comodulefppf} that $\mu$ is a comodule means that $(\mu \otimes \Id _{\calO _B [M]} )\circ \mu = ( \Id_{\calF} \otimes \Delta _{\calO _B [M]} ) \circ \mu $ and $ (\Id_ {\calF} \otimes \varepsilon  _{\calO _B [M]}) \circ \mu = \Id _{\calF} $.
The coassociativity condition 
\[
(\mu \otimes \Id _{\calO _B [M]} )\circ \mu = ( \Id_{\calF} \otimes \Delta _{\calO _B [M]} ) \circ \mu
\] 
reads as follows: the left-hand side is 
\[
(\mu \otimes \Id _{\calO _B [M]} )\circ \mu = \sum_{m,k \in M} (\mu_m \circ \mu_k) \otimes X^m \otimes X^k,
\] 
and the right-hand side is 
\[
( \Id_{\calF} \otimes \Delta _{\calO _B [M]} ) \circ \mu= \sum_{m \in M} \mu_m \otimes X^m \otimes X^m,
\] 
because the comultiplication 
\(\Delta_{\mathcal{O}_B[M]}: \mathcal{O}_B[M] \to \mathcal{O}_B[M] \otimes_{\mathcal{O}_B} \mathcal{O}_B[M]\) 
is defined by \(\Delta(X^m) = X^m \otimes X^m\). Comparing components in the direct sum 
\(\bigoplus_{m,k \in M} \mathcal{F} \otimes X^m \otimes X^k\) 
shows that equality holds if and only if 
\(\mu_m \circ \mu_k = 0\) for \(m \neq k\) and \(\mu_m \circ \mu_k = \mu_k\) if $k=m$. 
The counit condition 
\[
(\mathrm{Id}_{\mathcal{F}} \otimes \varepsilon_{\mathcal{O}_B[M]}) \circ \mu = \mathrm{Id}_{\mathcal{F}}
\]
reads as follows. Recall that the counit is defined by
\(\varepsilon_{\mathcal{O}_B[M]}(X^m) = 1\) for all \(m \in M\). Applying it, we obtain
\[
(\mathrm{Id}_{\mathcal{F}} \otimes \varepsilon_{\mathcal{O}_B[M]}) \circ \mu 
= \sum_{m \in M} \mu_m \cdot \varepsilon_{\mathcal{O}_B[M]}(X^m) 
= \sum_{m \in M} \mu_m.
\]
Thus, the counit axiom exactly gives
\[
\sum_{m \in M} \mu_m = \mathrm{Id}_{\mathcal{F}}.
\]

Conversely, any such collection \((\mu_m)_{m \in M}\) defines a comodule map 
\(\mu = \sum_{m \in M} \mu_m \otimes X^m\)
which satisfies both the coassociativity and the counit axioms. \xpf

\prop \label{hoschaff} Let $\calF$ be a quasi-coherent $A(M)_B$-$\calO_B$-module over $B_{spaces, fppf}$. Let $V \in B_{spaces , fppf}$. We have \[\underline{C^n(A(M)_B, \calF)}(V) = \Big(\bigoplus_{(m_1 , \ldots , m_n )\in M^n} \calF  X^{m_1}   \cdots  X^{m_n}\Big) (V) \]
where for every $m=(m_1 , \ldots , m_n) \in M^n$, $\calF   X^{m_1}   \cdots  X^{m_n} = \calF ^m$ is a formal copy of $\calF$. In particular if $V$ is quasi-compact we have 
\[\underline{C^n(A(M)_B, \calF)}(V) = \bigoplus_{m\in M^n}  \calF ^m   (V) \]

 Moreover for any $n \geq 0$, the linear map  $\partial : C^n \to C^{n+1} $ is given on the component $\calF X^{m_1 } \cdots X^{m_n}$ and on quasi-compact objects by the formula
 
 \begin{align*} f^{m_1 , \ldots, m_n}  X^{m_1}   \cdots  X^{m_n} \mapsto &\mu (f^{m_1 , \ldots, m_n}) X^{m_1}  \cdots X^{m_n} \\
 +  \sum _{i=1}^n (-1)^i  f^{m_1 , \ldots, m_n} & X^{m_1} \cdots  \Delta (X^{m_i} )\cdots  X^{m_n} \\
+ (-1)^{n+1}f^{m_1,\cdots, m_n}&  X^{m_1} \ldots  X^{m_n}  X^0
\end{align*} where $\Delta (X^{m_i}) = X^{m_i}X^{m_i}$.
\xprop

\pf Put $p : A(M)_V^n \to V$ and $\calF _V = \calF |_V $. The identifications \begin{align*} \underline{C^n (A(M)_B, \calF)}(V)&= Mor_{V_{spaces, fppf}} (A(M)_V^n, \calF_V) \\ 
[\calF_V= \calhom _{\calO _V} (\calO _V, \calF_V) ]~~~~&= Mor_{V_{spaces, fppf}} (A(M)_V^n, \calhom _{\calO _V} (\calO _V, \calF_V)) \\ 
\cite[\href{https://stacks.math.columbia.edu/tag/03AB}{Tag 03AB}]{stacks-project}~~~~ &=  \calhom _{\calO _V} (\calO _V, \calF_V) ( A(M)_V^n )\\
\cite[\href{https://stacks.math.columbia.edu/tag/04TT}{Tag 04TT}]{stacks-project}  \cite[\href{https://stacks.math.columbia.edu/tag/03CZ}{Tag 03CZ}]{stacks-project} ~~~~ &= \Hom _{\calO _V} (\calO _V , \calF_V \otimes _{\calO _V}  \calO_V [M] \otimes _{\calO _V} \cdots \otimes _{\calO _V} \calO _V [M])\\
 &= \Big( \calF _V \otimes _{\calO _V}   \calO_V [M] \otimes _{\calO _V} \cdots \otimes _{\calO _V} \calO _V [M] \Big) (V)\\
  &= \Big( \calF \otimes _{\calO _B}   \calO_B [M] \otimes _{\calO _B} \cdots \otimes _{\calO _B} \calO _B [M] \Big) (V)\\
&= \Big( \bigoplus_{(m_1 , \ldots , m_n )\in M^n} \calF X^{m_1}  \cdots X^{m_n} \Big) (V)
\end{align*}
 prove the first assertion. The second assertion is now a consequence of Definition \ref{partialdef} and the definition of $\partial$.
\xpf

\prop \label{hunzerodiag} Let $\calF$ be a quasi-coherent $A(M)_B$-$\calO_B$-module over $B_{spaces, fppf}$, then \[H^1_{monoid}(A(M)_B , \calF ) =0.\]
\xprop 

\pf \begin{sloppypar} It is enough to prove that $Z^1 _{monoid}(A(M)_B , \calF) \subset B^1 (A(M)_B , \calF)$. Let $U $ be a quasi-compact object in $B_{spaces, fppf}$. Let us use the identification given in Proposition \ref{hoschaff} and \cite[\href{https://stacks.math.columbia.edu/tag/0935}{Tag 0935}]{stacks-project}.
 Put $F = \Gamma (U , \calF) $ and $R = \Gamma ( U , \calO _B)$.
 Let $\xi_{U}= \sum _{m \in M} f_m  X^m $ be an element in $  \Gamma ( U , \calF \otimes \calO _S [M] )= \bigoplus _{m \in M } F X^m$. Assume that $\partial ( \xi_{U} ) =0$,  then 
\[
0 = \Big( \sum _{(k,m) \in M^2}  \mu _k (f_m)  X^{(k,m)} \Big)
- \Big(\sum _{l \in M }  f_l X^{(l,l)} \Big) +\Big(\sum _{n \in M } f_n X^{(n,0)} \Big).
\]
 So $ \mu _k (f_0) = -f_k $ for all $k \in M$ with $k \ne 0$ and
$ \mu_0 (f_0)  = 0 $.
 Put $e_U:= -f_0 \in  F= \Gamma (U , \calF)$. We have 
 \[\partial (e_U) = \Big( \sum _{k\in M } \mu _k (e_U )X^k \Big) -e_U X^0  = \sum _{0 \ne k \in M} f_k  X^k + f_0  X^0 = \xi_{U}. \]
 So $\xi_{U}$ belongs to the image of $\Gamma (U , \calF) \xrightarrow{\partial} \Gamma ( U , \calF \otimes \calO _B [M])$. Now let $\xi \in \Gamma ( B, \calF \otimes \calO_B [M] )$ and assume $\partial (\xi )=0$, i.e. $\xi \in Z^1_{monoid} (A(M)_B , \calF)$. Let us consider the set $\{U_i\}_{i \in I}$ of all quasi-compact objects in $B_{spaces, fppf}$, of course it provides an fppf covering of $B$. For each $i \in I $, $\partial (\xi | _{U_i} )=0$, so the above computation gives us elements $e_i \in \Gamma ( U_i , \calF)$ such that $\partial ( e_i ) = \xi | _{U_i} $. These $e_i$ are compatible under restrictions, since they are given by explicit compatible formula. So there exists an element $e \in \Gamma ( B , \calF )$ such that $e|_{U_i} = e_i $ for all $i \in I$ (\cite[\href{https://stacks.math.columbia.edu/tag/00VM}{Tag 00VM}]{stacks-project}). We have $\partial (e) = \xi $. So $\xi \in B^1 _{monoid}(A(M)_B , \calF)$. This finishes the proof. \end{sloppypar}\xpf 
 
\rema  We conjecture that our method given in the proof of Proposition \ref{hunzerodiag} will lead to a proof that $H^n_{monoid}(A(M)_B , \calF)=0$ also for $n \geq 2$. Note that \cite[Exp I. Théorème 5.3.3]{SGA3} proves that  $H^n(D(M)_S , \calF)=0$ ($n \geq 1$) under the assumption that $S$ is affine and that $\calF$ is quasi-coherent, using derived functors. Our computational method, highly inspired by formulas given --and unexploited-- in \cite{SGA3}, shows that the assumption that $S$ is affine is unnecessary in \cite[Exp. I Théorème 5.3.3]{SGA3} (at least in the case $n=1$). 
\xrema

\section{{Equivariant infinitesimal deformations}}
\label{deformation}

Using Hochschild cohomology, we prove the existence of morphisms filling diagrams required in the formal smoothness conditions. Here we focus on the existence of such morphisms in general settings.

Let $S$ be a scheme. Let $G$ be a group algebraic space over $S$.

\prop \label{existencerelev}\begin{sloppypar} Let $X,X',Y$ and $ Z$ be $S$-algebraic spaces endowed with $S$-actions of $G$. Let \[ \xymatrix{X\ar[r]^{a} \ar[d]^i& Y\ar[d]^{f} \\ X' \ar[r]^{\varphi}&Z  } \] be a commutative diagram of $G$-equivariant $S$-morphisms where the left vertical morphism $i$ is a first order thickening, cf. \cite[\href{https://stacks.math.columbia.edu/tag/05ZK}{Tag 05ZK}]{stacks-project}.
Assume that $Y \to Z$ is formally smooth and \[H^1_{monoid}(G, p_* \calhom _{\calO _X} (a^* \Omega _{Y/Z} , \mathcal{C}_{X/X'} ))=0,\]  where $p_*: Sh(X'_{spaces,fppf}) \to  Sh( S_{spaces, fppf}) $ (cf.  \cite[\href{https://stacks.math.columbia.edu/tag/00X6}{Tag 00X6}]{stacks-project} for $p_*$, § \ref{Hoschgene} for $H^1_{monoid} (-,-)$ and  \cite[\href{https://stacks.math.columbia.edu/tag/04CN}{Tag 04CN}]{stacks-project} for $\mathcal{C}_{X/X'}$). Then there exists a $G$-equivariant morphism $\phi: X' \to Y$ filling the diagram, i.e. such that $\phi \circ i = a$ and $f \circ \phi = \varphi.$\end{sloppypar}
\xprop

\pf \begin{sloppypar}
For an object $U'$ of $(X')_{spaces, fppf}$ with $U= X \times _{X'} U'$, consider morphisms $a': U' \to Y$ such that 

(1) $a'$ is a morphism over $Z$, and 

(2) $a' | _{U} = a |_{U}$.\\
Then the arguments of  \cite[\href{https://stacks.math.columbia.edu/tag/061A}{Tag 061A}]{stacks-project} show that the rule $U' \mapsto \{ a':U' \to Y' $ such that (1) and (2) hold.$\}$ defines a sheaf of sets $\mathcal{P}$  on $(X')_{spaces, fppf}$. Note that the condition (3) in \cite[061A]{stacks-project} is empty because $Y=Y'$ in our situation.  By \cite[\href{https://stacks.math.columbia.edu/tag/061C}{Tag 061C}]{stacks-project}, there is an action of the sheaf of abelian groups $\calhom _{\calO _X} (a^* \Omega _{Y/Z} , \mathcal{C}_{X/X'} )$ 
on the sheaf $\mathcal{P}$. 
Moreover, the action of $\calhom _{\calO _X} (a^* \Omega _{Y/Z} , \mathcal{C}_{X/X'} )$ on $\mathcal{P}$ is simply transitive
 for any object $U'$ of $(X')_{spaces, fppf}$ 
 over which the sheaf $\mathcal{P}$ has a section. Let $p:X'\to S$ be the structural morphism. By \cite[\href{https://stacks.math.columbia.edu/tag/00X6}{Tag 00X6}]{stacks-project}, $p$ induces a morphism of topoi $p_*$ \begin{align*} 
p_*: Sh(X'_{spaces,fppf}) &\to Sh( S_{spaces, fppf})\\
 \mathcal{F}~~~~~~& \mapsto \big( (T \to S)  \mapsto \mathcal{F} (X' \times _S T )\big). 
 \end{align*} 
 The action of $\calhom _{\calO _X} (a^* \Omega _{Y/Z} , \mathcal{C}_{X/X'} ))$ 
on $\mathcal{P}$ induces an action of the sheaf of abelian groups $p_* \calhom _{\calO _X} (a^* \Omega _{Y/Z} , \mathcal{C}_{X/X'} )$ 
on the sheaf $p_*\mathcal{P}$. We want to show that $G$ acts on the sheaf of abelian groups $p_*\calhom _{\calO _X} (a^* \Omega _{Y/Z} , \mathcal{C}_{X/X'} )$ and
on the sheaf $p_* \mathcal{P}$, in a compatible way.  Let $T$ be an algebraic space over $S$ and let $g \in G(T)$. 
An element $a' \in ( p_* \mathcal{P})(T)$ is by definition an $S$-morphism $X' \times _S T \to Y  $ such that $a'$ is an $S$-morphism over $Z$ and $a' \mid _{X \times _S T} = a \mid _{X \times _S T}$. In other words, an element $a' \in (p_* \mathcal{P})(T) $ is by definition a $T$-morphism $X' \times _S T \to Y \times _S T $
 such that $a'$ is a $T$-morphism over $Z \times _S T $ and $a' \mid _{X \times _S T} = a \mid _{X \times _S T}$. Let $g^{-1}$ be the inverse of $g$ in $G(T)$. As usual, $X'_T, X_T,Y_T $ and $Z_T$ denote $X' \times _S T , X\times _S T , Y \times _S T $ and $Z \times _S T $. The element $g^{-1}$ gives us elements $g^{-1}_{X'_ T } \in \Aut _{T} ( X' _ T)$ and $ g^{-1}_{X _ T } \in \Aut _{T} ( X _ T)$. The element $g^{}$ gives us an element $g^{}_{Y _ T} \in \Aut _T ( Y _ T)$ and an element $g^{}_{Z _T} \in \Aut _T ( Z_ T)$. For any $a' \in (p_* \mathcal{P})(T)$, we put $g\cdot a':= g _{Y _ T} \circ a' \circ g^{-1} _{X' _ T}$, this is a $T$-morphism $X' _ T \to Y _T$. Let us check that $g\cdot a' $ belongs to $ (p_* \mathcal{P})(T)$. Since $a$ is $G$ equivariant we have $a =  g _{Y _T} \circ a \circ g^{-1} _{X _ T}$, this shows that $g\cdot a' \mid _{X _ T} = a \mid _{X _ T}$.  Let us now show that $ g \cdot a' $ is a $T$-morphism over $Z \times _S T$. Let $f_T $ and $\phi _T$ denote the base change from $S$ to $T$ of $f$ and $\phi$. We have  
\begin{align*}
f_T \circ g_{Y _ T}  \circ a' \circ g^{-1}_{X' _ T}& = \\
 \text{because } f_T \text{ is equivariant} &= g_{Z _ T} \circ f_T \circ a' \circ  g^{-1}_{X' _ T} \\
\text{because }a' \in  (p_* \mathcal{P})(T)  &= g_{Z _ T} \circ \varphi_T \circ  g^{-1}_{X' _ T} \\
\text{because } \varphi _T \text{ is equivariant} &= \varphi _T.
\end{align*} This finishes the verification that $g \cdot a' \in (p_* \mathcal{P})(T)$. We now define the action of $G$ on \[p_* \calhom _{\calO _X} (a^* \Omega _{Y/Z} , \mathcal{C}_{X/X'} )\] giving the action of $g \in G(T)$ on $ v \in ( p_* \calhom _{\calO _X} (a^* \Omega _{Y/Z} , \mathcal{C}_{X/X'} ))(T)$. The diagram \[ \xymatrix{ Y_T \ar[d]^{f_T} \ar[r]^{g_{Y_T}}& Y_T \ar[d]^{f_T} \\ Z_T \ar[r]^{g_{Z_T}}  & Z_T}, \] whose horizontal arrows are isomorphisms, is commutative since $f_T$ is equivariant. So we obtain an automorphism $g_{\Omega _T }$ of $\Omega _{ Y_T / Z_T } $, e.g. by  \cite[\href{https://stacks.math.columbia.edu/tag/04CX}{Tag 04CX}]{stacks-project}. Note that $g_{\Omega _T} \circ h_{\Omega _T} = (hg ) _{\Omega _T}$, e.g. by \cite[\href{https://stacks.math.columbia.edu/tag/05Z7}{Tag 05Z7}]{stacks-project}. That induces similar automorphisms of $a_T^* \Omega _{Y_T / Z_T}$. Similarly, using \cite[\href{https://stacks.math.columbia.edu/tag/04CP}{Tag 04CP}]{stacks-project} and \cite[\href{https://stacks.math.columbia.edu/tag/04G2}{Tag 04G2}]{stacks-project}, we get for any $g,h \in G(T)$  automorphisms $g_{\mathcal{C}_T}$ of $\mathcal{C }_{X_T/X'_T}$ such that $g_{\mathcal{C} _T} \circ h_{\mathcal{C}_T} = (hg ) _{\mathcal{C} _T}$. Now let $v \in \Hom _{\calO _X} (a_T^* \Omega _{Y_T/Z_T} , \mathcal{C} _{X_T /X'_T} ) = (p_* \calhom _{\calO _X} (a^* \Omega _{Y/Z} , \mathcal{C}_{X/X'} ))(T)$ and $g \in G (T)$. We put $g\cdot v = g^{-1} _{\mathcal{C} _T} \circ v \circ g _{\Omega _T}$. The element $g \cdot v $ belongs to $(p_* \calhom _{\calO _X} (a^* \Omega _{Y/Z} , \mathcal{C}_{X/X'} ))(T)$. This defines a left action because 
\[ h \cdot ( g \cdot v) = h ^{-1} _{\mathcal{C} _T} \circ g^{-1} _{\mathcal{C} _T} \circ v \circ g _{\Omega _T}  \circ h _{\Omega _T} =  (g _{\mathcal{C} _T} \circ h _{\mathcal{C} _T})^{-1} \circ v \circ g _{\Omega _T}  \circ h _{\Omega _T} =  (hg)^{-1} _{\mathcal{C} _T} \circ v \circ hg _{\Omega _T}  = hg \cdot v.\]
Using  \cite[\href{https://stacks.math.columbia.edu/tag/0618}{Tag 0618}]{stacks-project},
 we obtain that the actions are compatible in the sense that \[g\cdot a' + g \cdot v = g \cdot (a' +v) \text{ for any } g\in G(T) , a' \in (p_* \mathcal{P} ) (T) , v \in  (p_* \calhom _{\calO _X} (a^* \Omega _{Y/Z} , \mathcal{C}_{X/X'} ))(T).\] 
 We now define an element $c$ in $Z^1 _{monoid} (G , (p_* \calhom _{\calO _X} (a^* \Omega _{Y/Z} , \mathcal{C}_{X/X'} ) ) )$.  Let us fix $a' \in  (p_* \mathcal{P} ) (S)$, this is possible because $Y \to Z$ is formally smooth.  For any $T\to S \in S_{spaces , fppf}$, put $c(T): G(T) \to  (p_*  \calhom _{\calO _X} (a^* \Omega _{Y/Z} , \mathcal{C}_{X/X'} ) ) (T)$, $g \mapsto g \cdot a'_T - a'_T $. Then $c$ is indeed a cocycle because for any $T \in S_{spaces ,fppf} $ and $g,g' \in G(T)$ we have 
  \[c(g)+g \cdot c(g') = g \cdot a' _T -a'_T + g \cdot ( g' \cdot a'_T -a' _T) =  gg' \cdot a'_T-a'_T = c(gg'). \]
By assumption, $H^1_{monoid}  (G , (p_* \calhom _{\calO _X} (a^* \Omega _{Y/Z} , \mathcal{C}_{X/X'} ) ) ) =0$, so there is an element \[v \in  (p_*  \calhom _{\calO _X} (a^* \Omega _{Y/Z} , \mathcal{C}_{X/X'} ) )(S)\] such that $c(g)=g\cdot v_T - v_T $ for all $T \in S_{spaces , fppf}$ and all $g \in G(T)$ (cf. Remark \ref{remzbun}).
So for any $T \in S_{spaces , fppf} $ and any $g \in G(T)$, we have 
$ g \cdot a'_T - a'_T = g \cdot v_T - v_T  $, so we have $ g \cdot ( a'_T - v_T ) = a'_T -v_T $. We now put $\phi:= a'-v  \in (p_* \mathcal{P})(S)= \mathcal{P} ( X') $. Then $\phi $ is a morphism from $X' $ to $ Y$ satisfying all the required properties, this finishes the proof. \end{sloppypar}
\xpf 
\rema The structure of the proof of Proposition \ref{existencerelev} is partly similar to the argument sketched in \cite[Exp. XII proof of Lemma 9.4]{SGA3} (unpublished), though frameworks are different. 
\xrema

\coro 
 \label{existencerelevcor}\begin{sloppypar}Let $M$ be a finitely generated abelian group.  Let $X,X',Y$ and $ Z$ be $S$-algebraic spaces endowed with $S$-actions of $D(M)_S$. Let \[ \xymatrix{X\ar[r]^{a} \ar[d]^i& Y\ar[d]^{f} \\ X' \ar[r]^{\varphi}&Z  } \] be a commutative diagram of $D(M)_S$-equivariant $S$-morphisms where the left vertical morphism $i$ is a first order thickening, cf. \cite[\href{https://stacks.math.columbia.edu/tag/05ZK}{Tag 05ZK}]{stacks-project}.
Assume that $Y \to Z$ is formally smooth and locally of finite type. Assume that $X' \to S$ is quasi-compact and quasi-separated. Then there exists a $D(M)_S$-equivariant morphism $\phi: X' \to Y$ filling the diagram, i.e. such that $\phi \circ i = a$ and $f \circ \phi = \varphi.$\end{sloppypar}
\xcoro
\begin{proof} \begin{sloppypar} We first prove that  $p_* \calhom _{\calO _X} (a^* \Omega _{Y/Z} , \mathcal{C}_{X/X'} )$  (where $p_*: Sh(X'_{spaces,fppf}) \to  Sh( S_{spaces, fppf}) $ is as in Proposition \ref{existencerelev}) is quasi-coherent. By  \cite[\href{https://stacks.math.columbia.edu/tag/04CV}{Tag 04CV}]{stacks-project} and \cite[\href{https://stacks.math.columbia.edu/tag/05ZE}{Tag 05ZE}]{stacks-project}, we know that $\Omega _{Y/Z}$ is quasi-coherent and of finite type. By \cite[\href{https://stacks.math.columbia.edu/tag/03GA}{Tag 03GA}]{stacks-project} and \cite[\href{https://stacks.math.columbia.edu/tag/03DO}{Tag 03DO}]{stacks-project}, $a^* \Omega _{Y/Z}$ is quasi-coherent and of finite type. Moreover, by  \cite[\href{https://stacks.math.columbia.edu/tag/04CN}{Tag 04CN}]{stacks-project}, we know that $\mathcal{C}_{X/X'}$ is quasi-coherent. 
So by \cite[\href{https://stacks.math.columbia.edu/tag/03M1}{Tag 03M1}]{stacks-project}, $\calhom _{\calO _X} (a^* \Omega _{Y/Z} , \mathcal{C}_{X/X'} )$ is quasi-coherent. So by \cite[\href{https://stacks.math.columbia.edu/tag/03M9}{Tag 03M9}]{stacks-project}, we conclude that  $p_* \calhom _{\calO _X} (a^* \Omega _{Y/Z} , \mathcal{C}_{X/X'} )$ is quasi-coherent. 
Now Corollary \ref{existencerelevcor} follows from Propositions \ref{existencerelev} and \ref{hunzerodiag}. \end{sloppypar}
\end{proof} 

\section{{Formal étaleness, formal smoothness and formal unramifiedness}} \label{fosmo}

Using the results of Section \ref{deformation}, in this section we prove several results on étaleness, smoothness, and unramifiedness in Algebraic Magnetism.

Let $S$ be a scheme and let $M$ be a finitely generated abelian group. Let $X \to Y$ be a $D(M)_S$-equivariant morphism of algebraic spaces over $S$. Let $N\subset M$ be a submonoid.

\prop \label{formsmooth}  Assume that $f:X \to Y$ is formally smooth and locally of finite presentation (i.e. smooth by \cite[\href{https://stacks.math.columbia.edu/tag/04AM}{Tag 04AM}]{stacks-project}), then $f^N:X^N \to Y^N$ is formally smooth as transformation of functors (cf. e.g. \cite[\href{https://stacks.math.columbia.edu/tag/049S}{Tag 049S}]{stacks-project}).
\xprop  

\pf 
Let $\iota_T :\overline{T} \to T $ be a first order thickening of affine schemes, over $S$. Let $D$:
\[\xymatrix{\overline{T} \ar[r]^{\overline{\phi}} \ar[d]^{\iota _T} & X^N \ar[d]^{f^N} \\ T \ar[r]^{\varphi} & Y^N }\] be a commutative diagram of $S$-functors.  We have to prove that there exists an $S$-morphism $\phi : T \to X^N$ such that  the diagram
\[\xymatrix{\overline{T} \ar[r]^{\overline{\phi}} \ar[d]^{\iota _T} & X^N \ar[d]^{f^N} \\ T \ar[ur]^{\phi} \ar[r]^{\varphi} & Y^N }\] commutes.
The diagram $D$ corresponds to the following data (i), (ii)
\begin{enumerate}
\item $\overline{\phi} $ is a $D(M)_{\overline{T}}$-equivariant $\overline{T} $-morphism $A(N)_{\overline{T}} \to X_{\overline{T}} $
\item $\varphi $ is a $D(M)_T$-equivariant $T$-morphism $A(N)_T \to Y_T$
\end{enumerate} 
such that $\varphi |_{\overline{T}} : A(N)_{\overline{T}} \to Y_{\overline{T}}$ equals the map $A(N)_{\overline{T}} \xrightarrow{\overline{\phi}} X_{\overline{T}} \xrightarrow{f_{\overline{T}}} Y_{\overline{T}}$. Let us now consider the commutative diagram \[
\xymatrix{A(N)_{\overline{T}} \ar[d] \ar[r]^{\overline{\phi}}& X_{\overline{T}} \ar[r] &X_T \ar[d]^{f_T} \\ A(N)_T\ar[rr]^{\varphi}& & Y_T}.\]
 The spaces $A(N)_{\overline{T}}$ and $X_{\overline{T}}$ are canonically endowed with actions of $D(M)_T$, moreover all arrows in the diagram are $D(M)_T$-equivariant. By \cite[\href{https://stacks.math.columbia.edu/tag/09ZU}{Tag 09ZU}]{stacks-project}, the left vertical morphism is a first order thickening. By Corollary \ref{existencerelevcor}, we get a $D(M)_T$-equivariant map $\phi : A(N)_T \to X_T $ such that  the diagram \[
\xymatrix{A(N)_{\overline{T}} \ar[d] \ar[r]^{\overline{\phi}}& X_{\overline{T}} \ar[r] &X_T \ar[d]^{f_T} \\ A(N)_T\ar[rr]^{\varphi} \ar[rru]^{\phi}& & Y_T} \] is commutative. The map $\phi$ corresponds to a morphism $T \to X^N$ and satisfies all the required properties to make the diagram \[\xymatrix{\overline{T} \ar[r]^{\overline{\phi}} \ar[d]^{\iota _T} & X^N \ar[d]^{f^N} \\ T \ar[ur]^{\phi} \ar[r]^{\varphi} & Y^N }\] commutative. 
\xpf

\coro \label{corosmoothxy}
 Assume that $X \to Y$ is smooth and that  $X^N$ and $Y^N$ are representable by algebraic spaces. Then $X^N \to Y^N$ is smooth as algebraic spaces.
\xcoro

\pf
Trivial by Propositions \ref{xypres} and \ref{formsmooth},  and \cite[\href{https://stacks.math.columbia.edu/tag/060G}{Tag 060G}]{stacks-project}.
\xpf

\prop \label{formalunra} Assume that $f: X \to Y$ is formally unramified, then $f^N : X^N \to Y^N$ is formally unramified as transformation of functors (cf. e.g. \cite[\href{https://stacks.math.columbia.edu/tag/049S}{Tag 049S}]{stacks-project}). 
\xprop 
\pf
Let $\iota_T :\overline{T} \to T $ be a first order thickening of affine schemes, over $S$. Let $D$:
\[\xymatrix{\overline{T} \ar[r]^{\overline{\phi}} \ar[d]^{\iota _T} & X^N \ar[d]^{f^N} \\ T \ar[r]^{\varphi} & Y^N }\] be a commutative diagram of $S$-functors.  We have to prove that there exists a most one $S$-morphism $\phi : T \to X^N$ such that  the diagram
\[\xymatrix{\overline{T} \ar[r]^{\overline{\phi}} \ar[d]^{\iota _T} & X^N \ar[d]^{f^N} \\ T \ar[ur]^{\phi} \ar[r]^{\varphi} & Y^N }\] commutes.  Let $\phi_a,\phi_b$ be two such morphisms, we have to prove that $\phi_a=\phi_b$.
The diagram $D$ corresponds to the following data (i), (ii)
\begin{enumerate}
\item $\overline{\phi} $ is a $D(M)_{\overline{T}}$-equivariant $\overline{T} $-morphism $A(N)_{\overline{T}} \to X_{\overline{T}} $
\item $\varphi $ is a $D(M)_T$-equivariant $T$-morphism $A(N)_T \to Y_T$
\end{enumerate} 
such that $\varphi |_{\overline{T}} : A(N)_{\overline{T}} \to Y_{\overline{T}}$ equals the map $A(N)_{\overline{T}} \xrightarrow{\overline{\phi}} X_{\overline{T}} \xrightarrow{f_{\overline{T}}} Y_{\overline{T}}$. Let us now consider the commutative diagram \[
\xymatrix{A(N)_{\overline{T}} \ar[d] \ar[r]^{\overline{\phi}}& X_{\overline{T}} \ar[r] &X_T \ar[d]^{f_T} \\ A(N)_T\ar[rr]^{\varphi}& & Y_T}.\] The morphisms $\phi_a,\phi_b$ provide two diagonals of the diagram above. Now using that $f_T$ if formally unramified, we have $\phi_a=\phi_b$.
\xpf

\coro \label{corounramifiedd}
 Assume that $X \to Y$ is formally unramified and that  $X^N$ and $Y^N$ are representable by algebraic spaces. Then $X^N \to Y^N$ is formally unramified as morphism of algebraic spaces.
\xcoro

\pf
Trivial by Proposition \ref{formalunra} and \cite[\href{https://stacks.math.columbia.edu/tag/04G7}{Tag 04G7}]{stacks-project}.
\xpf

\prop \label{propétfor} Assume that $f: X \to Y$ is formally étale and locally of finite presentation (i.e étale by  \cite[\href{https://stacks.math.columbia.edu/tag/0616}{Tag 0616}]{stacks-project}), then $f^N : X^N \to Y^N$ is formally étale as transformation of functors (cf. e.g. \cite[\href{https://stacks.math.columbia.edu/tag/049S}{Tag 049S}]{stacks-project}). 
\xprop 

\pf
This follows from \cite[\href{https://stacks.math.columbia.edu/tag/049S}{Tag 049S}]{stacks-project} and Propositions \ref{formalunra} and \ref{formsmooth}.
\xpf

\coro \label{coroétt}
 Assume that $X \to Y$ is étale and that  $X^N$ and $Y^N$ are representable by algebraic spaces. Then $X^N \to Y^N$ is étale as morphism of algebraic spaces.
\xcoro

\pf
Trivial by Proposition \ref{propétfor},  \cite[\href{https://stacks.math.columbia.edu/tag/04GC}{Tag 04GC}]{stacks-project} and \cite[\href{https://stacks.math.columbia.edu/tag/0616}{Tag 0616}]{stacks-project}.
\xpf

\prop \label{formsmoothface}  Assume that $f:X \to S$ is smooth, let $F$ be a face of $N$,  then $f^{N,F}:X^N \to X^F$ is formally smooth as transformation of functors (cf. e.g. \cite[\href{https://stacks.math.columbia.edu/tag/049S}{Tag 049S}]{stacks-project}).
\xprop

\pf 
Let $\iota_T :\overline{T} \to T $ be a first order thickening of affine schemes, over $S$. Let $R$ and $I \subset R$ such that $T= \Spec (R)$ and $\overline{T} = \Spec (R/I)$. Let $D$:
\[\xymatrix{\overline{T} \ar[r]^{\overline{\phi}} \ar[d]^{\iota _T} & X^N \ar[d]^{f^{N,F}} \\ T \ar[r]^{\varphi} & X^F }\] be a commutative diagram of $S$-functors.  We have to prove that there exists an $S$-morphism $\phi : T \to X^N$ such that  the diagram
\[\xymatrix{\overline{T} \ar[r]^{\overline{\phi}} \ar[d]^{\iota _T} & X^N \ar[d]^{f^{N,F}} \\ T \ar[ur]^{\phi} \ar[r]^{\varphi} & X^F }\] commutes.
The diagram $D$ corresponds to the following data (i), (ii)
\begin{enumerate}
\item $\overline{\phi} $ is a $D(M)_{\overline{T}}$-equivariant $\overline{T} $-morphism $A(N)_{\overline{T}} \to X_{\overline{T}} $
\item $\varphi $ is a $D(M)_T$-equivariant $T$-morphism $A(F)_T \to X_T$
\end{enumerate} 
such that $\varphi |_{\overline{T}} : A(F)_{\overline{T}} \to Y_{\overline{T}}$ equals the map $A(F)_{\overline{T}} \to A(N)_{\overline{T}} \xrightarrow{\overline{\phi}} X_{\overline{T}} $. We need to prove the following fact.
\fact \label{fact}
The push-out of the diagram $A(F)_T \leftarrow A(F)_{\overline{T}} \rightarrow A(N)_{\overline{T}}$, in the category of algebraic spaces over $T$, exists and is denoted \[A(F)_T \coprod _{A(F) _{\overline{T}} }A(N) _{\overline{T}}.\] Moreover $A(F)_T \coprod _{A(F) _{\overline{T}} }A(N) _{\overline{T}}$ is an affine scheme over $T$ and is in fact a push-out of \[{A(F)_T \leftarrow A(F)_{\overline{T}} \rightarrow A(N)_{\overline{T}}}\]in the category whose objects are algebraic spaces over $T$ and whose morphisms are $D(M)_T$-equivariant morphisms of algebraic spaces over $T$.
\xfact 
\pf
By \cite[\href{https://stacks.math.columbia.edu/tag/0ET0}{Tag 0ET0}]{stacks-project}, the push-out of $A(F)_T \leftarrow A(F)_{\overline{T}} \rightarrow A(N)_{\overline{T}}$ exists in the category of schemes and is given by 
\[ A(F)_T \coprod _{A(F) _{\overline{T}} }A(N) _{\overline{T}} := \Spec \big( \bigoplus _{n \in F} R.X^n \oplus \bigoplus _{n \in N \setminus F} (R/I). X^n \big). \] 
Since $D(M)_T$-equivariant morphisms of affine schemes correspond to $M$-graded morphisms, we obtain that $ A(F)_T \coprod _{A(F) _{\overline{T}} }A(N) _{\overline{T}}$ is endowed with an action of $D(M)_T$ and that
 \[{\nu _1 : A(F)_T \to A(F)_T \coprod _{A(F) _{\overline{T}}}A(N)_{\overline{T}}}\] and
  \[\nu _2 : A(N)_{\overline{T}} \to A(F)_T \coprod _{A(F) _{\overline{T}}}A(N)_{\overline{T}}\] are $D(M)_T $-equivariant. By \cite[\href{https://stacks.math.columbia.edu/tag/0ET0}{Tag 0ET0}]{stacks-project}, the push-out of \[D(M)_T \times _T A(F)_T \leftarrow D(M)_T \times _T A(F)_{\overline{T}} \rightarrow D(M)_T \times _T A(N)_{\overline{T}}\] exists in the category of schemes and the explicit formula allows us to identify it canonically with \[D(M)_T \times _T \big( A(F)_T \coprod _{A(F) _{\overline{T}} }A(N) _{\overline{T}} \big) .\]
  We now remark that push-outs in the category of schemes over $T$  give push-outs in the full category of algebraic spaces over $T$ by  \cite[\href{https://stacks.math.columbia.edu/tag/07SY}{Tag 07SY}]{stacks-project}. Now let $X$ be an algebraic space over $T$ with $D(M)_T$-equivariant morphisms $f_1: A(F)_T \to X$  and  $f_2: A(N)_{\overline{T}} \to X $ such that the diagram \[\xymatrix{ A(F)_{\overline{T}} \ar[r]^{f_1} \ar[d] & A(N)_{\overline{T}}^{f_2} \ar[d] \\ A(F)_T \ar[r] &X } \]
   commutes; we have to prove that the obtained morphism $f_1 \coprod f_2: A(F)_T \coprod _{A(F) _{\overline{T}} }A(N) _{\overline{T}}  \to X $ is $D(M)_T$-equivariant. For this consider the diagram\[
 \begin{tikzcd} D(M)_T \times _T A(F)_{\overline{T}} \ar[ddd, "Id \times \iota"]\ar[dddddd, " m_{A(F)_{\overline{T}}}"' near start ,bend right=60 ] \ar[rr, " Id \times i "]& & D(M)_T \times _T A(N)_{\overline{T}} \ar[dddddd, " m_{A(N)_{\overline{T}}}"' near start ,bend right=60 ] \ar[ddd, " Id \times a_2"] \ar[rrddd, "Id \times f_2" ]& & \\
 & & & & \\
 & & & &\\
 D(M)_T \times A(F)_T \ar[rrrr, "Id \times f_1"' near end, bend right =10] \ar[dddddd, " m_{A(F)_{{T}}}"'  ,bend right=60 ]\ar[rr, "Id \times a_1" ]& & D(M)_T \times _T \big( A(F)_T \coprod _{A(F) _{\overline{T}} }A(N) _{\overline{T}} \big) \ar[rr, "Id \times (f_1 \coprod f_2)" ]  \ar[dddddd, " m_{A(F)_T \coprod _{A(F) _{\overline{T}} }A(N) _{\overline{T}}}"' near end ,bend right=60 ]& & D(M)_T \times _T X \ar[dddddd, " m_X"   ]  \\
 & & & & \\
 & & & &\\ 
 A(F)_{\overline{T}} \ar[ddd, "\iota"] \ar[rr, "i"] && A(N)_{\overline{T}} \ar[ddd, "a_2" ]  \ar[rrddd, " f_2" ] & & \\
 & & & &\\
 & & & & \\
 A(F)_T \ar[rrrr, " f_1"' near end, bend right =10] \ar[rr , " a_1" ] & &A(F)_T \coprod _{A(F) _{\overline{T}} }A(N) _{\overline{T}}  \ar[rr , " f_1 \coprod f_2 "] & & X 
  \end{tikzcd}\]
  where $m_{\calS}$ denote the multiplication morphism for any $D(M)_T$-space $\calS$.  
 We have to show that $ m_X \circ \big( Id\times (f_1 \coprod f_2) \big) = (f_1 \coprod f_2) \circ m_{A(F)_T \coprod _{A(F) _{\overline{T}} }A(N) _{\overline{T}}}$. Since, as we noted before, $D(M)_T \times _T (A(F)_T \coprod _{A(F)_{\overline{T}}} A(N)_{\overline{T}})$ is the push-out of \[D(M)_T \times _T A(F)_T \xleftarrow{Id \times \iota} D(M)_T \times _T A(F)_{\overline{T}} \xrightarrow{Id \times i } D(M)_T \times _T A(N)_{\overline{T}},\] it is enough to prove that
 \begin{enumerate}
 \item $ m_X \circ \big( Id\times (f_1 \coprod f_2) \big) \circ ( Id \times a_1) = (f_1 \coprod f_2) \circ m_{A(F)_T \coprod _{A(F) _{\overline{T}} }A(N) _{\overline{T}}} \circ ( Id \times a_1 ) $
 \item $ m_X \circ \big( Id\times (f_1 \coprod f_2) \big) \circ ( Id \times a_2) = (f_1 \coprod f_2) \circ m_{A(F)_T \coprod _{A(F) _{\overline{T}} }A(N) _{\overline{T}}} \circ ( Id \times a_2 ) $.
 \end{enumerate}
 These identities are easy to check on the diagram using that $a_1,a_2,f_1,f_2$ are equivariant and that $f_1= (f_1 \coprod f_2 ) \circ a_1$, $ f_2 = (f_1 \coprod f_2) \circ a_2$.
\xpf

  We now continue the proof of Proposition \ref{formsmoothface}.
Let us consider the diagram 
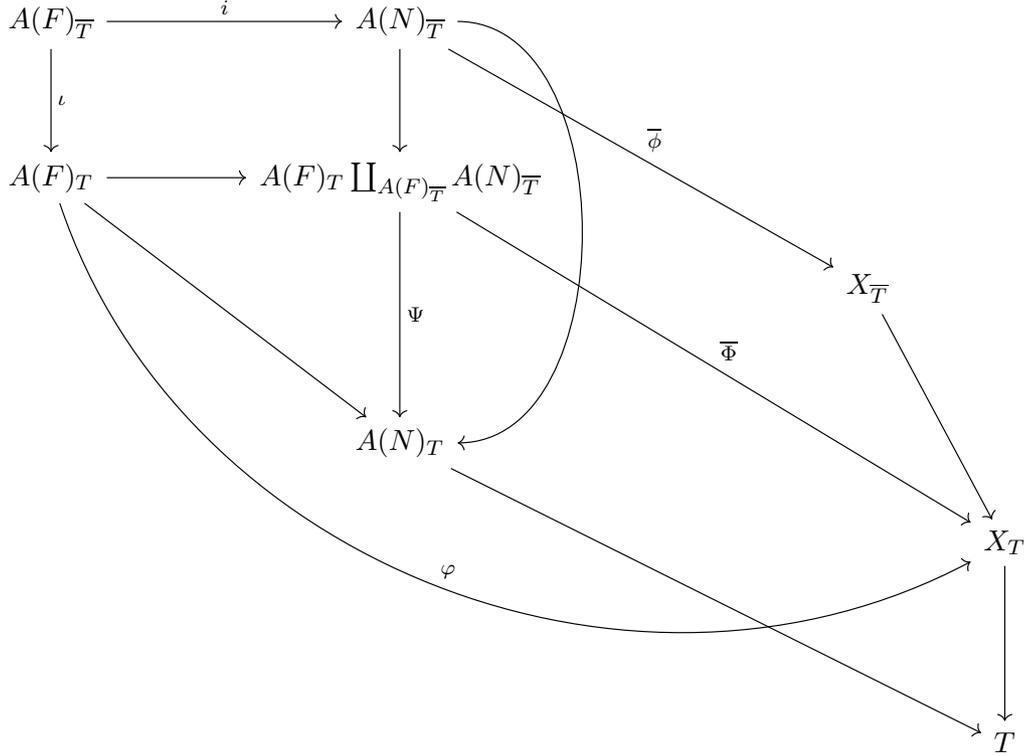
\[
\begin{tikzcd}A(F)_{\overline{T} }\ar[rr, "i"] \ar[dd, "\iota"] & & A(N) _{\overline{T}} \ar[dd] \ar[ddddd, bend left=90]\ar[rrrrddd, " \overline{\phi}"] & & & & &\\
 & & & & & & & \\
A(F)_T  \ar[rrddd]\ar[rrrrrrrdddd, bend right=50 ,"\varphi"] \ar[rr]& & A(F)_T \coprod _{A(F)_{\overline{T }}}A(N)_{\overline{T}} \ar[ddd, " \Psi"] \ar[rrrrrdddd, "\overline{\Phi}"] &  & & & &  \\
& & & &  & &X_{\overline{T}} \ar[dddr] &\\
& & & & & & &\\
& &  A(N)_T \ar[rrrrrdddd]& & & & &\\
& & & & & & &X_{T} \ar[ddd]   \\
& & & & & & &\\
& & & & & & &\\
& & & & & & &T \end{tikzcd}
\] 
where \begin{enumerate} 
\item $\overline{\Phi}$ is the map obtained from the universal property of pushout from $\overline{\phi}$ and $\varphi$.
\item $\Psi$ is the map obtained from the universal property of pushout from the canonical morphisms $A(N)_{\overline{T}} \to A(N)_T$ and $A(F)_T \to A(N)_T$. 
\end{enumerate} By Fact \ref{fact}, $i$ and $\iota$ are $D(M)_T$-equivariant and so $\Psi$ and $\overline{\Phi}$ are also $D(M)_T $-equivariant.

\lemm \label{cofot}
The morphism
$
A(F)_T \coprod_{A(F)_{\overline{T}}} A(N)_{\overline{T}} \longrightarrow A(N)_T
$
is a first order thickening.
\xlemm

\pf
The statement is local on $T$, so we may assume $T = \Spec(R)$ with $\mathcal{I} = I \subset R$ and $I^2 = 0$. Then $\overline{T} = \Spec(R/I)$.
By construction, the push-out
$
A(F)_T \coprod_{A(F)_{\overline{T}}} A(N)_{\overline{T}}
$
is affine and corresponds to the $R$-algebra
$
R[F] \oplus (R/I)[N \setminus F],
$
together with the natural morphism of $R$-algebras
$
R[N] \longrightarrow R[F] \oplus (R/I)[N \setminus F].
$
Let $\mathcal{K}$ denote the kernel of this morphism. Then an element of $\mathcal{K}$ is a finite $R$-linear combination supported on $N \setminus F$ whose coefficients lie in $I$, hence
$
\mathcal{K} = (I \cdot R[N \setminus F]).
$
Since $I^2 = 0$, it follows that $\mathcal{K}^2 = 0$.
Therefore the morphism is a closed immersion defined by a square-zero ideal, i.e.\ a first order thickening.
\xpf
 Now we apply Lemma \ref{cofot} and Corollary \ref{existencerelevcor} to get a $D(M)_T$-equivariant morphism $\phi : A(N)_T \to X_T $ such that the right lower quadrigone commutes. The morphism $\phi $ satisfies all required properties to make the diagram \[\xymatrix{\overline{T} \ar[r]^{\overline{\phi}} \ar[d]^{\iota _T} & X^N \ar[d]^{f^{N,F}} \\ T \ar[ur]^{\phi} \ar[r]^{\varphi} & X^F }\] commutative.
\xpf

\coro \label{123456789} Let $F $ be a face of a magnet $N$.
 Assume that $X \to S$ is smooth and that  $X^N$ and $X^F$ are representable by algebraic spaces. Then $X^N \to X^F$ is smooth as algebraic spaces.
\xcoro

\pf
Clear by Proposition \ref{formsmoothface}, and \cite[\href{https://stacks.math.columbia.edu/tag/060G}{Tag 060G}]{stacks-project}.
\xpf

\rema \label{flatvistoli}If $X\to Y$ is flat, then $X^N \to Y^N$ is not flat in general (as we know that flatness is not preserved by taking fixed-points, e.g. cf. Vistoli's answer in Mathoverflow ``Under what hypotheses are schematic fixed points of a flat deformation themselves flat?").
\xrema

\rema \label{cmvistoli}If $X $ is Cohen-Macaulay, then $X^N  $ is not Cohen-Macaulay in general (as we know that being Cohen-Macaulay is not preserved by taking fixed-points, e.g. cf. Vistoli's answer in Mathoverflow ``Are schematic fixed-points of a Cohen-Macaulay scheme Cohen-Macaulay?").

\xrema

\section{Topology and geometry of attractors} \label{topo}

In this section, we prove several results on the geometry and topology of attractors, in particular we prove a Białynicki-Birula decomposition in our setting.

Let $M$ be a finitely generated abelian group. Let $X$ be an $S$-algebraic space locally of
finite presentation endowed with an action of $D(M)_S$.

\prop \label{proptopo}  Let $N\subset M$ be a submonoid. Assume that $N$ is finitely generated and that there exists a $N^*$-FPR atlas of $X$. Then the map $X^N \to X^{N^*}$ (which is affine by Proposition \ref{representable}) has geometrically connected fibers and induces a bijection on the sets of connected components $\pi _0 (X^N) \simeq \pi _0 (X^{N^*})$ of the underlying spaces.

\xprop

\pf We adapt \cite[Cor. 1.12]{Ri16}. Using Proposition \ref{propequiv} we can assume that $N^{*}= 0$ and that $N$ is fine and sharp. Let $K$ be a field, and let $x: \Spec (K) \to X^{0}$ be a point. Let $X^N_x= X^N \times _{X^0 , x } \Spec (K)$. We claim that its underlying topological space $|X^N_x|$ is connected. Let $L$ be a field and let $y: \Spec (L) \to X^N_x$ be a point, and denote by $x_L$ the composition $\Spec(L) \to \Spec (K) \xrightarrow{x} X^0.$ Then $x_L$ and $x$ define the same point of $|X^N_x|$.
 Recall that we have a natural action of the monoid scheme $A(N)_S$ on the attractor $X^N$. The $A(N)_{L}$ orbit of $y$ defines a map $h: A(N)_L \to X^N_x$ with $h(1)= y$ and $h(0)=x_L$. Since $N $ is sharp, $A(N)_L$ is connected. So $x$ and $y$ lie in the connected set $|h|(|A(N)_L|)$. Since $y$ was arbitrary, this shows that $|X^N_x|$ is connected. So the continous map $|X^{N}| \to |X^0|$ has connected fibers, and the assertion on connected components follows from the existence of a continuous section $|X^0| \subset |X^{N}| $.
\xpf

\prop \label{opclosmoo}
Let $N \subset M $ be a magnet. Let $f:X\to Y $ be a $D(M)_S$-equivariant morphism of $S$-algebraic spaces such that $X^N $ and $Y^N$ are representable by $S$-algebraic spaces. Let $f^N:X^N \to Y^N$ be the canonical morphism obtained on attractors.  The following assertions hold.
\begin{enumerate}
\item If $f$ is an open immersion, then $f^N$ is an open immersion.
\item If the following conditions hold \begin{enumerate} \item $f$ is a closed immersion, \item  there exists a $N^*$-FPR atlas of $Y$, 
\item  $N$ is finitely generated or $X$ is separated, \end{enumerate} then $f^N$ is a closed immersion, moreover $X^N \cong Y^N \times _Y X$.
\end{enumerate}
\xprop

\pf \begin{enumerate} \item Since $f$ is an open immersion, it is smooth and a monomorphism. By Corollary \ref{corosmoothxy}, $f^N$ is smooth  and in particular locally finitely presented and flat (cf.  \cite[\href{https://stacks.math.columbia.edu/tag/04TA}{Tag 04TA}]{stacks-project}). By Fact \ref{monomono}, $f^N$ is a monomorphism. Therefore, by \cite[\href{https://stacks.math.columbia.edu/tag/05VH}{Tag 05VH}]{stacks-project}, $f^N$ is universally injective and unramified. Consequently, by \cite[\href{https://stacks.math.columbia.edu/tag/06LU}{Tag 06LU}]{stacks-project}, $f^N$ is étale. Finally, by \cite[\href{https://stacks.math.columbia.edu/tag/05W5}{Tag 05W5}]{stacks-project}, $f^N$ is an open immersion.
\item
Assume $f$ is a closed immersion. The canonical morphisms $X^N \to Y^N $ and $X^N \to X$ induces a canonical morphism $i: X^N \to Y^N \times _Y X$. We are going to prove that $i$ is an isomorphism. Let $U$ be a $N^*$-FPR atlas of $Y$. 
By Proposition \ref{clofpr}, the map $U \times _Y X \to X$ is a $N^*$-FPR atlas of $X$. So by Proposition \ref{fixxx} the canonical map $(U \times _Y X)^N \to X^N$ is étale and surjective. The canonical map $U^N \to Y^N $ (étale and surjective) induces an étale and surjective map $(U \times _Y X) \times _U U^N = X \times _Y U^N \to  X \times _Y Y^N $.
So we get a diagram 
\[
\xymatrix{(U \times _Y X)^N \ar[r] \ar[d] & X^N \ar[d]  \\
 (U \times _Y X) \times _U U^N = X \times _Y U^N \ar[r] & X \times _Y Y^N.}\]
The left arrow is an isomorphism by Lemma \ref{equiclosed} (write $U$ as a disjoint union of stable $S$-affine schemes $U_i$ and use that $U_i \times _Y X \to U_i \times _Y Y=  U_i $ is affine for all $i$). The horizontal arrows are étale and surjective. The diagram is a cartesian square because \[(X \times _Y U^N ) \times _{X \times _Y Y^N} X^N = U^N \times _{Y^N} X^N = (U \times _Y X)^N \] by Proposition \ref{produitNNN}. So by \cite[\href{https://stacks.math.columbia.edu/tag/03I2}{Tag 03I2}]{stacks-project} we obtain $X^N \cong X \times _{Y} Y^N$.\end{enumerate}
\xpf

\prop \label{retract} Let $N$ be a magnet of $M$ and assume $N$ is fine as monoid. Let $F$ be a face of $N$. Let $X^F \xrightarrow{\iota} X^N$ and $X^N \xrightarrow{r} X^F$ be the associated morphisms (cf. \ref{morphismsgeneral} and \ref{faceattractors}).
Then $\iota,r $ is a strong deformation retract. That is there exists a morphism $\delta: X^N \times _S \bbA ^1 _S \to X^N $ such that 
\begin{enumerate}
\item $\delta \circ j_0 = \iota \circ r $ where $j_0 : X^N  \to X^N \times _S \bbA ^1 _S$ is given by the section  $0$  of $\bbA^1_S$,
\item  $\delta \circ j_1 = Id_{X^N} $ where $j_1 : X^N \to X^N \times _S \bbA ^1 _S $ is given by the section $1$ of $\bbA^1_S$,
\item $\delta \circ ( \iota \times Id_{\bbA^1 _S}) = \iota \circ pr_1$ where $pr_1: X^F \times _S \bbA^1_S \to X^F$ is the projection.
\end{enumerate}
\xprop 
\pf
Let $f: A(N) _S \times _S \bbA ^1 _S \to A(N)_S $ be a map as in \cite[Prop. 3.3.1 (iii) p. 75]{Og}. Note that $f$ is $D(M)_S$-equivariant, i.e. $f (g \cdot x , y) = g \cdot f(x,y)$ by construction of $f$ in \cite[Proof of Prop. 3.3.1 (iii) p. 75]{Og}. We now define $\delta$ using $f$. Let $T/S$ be a scheme and let us take an element $x$ in $(X^N \times _S \bbA^1_S) (T)$. This element $x$ corresponds to a pair of morphisms $(b,c)$ where $b: A(N)_T \to X_T$ is $D(M)_T$-equivariant and $c: T  \to \bbA^1_T  $ belongs to $\bbA^1_T (T).$ 
Let $c'$ be the composition 
\[c': A(N)_T = A(N)_T \times _T T \xrightarrow{Id \times c} A(N)_T \times _T \bbA^1_T \xrightarrow{f} A(N)_T . \]
The morphism $c' : A(N)_T \to A(N)_T $ is $D(M)_T $ equivariant.
Now we define $\delta(x)$ as the morphism $T \to X^N$ corresponding to the $D(M)_T$-equivariant composition $\delta(x): A(N)_T \xrightarrow{c'} A(N)_T \xrightarrow{b} X_T $. This defines $\delta$ as functor.
Now \cite[Prop. 3.3.1 (iii) p. 75]{Og} implies that (i), (ii) and (iii) of Proposition \ref{retract} hold.
\xpf 

The following theorem generalizes the Bialynicki-Birula decomposition.

\theo \label{BB} Assume that $X$ is smooth over $S$ and that the action $a$ of $D(M)_S$ on $X$ is Zariski locally linearizable (cf. §\ref{Zariskicase}). Assume that $N$ is finitely generated as monoid or that $X$ is quasi-compact over $S$.
Then the morphism $X^N \to X^{N^*}$ is an affine space fibration.
\xtheo 

\pf
 Since the action $a$ is Zariski locally linearizable, we reduce to the case where $X$ is affine over $S$. Then $X^N$ and $X^{N^*}$ are affine over $S$. Using Proposition \ref{existfini} we reduce to the case where $N$ is finitely generated. Using Proposition \ref{propequiv} (to $M \to Z= M/ N^{*}$ and $Y=N/N^*$), we then reduce to the case where $N$ is fine and sharp, in particular $N^*=0$. 
By \ref{actionmono} $A(N)_S$ acts on $X^{N}$. So there exists a quasi-coherent $\calO _S$-algebra $A= \bigoplus _{n \in N } \calA _n $ graded by the fine and sharp monoid $N$ such that $X^{N} = \Spec _S ( \cal A)$. By \ref{NLNL}, we have $X^0= (X^N)^0$. Since $N $ is sharp, the ideal $\calI$ generated by $\{\calA _n | n \in N \setminus 0\}$ identifies with  $\bigoplus _{n \in N \setminus 0 } \calA _n $ and $\calA / \calI = \calA _0$. Then by Theorem \ref{representableaffine}, $X^{0} = \Spec _S ( \calA _0)$ and we have a closed immersion $X^0 \to X^N$. By \ref{corosmoothxy}, $X^N \to S $ and $ X^0 \to S $ are smooth, so by  \cite[\href{https://stacks.math.columbia.edu/tag/067U}{Tag 067U}]{stacks-project} $X^0 \to X^N$ is a regular closed immersion. So by \cite[\href{https://stacks.math.columbia.edu/tag/063M}{Tag 063M}]{stacks-project} the conormal sheaf $\calI / \calI ^2$ is finite locally free hence projective. The ideal $\calI$ is graded by definition. The definition of $\calI $ also implies that $\calI^2$ is graded. 
Therefore, the quotient
map $\calI \to \calI /\calI ^2$ is a homomorphism of graded $\calA$-modules. Hence, $\calA _n$ maps onto
$(\calI / \calI^2)_n$ for every nonzero $n \in N $. The $\calA _0$-module $(\calI / \calI^2)_n$ is projective because it is a direct summand of $\calI / \calI ^2$,
 and so the surjection $\calA _n \to (\calI / \calI^2)_n$ admits a section  $\fraks _n : (\calI / \calI^2)_n \to \calA_n.$ 
 Let $\fraks: \calI /  \calI ^2 \to \calA$ be the direct sum $\bigoplus _{n \in N}\fraks_n$, note that $\fraks(\calI / \calI ^2) \oplus \calI ^2 = \calI$.
Now $\fraks$ induces a morphism of $\calA _0$-algebras $\Phi: Sym _{\calA _0} (\calI /\calI^2) \to \calA$. It is enough to show that this morphism is an isomorphism. We note that $Sym _{\calA _0} (\calI /\calI ^2)$ is the sheafification of the presheaf $U \mapsto Sym _{\calA _0 (U)} (\calI (U) /\calI ^2(U))$, cf. e.g.  \cite[\href{https://stacks.math.columbia.edu/tag/01CG}{Tag 01CG}]{stacks-project}. So, using \cite[\href{https://stacks.math.columbia.edu/tag/007V}{Tag 007V}]{stacks-project}, we reduce to the case where $\calA$ is a ring. We first prove that $\Phi$ is surjective.
By \cite[Prop. 2.2.1 p. 34]{Og}, there exists a morphism of monoids $h:N \to \bbN $ such that $h^{-1}(0)=0$. The image of $Sym _{\calA _0} (\calI /\calI ^2) \to \calA$ is the $\calA _0$-subalgebra $\calA _0 [\fraks (\calI /\calI ^2)]$ of $\calA$  generated by $\fraks (\calI /\calI ^2)$. Let $c \in N$ and $a \in \calA _c$. We shall prove by induction on $h(c)$ that $a \in \calA _0 [\fraks (\calI /\calI ^2)]$. 
Certainly this is true if $c=0$. If $h(c)=1$, then $a$ belongs to $\fraks(\calI /\calI ^2)$ and so $a \in \calA _0 [\fraks (\calI /\calI ^2)]$. We now assume that $\calA _c \subset \calA _0 [\fraks (\calI /\calI ^2)]$ for all $c $ such that $h(c) \in \{ 0 , \ldots , k \}$ with $k \geq 1$. Let $a \in \calA_c$ with $h(c) =k+1$ and let us prove that $a $ belongs to $\calA _0 [\fraks (\calI /\calI ^2)]$. There exists $a' \in \fraks (\calI / \calI^2)$ such that
$a-a' \in \calI ^2$, and so we may suppose that $a \in \calI^2$. Write $a= \sum _i b_i d_i$ with $b_i$ and $d_i$ homogeneous elements of $\calI$. The equality still holds when we omit any terms $b_i d_i $ with $deg(b_i)+ deg (d_i) \not = c$. For the remaining terms, $h(deg(b_i)), h(deg(d_i))< h(c)$ and so $b_i, d_i  \in \calA _0 [\fraks (\calI /\calI ^2)]$. So $\Phi$ is surjective. We then prove that $\Phi$ is injective.
 By \cite[\href{https://stacks.math.columbia.edu/tag/062Z}{Tag 062Z}]{stacks-project} and \cite[\href{https://stacks.math.columbia.edu/tag/063M}{Tag 063M}]{stacks-project} (cf. also \cite[\href{https://stacks.math.columbia.edu/tag/063H}{Tag 063H}]{stacks-project}), the surjective morphism $\Theta: Sym _{\calA _0} (\calI /\calI^2) \to Gr _{\calA _o} (\calI) = \bigoplus _{k \geq } \calI^k / \calI ^{k+1}$ is an isomorphism. The induced morphism of algebras $\Psi: \bigoplus _{k \geq } \calI^k / \calI ^{k+1} \to \calA$ is explicitly given by $\Psi ([i_1 \ldots  i_k]_{\mathrm{mod} \calI ^{k+1}}) \mapsto \fraks([i_1]_{\mathrm{mod} \calI ^2})\ldots \fraks([i_k]_{\mathrm{mod} \calI ^2})$. Let $\Psi _k$ be the restriction of $\Psi $ to $\calI^k / \calI ^{k+1}$, it factors through $\calI^k$. Let $\pi _k$ be the reduction modulo $\calI _{k}$ for any $k \in \bbN$. Then for $d\geq k $, the composition $\pi _d \circ \Psi _k : \calI ^k / \calI ^{k+1} \to \calA / \calI^d$ equals $0$ if $d \leq k$. Moreover, the composition $\pi _{k+1} \circ \Psi _k : \calI ^k / \calI ^{k+1} \to \calI^k/ \calI^{k+1}$ is the identity. We now prove that $\Psi $ is injective. So let $X = \sum _{k=0}^{n} i_k$ with $i_k \in \calI^k/\calI ^{k+1}$ and $\Psi (X)=0$. We now prove by induction on $j$ that $i_0, i_1, ... , i_j$ equal zero. We have $0=\Psi (X) \mathrm{mod}\calI =\sum _{k=0} ^{n} \Psi (i_k) \mathrm{mod} \calI = \Psi (i_0) \mathrm{mod} \calI = i_0$. So $i_0=0$ and $X= \sum _{k=1}^{n} i_k$. Now assume that $i_0, i_1, ... , i_j$ equal zero. We have $0= \Psi (X) \mathrm{mod}\calI^{j+1} =\sum _{k={j+1}} ^{n} \Psi (i_k) \mathrm{mod} \calI^{j+1} = \Psi (i_{j+1}) \mathrm{mod} \calI^{j+1} = i_{j+1}$. This shows that $\Psi$ and $\Phi$ are injective and finishes the proof.
\xpf 
 We obtain the following corollary, in the spirit of a classical result on toric varieties over fields (cf. e.g. \cite[Theorem 1.3.12, Definition 1.2.16]{CLS11}). 
\coro \label{smoothANisvector}
Let $S$ be a scheme. Let $N$ be a fine and sharp monoid. If $A(N)_S \to S$ is smooth, then $A(N)_S \to S$ is a vector bundle.
\xcoro 
\pf
Let $M$ be $N^{gp}$, then $(A(N)_S)^N =A(N)_S$ and $(A(N)_S)^0 =S $. Now Corollary \ref{smoothANisvector} follows from the proof of Theorem \ref{BB}.
\xpf

\section{Attractors and dilatations}  \label{dilatationse}
In this section, we study the compatibility between attractors and dilatations.

Let $S$ be a scheme and let $S'$ be a closed locally principal subscheme of $S$. 
Let $X$ be an algebraic space over $S$ with a $D(M)_S$-action where $M$ is an abelian group. Put $D= X _{S'}$. 
 Let $Y$ be a closed subspace of $X_{S'}$. 
Then by \cite{Ma23} (or \cite{MRR20} if $X$ is a scheme), we get a space $\Bl _Y^D X$ called the dilatation of $X$ with center $Y$ along $S'$, and an affine morphism of spaces $\Bl_Y^DX \to X$. Let $N$ be a submonoid in $M$.
Assume that $Y$ is stable under the action of $D(M)_{S'}$ on $X_{S'}$. Then by Proposition \ref{opclosmoo} $Y^N \to D^N$ is a closed immersion. Moreover $D^N=(X_{S'})^N=(X^N)_{S'}$ (cf. Proposition \ref{basechangeee}) is a locally principal closed subscheme of $X^N$. So $\Bl_{Y^N}^{D^N} X^N$ is well-defined.

\prop \label{dilattra}
 Assume moreover that $X \to S$, $(\Bl_{Y}^DX )^N\to S$ and $ \Bl _{Y^N}^{D^N} X^N\to S$ are flat. Then $D(M)_S$ acts naturally on $\Bl_{Y}^DX$, moreover we get a canonical isomorphism  
 \[ \Theta: (\Bl _Y^D X )^N \cong \Bl _{Y^N}^{D^N} X^N .\] 
\xprop

\rema We refer to  \cite[Prop. 2.16]{MRR20} for conditions ensuring flatness of dilatations.
\xrema 

 \begin{proof} Remark first that since  $D(M)_S \to S$,  $X \to S$, $\Bl_Y ^D X \to S$ are flat, $D(M)_S$, $X$ and $\Bl_Y ^D X $ belong to $\mathrm{Spaces}_{S}^{S'\text{-}reg}$ (cf. \cite{Ma23} for the definition of  $\mathrm{Spaces}_{S}^{S'\text{-}reg}$), moreover by flatness any products of these objects in the category $\mathrm{Spaces}/S$ or $\mathrm{Spaces}_{S}^{S'\text{-}reg}$ coincide. 
 So to check that we have an action of $D(M)_S$ on $\Bl_Y ^D X$ it is enough to show that $D(M)_S(T) $ acts on $\Bl_Y ^D X(T)$ functorially for any $T \in \mathrm{Spaces}_{S}^{S'\text{-}reg}$. So let $T \in \mathrm{Spaces}_{S}^{S'\text{-}reg}$. 
 Let $(g,x) \in D(M)_S(T) \times \Bl_Z^D(X)(T)$, note that $x$ corresponds to a morphism $T \xrightarrow{x} X$ such that $T|_{S'}\to X_{S'}$ factors through $Y$. We define $g.x $ as the composition 
 
 \[ T \xrightarrow{(g,x)} D(M)_S \times X \xrightarrow{\text{action}} X .\]
 Then $g.x$ restricted to $S'$ factors through $Y$ and so $g.x \in \Bl_Z^DX (T)$. So $D(M)_S$ acts on $\Bl_Y^D(X)$.
 We obtain the following two diagrams 
\[ (\Bl _Y ^D X )^N \to X^N ,\]

 \[ \begin{tikzcd} (\Bl _Y^D X)^N | _{S'}  =  (\Bl _Y^D X | _{S'} )^N\ar[rr] \ar[rd] & & X^N |_{S'} = (X | _{S'}) ^N \\ & Y^N \ar[ru] & \end{tikzcd} .\]
 By the universal property of dilatations, we obtain a morphism  \[ \Theta: (\Bl ^D_Y X )^N \to \Bl ^{D^N}_{Y^N} X^N .\] We now prove that it is an isomorphism.
 Again, let $T \in \mathrm{Spaces}_S^{S'-reg}$ and let $T'$ be $T\times _S S'$. Then $T'$ is a closed and locally principal subspace of $T$.  Moreover, since $A(N)_T$ is flat over $T$, $A(N)_T $ belongs to $\mathrm{Spaces}_{T}^{T ' -reg}$. Since $\Bl _Y^D X \to S$ is flat, by \cite[§3.6]{Ma23}, we have $(\Bl _Y^D X)_T = \Bl _{Y_{T'}}^{D_{T'}}X_T$.
 So, on the one hand
 \begin{align*}
 ( \Bl_Y ^D X )^N (T) &= \Hom ^{D(M)_T} (A(N)_T , (\Bl _Y^D X)_T ) \\
 &= \Hom ^{D(M)_T} (A(N)_T , \Bl _{Y_{T'}}^{D_{T'}}X_T ) \\  
 & = \{ A(N)_T \xrightarrow{f} X_T |~ f \text{ is } D(M)_T\text{-equiv. and } A(N)_T|_{T'} \xrightarrow{f_{T'}} X|_{T'} \text{ factors through } Y_{T'} \},
 \end{align*}
 moreover, on the other hand,
 \begin{align*}
 (\Bl ^{D^N}_{Y^N}X^N) (T)&=\{T \to X^N |~~ T|_{S'}\to X^{N}|_{S'} \text{ factors through } Y^N \}\\
 &=\{A(N)_T \xrightarrow{f} X | ~f \text{ is }D(M)_T\text{-equiv. and }A(N)_T|_{T'} \xrightarrow{f_{T'}} X|_{T'} \text{ factors through } Y_{T'}\}.
 \end{align*}
 Now since $\Bl ^{D^N}_{Y^N}X^N$ and $( \Bl_Y ^D X )^N $ are flat over $S$, they belong to $\mathrm{Spaces}_S^{S'\text{reg}}$. So by Yoneda $\Bl ^{D^N}_{Y^N}X^N = ( \Bl_Y ^D X )^N $.  This finishes the proof.
 \end{proof}

 \coro \label{dilamagnecor}Let $S$ be a scheme and let $S'$ be an effective Cartier divisor on $S$. 
Let $X$ be a smooth scheme over $S$ with a $D(M)_S$-action where $M$ is an abelian group. Put $D= X _{S'}$. 
 Let $Y$ be a closed subscheme of $X_{S'}$ such that $Y \to S' $ is smooth. 
Let  $\Bl _Y^D X$ be the dilatation of $X$ with center $Y$ along $S'$. Let $N$ be a submonoid in $M$.
Assume that $Y$ is stable under the action of $D(M)_{S'}$ on $X_{S'}$. Then $D(M)_{S'}$ acts on $\Bl _{Y}^D X$ and we get a canonical isomorphism  
 \[ \Theta: (\Bl _Y^D X )^N \cong \Bl _{Y^N}^{D^N} X^N .\] 
 \xcoro

 \pf Since $X \to S ,D \to S'$ and $Y\to S$ are smooth,  by Corollary \ref{corosmoothxy} $X^N \to S, D^N \to S$ and $Y^N \to S$ are smooth. So by \cite[Proposition 2.16]{MRR20}, $\Bl_Y^DX \to S$  and $\Bl _{Y^N}^{D^N} X^N \to S$ are smooth. So using Corollary \ref{corosmoothxy} again, $(\Bl_Y^DX )^N \to S$ is smooth. Now since smooth implies flat, Corollary \ref{dilamagnecor} follows from Proposition \ref{dilattra}.
  \xpf

\rema  We note that the fact that dilatations commute with attractors may be used to study valued root data as in Bruhat-Tits theory. Indeed by Section \ref{sectionroot} root groups of reductive groups are examples of attractors (cf. also Section \ref{sectlie}) and dilatations allow to define filtrations.
\xrema

\section{Ind-algebraic spaces} \label{indspaces}
In this section, we study Algebraic Magnetism for Ind-algebraic spaces.

Let $S$ be a scheme and let $(Aff/S)$ be the category of affine schemes over $S$. Its objects are morphisms $\Spec(R) \to S$ from affine schemes to $S$. We use \cite{HR21} for the definition of ind-algebraic spaces.
\defi An ind-algebraic space (resp. ind-scheme) over $S$ is a functor $(Aff/S) \to Set$ which admits a presentation $X \cong \colim _{i \in I} X_i$ as a filtered colimit of $S$-algebraic spaces (resp. $S$-schemes) where all transition maps $\phi _{ij} :X_i \to X_j$, $i \leq j$ are closed immersions. The category of ind-algebraic spaces (resp. ind-schemes) over $S$ is the full subcategory of functors $(Aff/S ) \to Set$ whose objects are ind-algebraic spaces (resp. ind-schemes) over $S$.
\xdefi
Any algebraic space (resp. scheme) over $S$ is naturally an ind-algebraic space (resp an ind-scheme) over $S$. Any ind-scheme over $S$ is naturally an ind-algebraic space over $S$.
\rema \cite[§1.5]{HR21} If $X = \colim_i X_i$ and $Y = \colim_j Y_j$ are presentations of ind-algebraic spaces
(resp. ind-schemes) over $S$, and if each $X_i$ is quasi-compact, then as sets
$\Hom (X, Y ) = \lim _i \colim _j \Hom (X_i, Y_j ),$
because every map $X_i \to Y$ factors over some $Y_j$ by quasi-compactness of $X_i$. The categories of ind-algebraic spaces and ind-schemes are closed under fiber product. If $\mathbf{P}$ is a property of algebraic spaces (resp. schemes), then an $S$-ind-algebraic space
(resp. $S$-ind-scheme) $X$ is said to have ind-$\mathbf{P}$ if there exists a presentation $X = \colim_iX_i$ where
each $X_i$ has property $\mathbf{P}.$ A map $f : X \to Y$ of $S$-ind-algebraic spaces (resp. $S$-ind-schemes) is
said to have property $\mathbf{P}$ if $f$ is representable and for all schemes $T \to Y ,$ the pullback $f \times _Y T $ has
property $\mathbf{P}$. Note that every representable quasi-compact map of $S$-ind-schemes is schematic.
\xrema

\defi Let $G$ be a group scheme over $S$ and let $X$ be an ind-algebraic space over $S$. \begin{enumerate}
\item A categorical action of $G$ on $X$ is an action of $G$ on $X$ seen as a functor. Equivalently by Yoneda, a categorical action is a morphism in the category of ind-algebraic spaces over $S$ $\sigma :G \times _S X \to X$ satisfying the usual axioms, i.e.

\begin{enumerate}
\item $\sigma \circ ( \Id _G \times \sigma ) = \sigma \circ ( m \times \Id _X ) $ where $m: G \times _S G \to G$ is the group law of $G$,
\item $\sigma \circ (e \times \Id _X) = \Id _X $ where $e : S \to G $ is the identity section of $G$. 
\end{enumerate}

\item A collection of actions $\sigma _i : G \times _S X_i \to X_i $ such that for $i \leq j$ we have $ \phi _{ij} \circ \sigma _i = \sigma _j  \circ ( \Id _G \times \phi _{ij} )$ gives birth to a categorical action of $G$ on $X$. Such categorical actions are called ind-actions.
\end{enumerate}
\xdefi

Let $M$ be a finitely generated abelian group.

\prop \label{indprop}Let $X$ be an ind-algebraic space (resp. ind-scheme) over $S$. Assume that we have a presentation $X = \colim _i X_i $ with $X_i$ quasi-separated and locally finitely presented over $S$. Assume that $D(M)_S$ acts on $X$ via an ind-action on the presentation $X = \colim _i X_i $. Let $ N \subset M$ be a monoid. Assume that $N$ is finitely generated or that $X$ is separated. Then the attractor functor 
\[ (T \to S ) \mapsto \Hom ^{D(M)_T} ( A(N)_T , X_T) \] is representable by the ind-algebraic space (resp. ind-scheme) $\colim _i X_i ^N $. Moreover the natural morphism $X^N \to X$ is representable by algebraic spaces (resp. schemes). If $N=Z$ is a group the natural morphism $X^Z \to X$ is representable by a closed immersion. 
\xprop

\pf The first assertion is a direct corollary of Proposition \ref{opclosmoo}. To show that  $X^N \to X$, is representable by algebraic spaces, we notice that if $T$ is an affine scheme and $T \to X$ is
a morphism, then there exists $i$ such that this morphism is induced by a morphism $T \to X_i$, and
then we have
$ T\times _X X^N = T \times _{X_i} {X_i}^N.$ If $N=Z$ is a group, the last assertion follows from Corollary \ref{repgroup}.
\xpf

\section{Pure magnets}
\label{sectionefficient}

In this section, we study the behavior of attractors under variations of the underlying magnets. This leads to the notion of pure magnets and to a canonical invariant in Algebraic Magnetism.

Let $M$ be an abelian monoid. Let $X$ be a separated algebraic space over $S$ endowed with an action $a$ of $A(M)_S$. Let $m(a)$ be the set of magnets of $a$, i.e. the set of all submonoids of $M$. We also use the notation $m(M)$ to denote $m(a)$.

\prop \label{magnetequivlaence}Let $N \in m(a)$, the following conditions are equivalent.\begin{enumerate}
\item For any  $L\in m(a)$, $L  \subsetneq N \Rightarrow X^{L} \subsetneq  X^{N}$.
\item For any $L \in m(a)$, $X^N = X^L \Rightarrow N \subset L$.
\end{enumerate}
\xprop 

\pf Let us prove (i) $\Rightarrow $ (ii). So let $L \in m (M) $ and assume $X^N = X^L$. Then $X^N = X^{N \cap L}$ by Proposition \ref{NLNL}.
 Since $N \cap L \subset N $, (i) implies that $N \cap L = N$. So $N \subset L$.  Reciprocally assume (ii) holds. Let $L \in m(a)$ such that $L \subsetneq N$. By Proposition \ref{monosep}, we have $X^L \subset X^N$. It remains to prove that $X^L \ne X^N$.  This is clear, indeed otherwise (ii) implies that $N \subset L $ and so $N \subsetneq  N$ which is absurd.
\xpf 

\defi We use the following terminology.
 \begin{enumerate}
\item A pure magnet for the action $a$ is a magnet $N\in m(a)$ satisfying the equivalent properties of Proposition \ref{magnetequivlaence}.
\item The set of all pure magnets of the action $a$ is denoted $\mho (a)$.
\end{enumerate}
\xdefi

\theo \label{efficienttheo} We have a canonical bijection  between $\mho(a)$ and the set \[\{ X^N \subset X | N \in m (M) \}\] with the convention that we identify $X^N$ and $X^L$ if and only if $X^N (T)= X^L(T) \subset X(T)$ for all $T/S$ (cf. Proposition \ref{monoX}). The bijection sends a pure magnet $N \in \mho (a)$ to the attractor $X^N$. The reciprocal bijection sends $Y \in \{ X^N \subset X | N \in m (M) \}$ to the pure magnet \[E(Y):= \bigcap  \big\{ N \in m (M) | X^N = Y \big\}.\] 
\xtheo

\pf  
 Let $Y \in \{ X^N \subset X | N \in m (M) \}$ and let $E(Y) $ be the monoid defined in the statement, i.e. $E(Y)=  \bigcap  \big\{ N \in m (M) | X^N = Y \big\}$. By Proposition \ref{intersectioninfini}, we have $Y= X^{E(Y)}$. Let us prove that $E(Y)$ is a pure magnet. Let $L \in m (M) $ such that $X^L = X^{E(Y)} =Y $, by definition of $E(Y)$, we have $E(Y) \subset L$ and so $E(Y)$ is a pure magnet.  It remains to prove that for all $N \in \mho (a)$ and all $Y \in \{ X^N \subset X | N \in m (M) \}$, we have $E(X^N)=N$ and $Y= X^{E(Y)}$. So let us first take $N \in \mho (a) $. It is obvious that $E(X^N) \subset N$. On an other hand, by Proposition \ref{intersectioninfini}, we have $X^N = X^{E(X^N)}$ and so $ N \subset E(X^N)$ because $N$ is a pure magnet. Now let us take $Y \in  \{ X^N \subset X | N \in m (M) \}$, Proposition \ref{intersectioninfini} implies that $Y =X^{E(Y)}$. 
\xpf

\theo \label{strumagnet}Assume moreover that $X$ is finitely presented over $S$. Assume additionally that one of the following conditions is satisfied 
\begin{enumerate}
\item $X$ is affine over $S$,
\item $M$ is a group and there exists an $S$-affine strongly-FPR atlas for $X$ (e.g. the action is Zariski locally linearizable),
\end{enumerate}
then  $\mho (a) $ is a finite poset.
\xtheo

\pf The set $\mho (a)$ is a poset for the inclusion. It remains to show that $\mho(a)$ is finite. Assume (i) holds.
   Write $X= \Spec _S (\calA)$, then we have an $M$-grading $\calA = \bigoplus _{m \in M } \calA _m$. Because of $X$ is finitely presented over $S$, there exists a finite subset $E \subset M$ such that $\calA$ is generated locally by homogeneous elements of degree in $E$. Now for any monoid $N$, we have $X^N = X^{[N \cap E \rangle}$. This implies that $\mho (a)$ is finite. Case (ii) follows from case (i) and the following result:
\prop \label{lemmmmm}
Let $U  \to X$ be an $S$-affine strongly-FPR atlas. Let $a_U $ denote the action of $D(M)_S$ on $U$. Then for all submonoids $N,L $ of $M$, we have $U^N = U^L \Rightarrow X^N= X^L $; moreover  \[  \mho (a_X)\subset \mho (a_U) .\]
\xprop 
\pf Note that $X^P$ is representable for any magnet $P$ by Proposition \ref{representable}.
Let $N,L $ be submonoids of $M$ and assume $U^N=U^L$. By Proposition \ref{NLNL} or Proposition \ref{inter}, we have $U^N =U^{N \cap L} = U^{L}$. So it is enough to prove that $X ^{N \cap L }= X^{L}$, in other words we can change notation and assume $N \subset L$. Let us remark that $U \to S$ is of finite presentation, indeed $U \to S$ being affine, it is quasi-compact and quasi-separated and $U \to X \to S$ is locally of finite presentation as compositions of two such morphisms. The map $U^N \to X^N$ is étale and surjective. Indeed since $u$ is strongly-FPR, by Proposition  \ref{fixxx} we get a diagram
\[
\begin{tikzcd}
U^N \ar[r] \ar[d] & U^{N^*} \ar[r] \ar[d] & U \ar[d] \\
X^N \ar[r] & X^{N^*} \ar[r] & X
\end{tikzcd}
\] with Cartesian squares and $U \to X$ is étale and surjective. We get a commutative triangle \begin{tikzcd}U^{N} \ar[rr, "f"] \ar[rd, "p"] & &X^N \ar[dl,"q"]\\ & X^L &\end{tikzcd} where $f$ is surjective and étale, $p$ is étale by Corollary \ref{coroétt} and $q$ is locally of finite presentation (e.g. by \cite[\href{https://stacks.math.columbia.edu/tag/05WT}{Tag 05WT}]{stacks-project} and Theorem \ref{representable}). So $q$ is étale by  \cite[\href{https://stacks.math.columbia.edu/tag/0AHE}{Tag 0AHE}]{stacks-project} and in particular unramified. Now since $X$ is separated, $q$ is a monomorphism by Proposition \ref{monosep}. Being an unramified monomorphism, $q$ is universally injective by \cite[\href{https://stacks.math.columbia.edu/tag/05W6}{Tag 05W6}]{stacks-project}. So $q$ being étale and universally injective,  by  \cite[\href{https://stacks.math.columbia.edu/tag/05W5}{Tag 05W5}]{stacks-project}, it is an open immersion. Now since $p$ is surjective, we get that $q$ is surjective. So $q$ being a surjective open immersion, it is an isomorphism. So we proved that $ X^N=X^L$. Now let $N \in \mho (a_X)$  and let us prove that $N \in \mho (a_U)$. So let $L \in m (M)$ such that $U^N = U^L$. We proved that $X^N = X^L$ and so $L \subset N$ because $N \in \mho (a_X)$. This finishes the proof.
\xpf 

\xpf

\rema Theorem \ref{strumagnet} shows in particular that, under the same assumptions, the set $ \{ X^{D(M/Z)_S}\subset X \mid Z \subset M$ is a subgroup of $M\}$ is finite. We did not know a reference for this fact.  
\xrema

\conj \label{conj2}
We conjecture that the conclusion of Theorem \ref{strumagnet} remains true without the additional assumption (i),(ii) or (iii).
\xconj
\rema Conjecture \ref{conjecture} implies Conjecture \ref{conj2}. \xrema

 If $N$ is a monoid, we denote by $gen(N)$ the set of all subsets of $N$ that generate $N$ as monoid. The emptyset generates the zero monoid.
\defi The rank $mk(N)$ (possibly infinite) of a monoid $N$ is the cardinal defined by
\[
mk(N)= \min \{\#E | E \in gen(N)\}.
\]
\xdefi
Understanding combinatorial aspects of $\mho (a)$, including ranks of pure magnets, is a fundamental invariant of the action $a$. 

\exam \label{exammh} ${}$ \begin{enumerate} \item Let $a$ be the trivial action of $A(M)$ on a space $X$, then $\mho(a)=\{0\}$.
 \item Let $a$ be the action of $D(\bbZ^n)$ on itself by multiplication. Then $\mho(a)=\{0, \bbZ^n \}$. 
 \fact
We have $mk(\bbZ^n)=n+1$.
 \xfact  \pf Indeed, as monoid, $\bbZ^n$ is generated by $\{\{e_i\}_{i\in \{1 , \ldots , n \}},  - \sum _{i =1}^n e_i \}$.  Assume that $\bbZ^n$ is generated by $n$ elements $f_1, \ldots , f_n$, then there exists $p_1, \ldots , p_n >0$ such that $-(f_1+ \ldots +f_n) = p_1f_1+\ldots + p_n f_n$, this implies $(p_1+1) f_1+ \ldots + (p_n+1) f_n =0$. So $f_1, \ldots , f_n$ are linked and can not generate $\bbZ^n$ as group. This is absurd. So $\bbZ^n$ can not be generated by only $n$ elements.\xpf
 
 \item Let $a$ be the action of $A(\bbN)$ on itself by multiplication. Then $\mho(a)= \{0 , \bbN \}$,  note that $mk(\bbN )=1$.
 \item Let $a$ be the action of $D(\bbZ)$ on $A(\bbN)$ by multiplication. Then $\mho(a)= \{0, \bbN \}$.
 \item Let $a$ be the action of $D(\bbZ)$ on $A(\bbN)$ given by $\lambda\cdot x= \lambda ^2 x$. Then $\mho(a)= \{0, 2\bbN \}$, note that $mk(2\bbN )=1$.
 \item Let $a$ be the action of $D(\bbZ)$ on $\bbP^1$, then $\mho(a)= \{ 0 , \bbN , -\bbN , \bbZ \}$.
 \item Let $a$ be the action of $D(\bbZ^2)$ on $A (\bbN^2)$ by multiplication, then $\mho(a) = \{0\times 0,\bbN \times 0 , 0 \times \bbN , \bbN \times \bbN \}$.
 \item Let $a$ be the action of $D (\bbZ / 6 \bbZ)$ on $D(\bbZ /6\bbZ ) \times D (\bbZ /6 \bbZ ) \times D (\bbZ / 6 \bbZ ) $ given by \[\lambda \cdot (x,y,z)= (\lambda x , \lambda^2 y, \lambda ^3 z).\] Then $\mho (a) = \{ \bbZ / 6 \bbZ , 2 \bbZ / 6 \bbZ , 3 \bbZ /6 \bbZ , 6 \bbZ / 6 \bbZ  \}$.
 \item Let $a$ be the action of $D(\bbZ)= \bbG _m$ on $A(\bbN \setminus \{1\})$ where $\bbN \setminus \{1\}= [2,3\rangle$ is the submonoid of $(\bbN,+)$ generated by $2$ and $3$, then $\mho (a) = \{0,\bbN \setminus \{1\}, 2\bbN ,3\bbN\}$. This shows that Algebraic Magnetism refines strictly some aspects of the theory of $\bbG_m$-attractors. 
 \item Let $a$ be the adjoint action of a maximal split torus $T=D(M)$ on a given reductive group scheme $G$. Then $N \mapsto N \cap \Phi$ provides a bijection between $\mho(a) $ and the set of subsets of $\Phi (G,T)$ that are closed under addition. Pure magnets of rank $1$ correspond to roots (cf. §\ref{sectroots} for more details about the example of reductive groups).
 \item \cite[Proposition 3.1]{Ma25} Let $M$ be a cancellative monoid. Let $E$ be a subset of $M$. Let $[E\rangle$ be the submonoid of $M$ generated by $E$. Let $a$ be the action of $A(M)_S$ on $X=A([E\rangle)_S$ corresponding to the inclusion $[E\rangle \subset M$.
  Assume that \begin{enumerate} \item for all $E' \subsetneq E$, we have $[E'\rangle \subsetneq [E\rangle$ (i.e. $E'$ does not generate the monoid $[E\rangle$), and
\item $[E\rangle^*=0$,
\end{enumerate}
 then $\mho (a) \cong \{$subsets of $E\}$. More generally, cf. \cite{Ma25} for the pure magnets of self-actions of $A(L)_S$ with $L$ not necessarily sharp. 
\end{enumerate}
\xexam

\section{Complements, applications and examples} \label{examplesse}

This section collects several complementary results and examples, which are illustrations of the theory and have independent interest.

\subsection{Tangent spaces and attractors}

Let $S$ be a scheme. We denote by $I_S$ the scheme of dual numbers over $S$ as in \cite[Exp. II Définition 2.1]{SGA3}. For any scheme $T$ over $S$, we have $I_T= I_S \times _S T $.
 Explicitly, $I_S = \Spec (\mathbb{Z}[\varepsilon ]) \times _{\mathbb{Z}} S$ where ${\mathbb{Z}[\varepsilon] \overset{\varepsilon= [T]}{=} \mathbb{Z}[T]/(T^2)}$. Let $X$ be a functor over $S$. Let $T_{X}$ be the tangent space of $X$ as in \cite[Exp. II Définition 3.1]{SGA3}. This is a functor from $\Sch _S $ to $ Set$ sending a scheme $R$ over $S$ to 
 \[ \Hom _R ( I_R , X_R) .\]
 \rema \label{tangentfibrepro}${}$ \begin{enumerate} \item For $S$-morphisms $X \to Y$ and $Z \to Y$, we have $T_{X \times _Y Z } \cong T_X \times _{T_Y } T_Z $.
 \item If $X\to Y$ is a monomorphism, then $T_X \to T_Y $ is a monomorphism.
 \item We have a canonical identification $T_S=S$. \end{enumerate}
 \xrema
Let $M$ be an abelian group and assume that  $D(M)_S$ acts on $X$. Then $D(M)_S$ acts naturally on $T_{X}$ using the definition of $T_X$. Let $N$ be a submonoid in $M$.

\prop \label{tangent}We have a canonical isomorphism
\[ (T_{X})^N \cong T_{X^N} .\]
\xprop

\begin{proof} It is enough to show that $(T_{X})^N (S') \cong T_{X^N}(S')$ for any $S$-scheme $S'$. 
We have 
\begin{align*}
(T_{X})^N(S')& = \Hom _{S'} ^{D(M)_{S'}} ( A(N)_{S'} , T_{X_{S'}}) \\
& = \Hom _{S'}^{D(M)_{S'}}(A(N)_{S'} \times _{S'} I_{S'} , X _{S'} ) \text{ where } D(M)_{S'} \text{ acts trivially on } I_{S'}\\
&= \Hom _{I_{S'}}^{D(M)_{I_{S'}}} (A(N) _{I_{S'}} , X_{I_{S'}})\\ 
& = \Hom _{S'}  ( I_{S'} , {X_{S'}^N}) \\
&= (T_{X^N})(S').
\end{align*}
This shows that we have a canonical isomorphism.
\end{proof}

\subsection{Lie algebras and attractors} \label{sectlie}

Let $G/S$ be a group functor over a scheme $S$. Recall that in this case $T_{G}$ is a group functor over $S$ and we have two canonical morphisms of group functors $G \to T_{G } $ and $T_{G} \to G$ by \cite[Exp. II]{SGA3}. Let $e_S$ be the trivial group over $S$, as $S$-scheme we have $e_S=S$.  The Lie algebra of $G$ is defined as the fiber product \[
\xymatrix{ \Lie (G) = e_S \times _G T_G \ar[r]\ar[d] & T_G \ar[d] \\ e_S \ar[r] & G } \]
where $ e_S \xrightarrow{} G $ is the canonical morphism of group functors from $e_S$ to $G$. As in the previous section, let $M$ be an abelian group and assume that $D(M)_S $ acts on $G$. We assume moreover that this action is compatible with the group structure on $G$, i.e $D(M)_S$ acts by automorphisms on $G$. Then the induced action of $D(M)_S$ on $T_G$ is by group automorphisms. We thus obtain an action of $D(M)_S$ on $\Lie (G)$ by group automorphisms.

\rema \label{liefibrepro} ${}$\begin{enumerate}
\item For $S$-group functors $G,K,H$ with morphisms $G \to K, H \to K$, we have a canonical isomorphism $\Lie (G \times _K H ) \cong \Lie (G ) \times _{\Lie (K) }\Lie (H )$.
 \item If $G\to H$ is a monomorphism, then $\Lie (G ) \to \Lie (H) $ is a monomorphism.
 \item We have a canonical identification $\Lie ({e_S})=e_S$. \end{enumerate}
 \xrema

\prop \label{lieN}We have a canonical isomorphism of group functors over $S$
\[ \big(\Lie (G) \big)^N \cong \Lie (G^N) \]
\xprop
\pf Using Proposition \ref{produitNNN} and Proposition \ref{tangent}, we have
\[ \big(\Lie (G) \big)^N = \big( S \times _G T_G \big)^N \cong  S^N \times _{G^N} {T_G }^N \cong S \times _{G^N} T_{G^N} = \Lie (G^N).  \]
\xpf
Let us fix now a group functor $H$ over $S$ and a monomorphism $H \to G^{N^*} $ preserving the group structures. Recall that $G^N _H$ is the attractor with prescribed limit $H$ as in Definition \ref{definitionfonct}.

\prop \label{isoooooo} We have a canonical isomorphism of group functors over $S$
\[ \big(\Lie (G) \big)^N_{\Lie (H)} \cong \Lie (G_H^N) .\]
\xprop

\pf Using Proposition \ref{lieN} and Remark \ref{liefibrepro} we have
\begin{align*}
\Lie (G)^N _{\Lie (H)} &= \Lie (G) ^N \times _{ \Lie (G)^{N^*}} \Lie (H) \\
&\cong \Lie (G^N) \times _{\Lie (G^{N^*})} \Lie (H) \\
& \cong \Lie ( G^N \times _{G^{N^*} } H) \\
 &= \Lie ( G_H^N ) 
\end{align*}
This shows that we have the desired canonical isomorphism.
\xpf

\subsection{Semidirect products}

Let $G/S$ be an algebraic group in the sense of \cite{Mi17}. Here $S= \Spec (k)$ is a field. Let $M$ be a finite type abelian group and let $D(M)_S$ be the associated diagonalizable group scheme. Assume that $D(M)_S$ acts on $G$ by group automorphisms. Let $N \subset M$ be a submonoid. Recall that $G^N $ and $G^{N^*}$ are sub-algebraic-groups of $G$ over $S$ (e.g. combine \ref{representable}, \ref{monomono}, \ref{xyxyxyxyx} and \ref{groupmonoidstru}) and that we have a canonical morphism $G^N \to G^{N^*}$ of algebraic group over $S$ (cf. \ref{faceattractors}). The kernel of this morphism is $G^{N}_{e_G}$, the attractor with prescribed limit $e_G$ (cf. §\ref{prescribed}).

\prop \label{semidirect} The attractor $G^N$ is the semidirect product of  $G_{e_G}^{N} $ and $ G^{N^*}$. In particular, the multiplication map $G_{e_G}^{N} \rtimes G^{N^*} \to G^N$ is an isomorphism.
\xprop 
\pf
Note that $ G^{N^{*}}$ and $G^{N}_{e_G}$ are sub-algebraic-groups of $G^{N}$. There is a morphism $G^{N} \to G^{N^*}$ whose restriction to $G^{N^*}$ is the identity and whose kernel is $G^{N}_{e_G}$. Then Proposition \ref{semidirect} follows from \cite[2.34]{Mi17}.
\xpf 

\rema
Proposition \ref{semidirect} extends \cite[Theorem 13.33 (b)]{Mi17}.
\xrema 

\subsection{$\bbG _m$-attractors} \label{remarkclassicalvs} Let $X$ be an $S$-algebraic space with and action of $D(M)_S$. Let $N $ be a submonoid of the group $ M $. Then in some cases the attractor space $X^N$ can be obtained as a succession of attractors under $\bbG _m$ (cf. the introduction and \cite{Ri16}) and fixed-points. Let us gives two examples. 
\begin{enumerate} 
\item  Let $\alpha \in M$ and let us consider $X^{[\alpha\rangle}$. Under mild assumptions, by Proposition \ref{NLNL} we have $X^{[\alpha\rangle} = \big( X^{(\alpha)} \big) ^{[\alpha\rangle}$. By Proposition \ref{Group} $X^{(\alpha)} $ identifies with the fixed-points space $X^{D(M/(\alpha))_S}$. So the operation $X \rightsquigarrow X^{(\alpha )}$ can be realized as taking fixed-points. By Remarks \ref{actionmono} and \ref{actiongroup}, $D((\alpha))_S$ and $D(M)_S$ act on $X^{(\alpha)}$, and using Proposition \ref{propequiv}, we have $\big(X^{(\alpha)} \big) ^{[\alpha\rangle \subset (\alpha)} = \big(X^{(\alpha)} \big) ^{[\alpha\rangle \subset M}$ (cf. Remark \ref{Notationavecindicegroup} for the notation ${X}^{N \subset M}$). But $\big([\alpha\rangle \subset (\alpha) \big)\simeq \big(\mathbb{N} \subset \mathbb{Z}\big)$, so the operation $X^{(\alpha)} \rightsquigarrow  \big(X^{(\alpha)}\big)^{[\alpha\rangle} $ can be realized as taking attractor under $\mathbb{G}_m = D(\mathbb{Z}) _S$. So $X \rightsquigarrow X^{[\alpha\rangle}$ can be realized as fixed-points followed by taking the attractor under an action of $\bbG _m$.
\item Assume that $M= \mathbb{Z} \times \mathbb{Z}$ and that $N = \mathbb{N} \times \mathbb{N}$. Using Propositions \ref{NLNL} and \ref{propequiv}, we have
\[ X^{N}= \big(X^{\bbN \times \bbZ \subset \bbZ \times \bbZ  } \big) ^{\bbZ \times \bbN \subset \bbZ \times \bbZ } = \big(X^{\bbN \times 0 \subset \bbZ \times 0  } \big) ^{0 \times \bbN \subset 0 \times \bbZ } .\] This shows that $X \rightsquigarrow X^N$ can be realized as two stages of $\bbG _m$-attractors. \item If $M= \bbZ ^r $ and $N=\bbN ^r$ then $X \rightsquigarrow X^N$ can be realized as $r$ stages of $\bbG _m$-attractors.
\end{enumerate}

 \subsection{Magnetic point of view on reductive groups} \label{sectroots}
 This section is devoted to the observation that parabolic and Levi subgroups of reductive group schemes are easily described using attractors. We work with a split reductive group over a field for simplicity and accessibility but similar results hold more generally (e.g. cf. Proposition \ref{conradrootSGA} and \cite[§6.3]{ALRR22}). In fact, we expect attractors theory provide a natural framework to study some aspects of the advanced theory of group schemes from the beginning, but this is not the purpose of the present work.
So let $G$ be a split connected reductive group scheme over a field $R$. Let $T$ be a maximal split torus and choose a Borel $B$ containing $T$.
Let $\Phi =\Phi (G,T) \subset X^*(T)$ denote the set of roots associated to $(G,T)$ and $\Phi^+=\Phi (B,T)$ the roots in $B$. Let $\mathcal{B}$ be the basis of $\Phi$ determined by $\Phi ^+.$  For $\alpha \in \Phi $, let $U_{\alpha} $ be the associated unipotent root group and $\mathfrak{u} _{\alpha}$
 be the root group associated to $\alpha$ in the Lie algebra of $G$. We refer to \cite[Exp. XXII]{SGA3} for the definition of $U_{\alpha}$ and $\mathfrak{u} _{\alpha}$.
 Let $U$ be the unipotent radical of $B$ and let $\mathfrak{u} \subset \mathfrak{b} \subset \mathfrak{g}$ be the Lie algebras of $U,B$ and $G$.
 Consider the adjoint action of $T$ on $G,U$ and $\mathfrak{u}$.
 \prop \label{Uuequiv}
 There exists a $T$-equivariant isomorphism of $R$-schemes $\mathfrak{u} \simeq U$.
 \xprop
 \pf This is a direct consequence of \cite[Exp. XXII Th. 1.1, Exp. XXVI Prop. 1.12]{SGA3}, indeed these results imply the following assertions. For each root $\alpha \in \Phi ^+$, we have a $T$-equivariant isomorphism $U_{\alpha} \simeq \mathfrak{u} _{\alpha}$. 
 We have $T$-equivariant isomorphisms of schemes $\mathfrak{u} = \Pi _{\alpha \in \Phi ^+} \mathfrak{u}_{\alpha}$ and $U = \Pi _{\alpha \in \Phi ^+} U_{\alpha}$. This finishes the proof.  
 \xpf
 Let us now fix $\alpha \in \Phi ^+$. Recall that $[ \alpha \rangle \subset X^*(T)$ is the submonoid generated by $\alpha $.
 
 \prop \label{rootunipotent}
  We have a canonical isomorphism $\mathfrak{u}_{\alpha} \simeq \mathfrak{u}^{[ \alpha\rangle}$ (resp. $U_{\alpha} \simeq U^{[\alpha\rangle}$), between root group and $[\alpha\rangle $-attractor for the action of $T= \Spec (R[X^*(T)])$ on $\mathfrak{u}$ (resp. $U$).
 \xprop

 \pf  By Proposition \ref{Uuequiv}, we have a $T$-equivariant isomorphism $\mathfrak{u} \simeq U$, so it is enough to prove the statement for $\mathfrak{u}$. Since $G$ is split, $\Phi$ is reduced and Proposition \ref{diagonalizableattra} finishes the proof.
 \xpf
Recall that we have a bijection between parabolic subgroups of $G$ containing $B$ and subsets of $\mathcal{B}$, cf e.g. \cite[Page 35, lines 4-5]{Co14} or \cite[2.2.8]{CGP10}.
If $\Sigma$ is a subset of $X^*(T)$, we also use the notation $N_{\Sigma}$ to denote $[\Sigma \rangle$, the monoid generated by $\Sigma$ in $X^*(T)$. 
 
\prop \label{parabolicprop}  Let $\zeta \subset \mathcal{B}$.  Let $\Theta$ be $\zeta \cup -\zeta $. Let $\Sigma$ be $\mathcal{B} \cup - \zeta $.
\begin{enumerate} 
 \item  The attractor $G^{N_{\Theta}}$ is the Levi subgroup $L_{\Theta}$ such that $\Phi (L_{\Theta} , T) = N_{\Theta } \cap \Phi$.
 
 \item  The attractor $G^{N_{\Sigma}}$ is the associated parabolic subgroup, moreover $L_{\Theta}$ is a Levi component of $P$.
 
 \item Let $\xi \subset \zeta$. Let $\Gamma = \mathcal{B} \cup - \xi$. Let $N $ be a submonoid of $N_{\Sigma}$ such that $N \cap \Sigma = \Gamma $. Then $G^{N} = G^{N_{\Gamma  }}$.
 
\end{enumerate}

 \xprop
 
 \pf  Let $P_{\Sigma} $ be the parabolic corresponding to $\zeta $.
By \cite[2.2.8, 2.2.9]{CGP10}, there exists a $\lambda \in X_*(T)$ such that $P_{\Sigma}$ is the attractor associated to the monoid $\bbN$
relatively to the action of $\bbG _m=D(\bbZ)_S$ on $G$ via $x . g= \ad _{\lambda(x)}g$ and such that $\lambda ( \beta ) \geq 0 $ for all $\beta \in \Sigma$ and $\lambda (\beta) =0 $ for all $\beta \in \Theta$.
The Levi subgroup $L_{\Theta}$ corresponding to $\Theta $ is the fixed space in $G$ under the action of $\lambda $ by conjuguation, i.e $L_{\Theta}=G^0$. 
Now we prove the Proposition. \begin{enumerate}
\item Assume first that $\zeta = \mathcal{B}$. Then $L_{\Theta} = G$ and $N_{\Theta}= N_{\Phi (G,T)}$.
Using Prop. \ref{rootunipotent}, we deduce that the big cell
$\Omega = \Pi _{\alpha \in \Phi ^-} U_{\alpha} \times T \times \Pi _{\alpha \in \Phi ^+} U_{\alpha}$ is in $G^{N_{\Theta}}$. Now we have inclusions $\Omega \subset G^{N_{\Theta}} \subset G$ with $\Omega $ dense in $G$ and $G^{N_{\Theta}}$ closed in $G$. This implies $G^{N_{\Theta}} = G$.
Let us now prove the general case, let $\zeta \subset \mathcal{B}$. 
By Proposition \ref{propequiv}, we have $G^{f^{-1}(0) }= G^0=L_{\Theta}$. We have $\Theta \subset f^{-1} (0)$, and so $N_{\Theta } \subset f^{-1} (0)$,
consequently $G^{N_{\Theta}} \subset G^{f^{-1}(0)}.$ So we have proved that $G^{N_{\Theta}} \subset L$ and let us now prove that this is an equality. We remark that $N_{\Theta} = N_{\Phi (L_{\Theta},T)}$.
Now since $L \subset G $ and using the first case done before, we have $L=L^{N_{\Phi (L_{\Theta},T)}} \subset G ^N_{\Theta}$. This finishes the proof.
      \item 
Recall that $\lambda: \bbG _m \to T$ corresponds to the morphism of abelian groups $f:X^*(T) \to \bbZ , \chi \mapsto (\lambda ,\chi)$. 
Now we see $G$ as a $\bbG_m$-scheme and as a $T$-scheme. By Proposition \ref{propequiv}, we have $ G^{\bbN} = G^{f^{-1}(\bbN)}$. 
Since, for all $\beta \in \Sigma$,  $f( \beta) =\lambda (\beta ) \geq 0$, we have $f(\beta ) \in \bbN$ and so $\Sigma \subset f^{-1} (\bbN)$, and so $N_{\Sigma} \subset f^{-1}(\bbN)$.
Consequently, $P_{\Sigma} = G^{\bbN} = G^{f^{-1}(\bbN)} \supset G^{N_{\Sigma}}$. Let us prove that
 $P_{\Sigma} \subset G^{N_{\Sigma}}$.  
 We have $P_{\Sigma} = L_{\Theta} \times \mathrm{R _u} (P_{\Sigma})$ where $\mathrm{R _u} (P_{\Sigma}) $ is the unipotent radical of $P_{\Sigma}$. 
 Using (i), we have $L_{\Theta} = L_{\Theta} ^{N_{\Theta}} \subset L_{\Theta} ^{N_{\Sigma}} \subset L_{\Theta} $, and so $L_{\Theta}^{N_{\Sigma}} = L_{\Theta}$. Using Proposition \ref{rootunipotent}, one has $( \mathrm{R _u} (P_{\Sigma}))^{N_{\Sigma}} =  \mathrm{R _u} (P_{\Sigma})$.
 So $P_{\Sigma}^{N_{\Sigma}}= P_{\Sigma}$, and so $P_{\Sigma} \subset G^{N_{\Sigma}}$. 
\item We have a canonical closed immersion $G^{N_{\Gamma}} \subset G^N$.  We have $N \cap N_{\Theta}= N_{\zeta \cup - \xi}$ and $N \cap N_{\Theta}$ it is thus included in $N_\Gamma$. By Proposition \ref{inter} we have $G^N= (G^{N_\Sigma})^N$. 
By the previous assertions we have $ G^{N_{\Sigma}} = G^{N_{\Theta}} \times \mathrm{R _u} (G^{N_{\Sigma}}) $. Now by \ref{inter} we have \[ \Big( G^{N_{\Theta}} \times \mathrm{R _u} (G^{N_{\Sigma}})  \Big) ^N = ( G^{N_{\Theta}})^N \times (\mathrm{R _u} (G^{N_{\Sigma}})  )^N= G^{N_{\zeta \cup - \xi}}\times \mathrm{R _u} (G^{N_{\Sigma}})    .\] This implies that $G^{N_{\Gamma}} \supset G^N$ and finishes the proof.
\end{enumerate}  \xpf

 \rema \label{cartesian-para}Proposition \ref{parabolicprop} implies that any parabolic or Levi subgroup of $G$  containing $T$ can be obtained as an attractor under the conjugation action of $T$ on $G$. Moreover, assume that $\mathbf{B} $ is a parabolic subgroup in a parabolic $\mathbf{P}$ and $\mathbf{L}$ is a Levi component of $\mathbf{P}$. We assume that $\mathbf{B}, \mathbf{P}$ and $\mathbf{L}$ contain $T$. Let $\mathbf{M}$ be $ \mathbf{B} \cap \mathbf{L}$, this is a parabolic subgroup in $L$. Then one has a cartesian square \[\begin{tikzcd} & \mathbf{B} \ar[rd] \ar[dl] & \\ \mathbf{P} \ar[rd] & & \mathbf{M} \ar[dl] \\& \mathbf{L} & \\ \end{tikzcd}. \]
 This square can be obtained using Proposition \ref{cartesian}. Indeed let $L' $ be the submonoid generated by $ {\Phi (\mathbf{P},T)}$, let $L$ be the submonoid generated by ${\Phi ( \mathbf{L},T)} $ and $N $ be the submonoid generated by ${\Phi ( \mathbf{M},T)}$. Using Proposition \ref{abelianimplynice}, we deduce that $L $ is a face of $ L'$. Now let $N' $ be $L' \setminus ( L \setminus N ) $. We have $G^{L'}= \mathbf{P}, G^{N}= \mathbf{M}$ and $G^{L}= \mathbf{L}$. Using Proposition \ref{parabolicprop} we have $G^{N'}= \mathbf{B}$.
 \xrema
For any root $\alpha \in \Phi^+$ we denote by $H_{\alpha} \subset G$ the semidirect product $T \ltimes U_{\alpha}$, this is a group scheme whose unipotent radical equals $U_{\alpha}$.
 
 \prop \label{root} We have a canonical isomorphism $H_{\alpha}\simeq G^{[\alpha \rangle}$.
 \xprop 
\pf Since $[ {\alpha }\rangle \subset \Sigma _{\Phi^+ }$ and by Proposition \ref{parabolicprop},  we have a closed immersion $G^{[ \alpha \rangle}\subset B=G^{\Sigma _{\Phi^+ }}$ where $B$ is the Borel subgroup.
So by Lemma \ref{equiclosed}, we get a closed immersion  $G^{[ \alpha\rangle}\subset B^{[\alpha \rangle}$ and thus an equality $G^{[\alpha \rangle}= B^{[ \alpha \rangle}$.
Now we have a $T$-equivariant isomorphism of schemes $B \simeq T \times U $. 
Using Proposition \ref{rootunipotent}, we get \begin{equation} \label{proofrooteq2}  U _{\alpha} = U^{[ \alpha \rangle}.\end{equation}
It is obvious that \begin{equation} \label{proofrooteq3} T=  T^{[ \alpha \rangle } .\end{equation}
Now equations $ \ref{proofrooteq2}$ and $\ref{proofrooteq3}$ and Proposition \ref{produitNNN} (ii) imply that $B^{[ \alpha \rangle} = H_{\alpha}$. This finishes the proof.
\xpf

\coro  \label{prescriroots} ${}$\begin{enumerate} \item  Let $e_G $ be the closed subscheme of $G$ corresponding to the unit section. Then the attractors $G_{e_G}^{[\alpha\rangle}$ with prescribed limit $e_G$ equals $U_{\alpha}$.

 \item We have a canonical isomorphism $U_{\alpha} = \mathrm{R_u}( G^{[ \alpha \rangle} )$ where $\mathrm{R_u}$ means the unipotent radical.
\end{enumerate}

\xcoro

\begin{proof} Clear from Proposition \ref{root} and Section \ref{prescribed}.
\end{proof}

\subsection{Root groups} \label{sectionroot}

Let $S$ be a scheme, let $G$ be a group scheme over $S$. Let $M$ be an abelian group. Assume that $D(M)_S $ acts on $G$ by group automorphisms. Then we get an action of $D(M)_S$ on $\Lie (G)_S$ by group automorphisms. Let $N$ be a submonoid of $M$. Let $G^N_{e_S}$ be the attractor associated to $N$ with prescribed limit the trivial subgroup $e_S/S$ relatively to the face $N^*$.

\defi Let $\alpha \in M$.  We call $ G^{[\alpha\rangle}_{e_S}$ the root group associated to $\alpha$ under the action of $D(M)_S $ on $G$. We call $G^{[\alpha\rangle}$ the non prescribed root group associated to $\alpha$.
\xdefi

\prop Let $\alpha \in M$. \begin{enumerate}

\item If $G$ is affine over $S$, we have closed immersions of group schemes over $S$\[ G_{e_S}^{[\alpha\rangle} \subset G^{[\alpha\rangle} \subset G .\]
\item If $G/S$ is smooth, then $G^{[\alpha\rangle}_{e_S} $ and $G^{[\alpha \rangle}$ are smooth over $S$.
\item We have canonical isomorphisms $\Lie (G^{[\alpha\rangle}_{e_S}) = \big( \Lie (G) \big) ^{[\alpha\rangle}_{e_S} $ and $\Lie (G^{[\alpha \rangle}\rangle=\Lie (G) ^{[\alpha\rangle} . $
\end{enumerate}
\xprop
\pf \begin{enumerate} \item We know that $G^{[\alpha\rangle}$ is a closed subgroup of $G$ by Theorem \ref{representableaffine}. We have $ G^{[\alpha\rangle}_{e_S} = G^{[\alpha\rangle} \times _{G^0} e_S$ and so $ G^{[\alpha\rangle}_{e_S}$ is a closed subgroup of $G^{[\alpha\rangle}$. 
\item By Corollary \ref{123456789}, $G^{[\alpha\rangle} \to S$ and $G^{[\alpha\rangle} \to G^0$ are smooth, this implies the claim.
\item These are special cases of Propositions \ref{isoooooo} and \ref{lieN}. \end{enumerate}
\xpf
We observe that our definition is compatible with the definition given by Conrad and SGA3.
\prop \label{conradrootSGA}Let $G \to S$ be a reductive group scheme over a non-empty scheme $S$, $ T \cong D_S(M)$ a split maximal torus, and $\alpha \in M $ a root. Let $\mathrm{exp}_\alpha (\mathbf{W} (\mathfrak{g} _\alpha )) \subset G$ be the $\alpha$-$root$ $group$ for $(G,T,M)$ considered by B. Conrad in \cite[Theorem 4.1.4]{Co14}, then 
\[ \mathrm{exp}_\alpha (\mathbf{W} (\mathfrak{g} _\alpha )) = G^{[\alpha\rangle}_{e_S}\]
\[ \mathrm{exp}_\alpha (\mathbf{W} (\mathfrak{g} _\alpha )) \times _S T = G^{[\alpha\rangle}.\]
\xprop

\pf We remark that we have $( \mathbf{W} (\mathfrak{g} _\alpha ) \times _S T)^{[\alpha\rangle}= \mathbf{W} (\mathfrak{g} _\alpha ) \times _S T$. By \cite[Theorem 4.1.4]{Co14} we have a closed immersion $\mathbf{W} (\mathfrak{g} _\alpha ) \times _S T\to G$. So Lemma \ref{equiclosed} gives us a closed immersion ${f: \mathbf{W} (\mathfrak{g} _\alpha ) \times _S T \to G^{[\alpha\rangle}}$. We now prove that $f$ is an open immersion. Since $f$ is of finite presentation by \cite[\href{https://stacks.math.columbia.edu/tag/02FV}{Tag 02FV}]{stacks-project}, using \cite[Proposition 14.18]{GW} it is enough to prove that $f$ is flat. Since $\mathbf{W} (\mathfrak{g} _\alpha ) \times _S T   \to S$ is flat, it is enough to prove that for any $s \in S $, the morphism $f_s :\big( \mathbf{W} (\mathfrak{g} _\alpha ) \times _S T \big)_s \to \big( G^{[\alpha\rangle} \big)_s $ is flat (cf. \cite[Proposition 14.25]{GW}). Since all involved constructions are compatible by base change, $f_s$ is an isomorphism (in particular flat) for any $s \in S$ by Corollary \ref{prescriroots}. Now $f$ is an open and closed immersion that gives isomorphisms on fibers. This implies that $f$ is an isomorphism and proves the second equality. Now the first equality is clear because $[\alpha\rangle^*=0$, $G^0=T$ and so  \[G^{[\alpha\rangle} _{e_S} = G^{[\alpha\rangle}\times _{G^0} e_S = \big( \mathrm{exp}_\alpha (\mathbf{W} (\mathfrak{g} _\alpha )) \times _S T \big) \times _{T} e_S = \mathrm{exp}_\alpha (\mathbf{W} (\mathfrak{g} _\alpha )).\]
\xpf

\rema Let us remark that, 
$G^{[\alpha\rangle}=\big(G^{D(M/(\alpha))_S} \big)^{[\alpha\rangle \subset (\alpha)}$ (e.g. by §\ref{remarkclassicalvs}). 
So $G^{[\alpha\rangle}$ can be obtained as a first stage of fixed-points followed by a stage of $\bbG _m$-attractor. 
Similarly $G^{[\alpha\rangle}_{e_S}=\big( G^{(\alpha)}_{T_{\alpha}} \big)^{[\alpha\rangle}_{e_S}$ where $T_{\alpha}$ is defined in \cite[Lemma 4.1.3]{Co14} or \cite{SGA3}, this explains why \cite{Co14} works with small semisimple groups of rank one as a first stage in order to build root groups using $\bbG _m $-attractors. Note that \cite{SGA3} also uses semisimple groups of rank one as a first stage \cite[Exp. XX]{SGA3} before defining root groups in the general case \cite[Exp. XXII]{SGA3}.
\xrema

\subsection{Monoschemes and toric schemes} Let $M$ be an abelian group and let $N$ and $L$ be submonoids of $M$. Let us consider the attractors $A(N)^L$ associated to the monoid $L$ under the action of $D(M)$ on $A(N)$. By Theorem \ref{representableaffine} $A(N)^L$ equals $\Spec \big( \bbZ [N] / (X^n | n \in N \setminus (N \cap L) ) \big)$. The ideal $ \calI := (X^n | n \in N \setminus (N \cap L) )$ of $\bbZ [N]$ equals $\bigoplus _{i \in I} \bbZ X^i$ where $I$ is the ideal of $N$ generated by $N \setminus (N \cap L)$. If $I$ is a prime ideal of $N$ (cf. \cite{Og}) then $N \setminus I$ is a submonoid of $N$ (and necessarily a face) and $\bbZ [N]/ \calI  = \bbZ [N \setminus I]$. In this case, $A(N)^L$ is also a scheme associated to a submonoid of $M$. For example if $N \cap L$ is a face of $N$, then $N \setminus I = N \cap L $ and $A(N)^L = A(N \cap L )$. In general $I$ is not a prime ideal and so $N \setminus I $ is not a submonoid of $N$ (e.g. take $M = \bbZ ^2, N= [(1,1),(1,-1),(1,0) \rangle $ and $L = [(1,0)\rangle $, then $N \setminus I = \{ (0,0) , (1,0) \}$).

More generally, let $\calN$ be a toric monoscheme whose associated finitely  generated abelian group $\Gamma$ is $M$ (cf. \cite[II §1.9]{Og}). Let $A(\calN)$ be the scheme asssociated to $\calN$ (cf. \cite[II Prop. 1.9.1]{Og}), this is a toric scheme. Let $\{ \mathrm{spec} (N_\tau) \}_{\tau \in \calA}$ be an open affine covering of $\calN$ ($N_\tau \subset M$ for all $\tau$). Then $\{ A(N_\tau ) \}_{\tau \in \calA}$ is an open affine covering of $A(\calN)$.  Then $\coprod _\tau A(N_\tau) \to A(\calN)$ gives an FPR atlas of $A( \calN)$ (recall that a $D(M)$-equivariant open immersion is $Z$-FPR for any subgroup $Z \subset M$). Let $L$ be a submonoid of $M$, then we obtain that $\{ A(N _{\tau} )^L \}$ is an affine open covering of $A(\calN)^L$. 

\ackn   A. Mayeux thanks A.-M. Aubert, M. Brion, S. Brochard, S. Cotner, V. Dobrev, A. Dubouloz, P. Gille, T. Haines, J. Lourenço, T. Richarz, S. Riche, J. P. dos Santos and the anonymous referee for useful feedback on this work and/or for the opportunity to present it at seminars in Saint-Étienne (2022), Montpellier (2022), Poitiers (2023), and Grenoble (2023); at the XV International Workshop “Lie Theory and Its Applications in Physics” (2023); at the Oberwolfach Workshop “Algebraic Groups” (2025); and in the “Moduli Spaces in Algebraic Geometry” section at AAG 2025. This work has received funding from ERC grant No 101002592 (Laboratoire de Mathématiques Blaise Pascal, Université Clermont Auvergne) and ISF grant 1577/23 (Einstein Institute of Mathematics, The Hebrew University of Jerusalem).
\xackn

\appendix  
\section{Fixed-point-reflecting morphisms, after Drinfeld and Alper-Hall-Rydh. \\by Matthieu Romagny}
\label{appendixx} 

Let $U\to X$ be an equivariant morphism of $S$-algebraic
spaces endowed with actions of an affine $S$-group scheme $G$
whose function algebra is free over $\calO_S$ (e.g. a
diagonalizable $S$-group scheme). It is classical that
if $U/S$ is separated
then the fixed point functor is a closed subspace $U^G\to U$
(\cite[Exp.~VIII, \S~6]{SGA3}). By \'etale descent, if
$U\to X$ is \'etale, surjective and reflects the fixed points,
then the same conclusion follows for $X^G\to X$.
This need not always hold: if $X$ is the affine
line with doubled origin over the field of complex numbers,
and $V,W\subset X$ are the two glued copies of the affine line,
then the action of $\bbZ/2\bbZ$ that permutes~$V$ and $W$
induces an action on $X$ whose fixed point scheme is the
complement of the two origins.

Thus there is a close relationship between the existence of
FPR atlases and the closedness of fixed points. In this appendix
we present two situations where these properties occur.
The first is a useful generalization of an argument of Drinfeld, who  
in \cite[Prop. 1.2.2]{Dr15} considers the case where $S$ is the spectrum of a  
field and $G$ is the multiplicative group $\bbG_m$ (cf. also \cite[Lemma 1.10]{Ri16} and Propositions \ref{repgroup} and \ref{lemmpf}).

\theo \label{G_connected_implies_fixed_pts_closed}
Let $G$ be a flat, finitely presented $S$-group scheme with
{\em connected fibres}. Let $X,Y$ be locally finitely presented
$S$-algebraic spaces endowed with $G$-actions. Assume either
\begin{itemize}
\item[\rm (i)] $S$ is locally noetherian, or
\item[\rm (ii)] $X,Y$ are quasi-separated.
\end{itemize}
Then $X^G\to X$, $Y^G\to Y$ are closed immersions of finite
presentation and all \'etale equivariant morphisms $X\to Y$
are fixed point reflecting.
\xtheo

Interestingly this gives a case where closedness of
$X^G\to X$ is ensured by an assumption on the group
(connectedness) rather than on the space (separation).
The second result, Theorem \ref{adapted-atlas-app}, is essentially a corollary
of results of Alper, Hall and Rydh \cite{AHR21} and the present version was formulated by Mayeux (cf. Theorem \ref{adapted-atlas}).

\theo \label{adapted-atlas-app}
Let $X$ be a quasi-separated $S$-algebraic space locally of
finite presentation endowed with an action of a finitely
presented diagonalizable group scheme $G=D(M)_S$.
Let $H=D(M/Z)$ be a subgroup scheme and assume that one of the
following assertions holds:
\begin{enumerate}
\item $X$ is separated over $S$,
\item $H$ has connected fibres.
\end{enumerate}
Then there exists a $Z$-FPR atlas $U \to X$ (cf. Definition \ref{deffpratl}), which may be chosen
quasi-compact if $X\to S$ is.
\xtheo
Let us get the proofs started. We start with Drinfeld's
result on the closedness of fixed points. We shall prepare
the discussion with three preliminary lemmas.
Let $F=(F,e:S\to F)$ be a pointed $S$-algebraic space.
An {\em action} of $F$ on~$X$ is a morphism of $S$-algebraic
spaces $F\times X\to X$, $(f,x)\mapsto fx$ through which the
section~$e$ acts trivially. Examples are the projection
$\pr_2:F\times X\to X$, called the {\em trivial action}, and
the action induced when $F$ is a subscheme of an $S$-group
scheme~$G$ acting on $X$ in the usual sense.
To any action is associated
its {\em stabilizer}, the pointed sub-$X$-algebraic space of
$F\times X$ defined by pulling back the map
$F\times X\to X\times X$, $(f,x)\mapsto (x,fx)$ along the diagonal:
\[
\begin{tikzcd}
\Stab \ar[r,"\varphi"] \ar[d] & F\times X \ar[d] \\
X \ar[r,"\Delta"] & X\times X.
\end{tikzcd}
\]

\lemm \label{fixed_points_ff_subscheme}
Let $F=(F,e)$ be a finite locally free, infinitesimal pointed
$S$-algebraic space acting on the $S$-algebraic space $X$. Then
the functor of $F$-fixed points, whose values over an $S$-algebraic
space $T$ are the $F$-equivariant maps $u:T\to X$ where $T$
is endowed with the trivial action, is representable by a closed
subspace $X^F\subset X$.
\xlemm

\pf
Note that $X^F$ is the Weil restriction of
$\varphi:\Stab\to F\times X$ along the projection
$\pr_2:F\times X\to X$.
In particular, its representability by an algebraic space is
standard, for all $F$ finite locally free (see e.g.
\cite[Th.~3.7]{Ry11}). To prove the closed immersion property
in the infinitesimal case, we argue as follows. Because the
stabilizer is sandwiched as
\[
X\hookto \Stab \hookto F\times X\]
between two finite infinitesimal $X$-spaces, then also
$\Stab\to X$ is finite infinitesimal and~$\varphi$
is a closed immersion. Let $\calA$
be the $\calO_X$-algebra of functions of $F\times X$ and
$\calI\subset\calA$ the defining ideal sheaf of $\varphi$.
Let $u:T\to X$ be a morphism. Saying that $F$ acts trivially on
$T$ means that~$\varphi$ restricts to an
isomorphism above $F\times T$. The latter assertion means that
after the base change $u^\sharp:\calO_X\to u_*\calO_T$, the map
$\calI\to\calA$ becomes the zero map. The $\calO_X$-algebra
$\calA$ being finite locally free, the formation of its linear dual
commutes with base change; hence it is equivalent to say that
$\calI\otimes\calA^\vee\to\calO_X$ becomes the zero
map. Equivalently the $\calO_X$-ideal
$\calJ:=\text{im}(\calI\otimes\calA^\vee\to\calO_X)$
is contained in $\ker(u^\sharp)$. This shows that the functor of
fixed points is representable by the closed subscheme
$V(\calJ)\subset X$.
\xpf

\lemm \label{lemm_cofinal_thickenings_1}
Let $R\to A$ be a ring map and $x_1,\dots,x_d$ elements of $A$.
Let $I_n$ be the ideal of $A$ generated by the powers
$(x_1)^n,\dots,(x_d)^n$, for each $n\ge 1$. If $A/I_1$ is finite
over $R$, then $A/I_n$ is finite over $R$ for all $n\ge 1$.
\xlemm

\pf
Let $(I_1)^n$ be the powers of the ideal $I_1$. By induction,
using that $(I_1)^n/(I_1)^{n+1}$ is a finite $A/I_1$-module hence
finite over $R$, we see that $A/(I_1)^n$ is finite over $R$
for all $n$. The containments
$(I_1)^{d(n-1)+1}\subset I_n\subset (I_1)^n$ show that the
sequences of ideals $\{I_n\}_{n\ge 1}$ and $\{(I_1)^n\}_{n\ge 1}$
are cofinal, and the claim follows.
\xpf

\lemm \label{lemm_cofinal_thickenings_2}
Let $G$ be a flat, finitely presented $S$-group scheme with
connected fibres. Then for each point $s\in S$ there is an
\'etale neighbourhood $S'\to S$ and a sequence of finite
locally free infinitesimal neighbourhoods of the
unit section of $G\times_SS'$ which is cofinal to the canonical
sequence of $n$-th order thickenings.
\xlemm

\pf
Let $e=e(s)$ be the unit section
of the fibre $G_s$. Since the local ring $\calO_{G_s,e}$ is
Cohen-Macaulay (\cite[Exp.~VII$_B$, Corollaire 5.5.1]{SGA3}),
it admits a regular sequence $\bar x=(\bar x_1,\dots,\bar x_d)$
of length $d=\dim(\calO_{G_s,e})$. Let~$U$ be an open subscheme
of $G$ over which the germs $\bar x_i$ extend to local functions
$x_i$ belonging to the augmentation ideal
$\ker(e^\sharp:\calO_U\to\calO_S)$.
Because $G\to S$ is flat and the sequence $\bar x$ is
regular, the closed subscheme $F\subset U$ cut out by
$x_1,\dots,x_d$ is {\em flat} and
quasi-finite over $S$ at the point $e$. Thus shrinking $U$ is
necessary, we may assume that $F\to S$ is flat and quasi-finite;
it is furthermore finitely presented. By
\cite[IV, Th.~18.12.1]{EGA} there is an étale extension
$S'\to S$ and an open neighbourhood $V'$ of $e'=e(s')$ in
$G\times_S S'$ such that $F':=(F\times_SS')\cap V'$
is {\em finite} over $S'$, and therefore locally free. For
each $n\ge 1$, the sequence $(x_1)^n,\dots,(x_d)^n$ is again
regular and it now follows from
Lemma~\ref{lemm_cofinal_thickenings_1} that its vanishing
locus $F'_n$ in $V'$ is also finite locally free over $S'$.
\xpf

\rema
If $G\to S$ is smooth, the $n$-th order thickenings of the
unit section are finite locally free, hence they satisfy the required property.
\xrema

We are ready to prove that ``$X^G\to X$ is closed when $G$
is connected''.

\begin{proof}[Proof of Theorem~\ref{G_connected_implies_fixed_pts_closed}]
We start with the statement about $X^G\to X$.
The assumptions and conclusion being local over $S$, we may
assume that~$S$ is affine. Because $G\to S$ is flat and
finitely presented, the orbit in $X$ of an open quasi-compact
subspace $W\subset X$ is open and quasi-compact. We may
replace~$X$ by one such orbit and hence assume that~$X$ is
quasi-compact. In the case when $X$ is assumed quasi-separated,
it is then of finite presentation. By standard results on
limits (as in \cite[\href{https://stacks.math.columbia.edu/tag/07SJ}{Tag 07SJ}]{stacks-project}), the space $X$ and
the $G$-action then come from a finitely presented algebraic
space $X_0\to S_0$ with $G_0$-action by a base
change $S\to S_0$ where $S_0$ is of finite type over~$\bbZ$.
Thus in all cases we may assume that $S$ is noetherian
and $X$ is locally noetherian.

The claims about $X^G\to X$ are \'etale-local over $S$
so using
Lemma~\ref{lemm_cofinal_thickenings_2} we can assume that
there exists a sequence $\{F_n\}$ of finite locally free
infinitesimal neighbourhoods of the unit section of $G$ which
is cofinal to the canonical sequence of $n$-th order thickenings.
By Lemma~\ref{fixed_points_ff_subscheme} the fixed points
$X^{F_n}$ are closed subspaces of $X$. Let $X_0$ be the closed
subspace of $X$ equal to their intersection. From the obvious
inclusions $X^G\subset X^{F_n}$ follows that $X^G\subset X_0$.
In order to prove the opposite inclusion and conclude, it is
enough to prove that the top map $\varphi_0$ in the pullback
diagram
\[
\begin{tikzcd}
\Stab_0 \ar[r,"\varphi_0"] \ar[d] & G\times X_0 \ar[d] \\
\Stab \ar[r,"\varphi"] & G\times X
\end{tikzcd}
\]
is an isomorphism. We view $\varphi_0$ as a map of
$X_0$-group schemes. Note that $\varphi_0$ is of finite
presentation. We claim that $\varphi_0$ is formally étale
along the unit section. Since the target is locally
noetherian, argueing as in \cite[IV, Prop.~17.14.2]{EGA}
it is enough to prove that each diagram
\[
\begin{tikzcd}
\Stab_0 \ar[r,"\varphi_0"] & G\times X_0 \\
\Spec(A) \ar[r] \ar[u] \ar[ru,dashed] & \Spec(A') \ar[u]
\end{tikzcd}
\]
can be filled as indicated, where $\Spec(A)\to\Spec(A')$ is
a square-zero thickening of artinian local schemes.
Base-changing along $\Spec(A')\to X_0$ we reduce to the case
where $X_0$ is local artinian, in which case the group
schemes involved are separated
(\cite[Exp.~VI$_A$, \S~0.3]{SGA3}) and
$\Spec(A')\to G\times X_0$ is a closed immersion.
Since $\Spec(A')$ is artinian, it is included in one of
the $F_n\times X_0$, and then the lifting result follows
from the definition of $X_0$.

This proves that the maximal open subscheme $U\subset \Stab_0$
where $\varphi_0$ is étale contains the unit section.
Clearly $U$ is stable by multiplication by local sections
of $U$, and also by inversion. Hence~$U$ is an open
subgroup scheme of $\Stab_0$ and also (because
${\varphi_0}_{|U}$ is étale) of $G\times X_0$. Since
$G\to S$ has connected fibres, its only open subgroup
scheme is itself, whence $U=G=\Stab_0$ and the conclusion.

We now consider the statement about \'etale equivariant
morphisms $X\to Y$. Proceeding as before we reduce to the
situation where we have a family $\{F_n\}$ with equalities
$X^G=\cap_{n\ge 0} X^{F_n}$ and $Y^G=\cap_{n\ge 0} Y^{F_n}$.
Since the
intersections commute with pullback along $X\to Y$,
it is enough to prove that the natural
map $X^{F_n}\to X\times_Y Y^{F_n}$ is an isomorphism.
For this it is enough to prove that the action of $F_n$
on the subspace $X_0:=X\times_Y Y^{F_n}$ of $X$ is trivial.
Consider the diagram:
\[
\begin{tikzcd}[column sep=40]
X_0 \ar[r,"{e\times\id}"] \ar[d]
& F_n\times X_0 \ar[d] \ar[ld,dashed] \\
X \ar[r,"\Delta"'] & X\times_Y X.
\end{tikzcd}
\]
Here the right vertical map is $(f,x)\mapsto (x,fx)$
and the bottom arrow is the diagonal. The diagram is
commutative by the very definition of $X_0$. The
top arrow is a homeomorphism and the diagonal is an open
immersion, since $X\to Y$ is \'etale. Therefore we obtain
a diagonal filling as indicated. This proves our claim.
\end{proof}

\rema
If $G\to S$ is infinitesimal (understood, finite locally free)
then as in the proof of Lemma~\ref{fixed_points_ff_subscheme}
we see that the stabilizer of the action is finite infinitesimal
over $X$. Then there is a quotient algebraic space
$q:X\to Y=X/G$ such that $q$ is affine (\cite[Th.~5.3]{Ry13}).
Let $V\to Y$ be an étale surjective map whose source is a disjoint
sum of affine schemes; then $U=V\times_YX\to X$ has the same
properties. Moreover $U$ is endowed with a $G$-action such that
$U\to X$ is equivariant, and from the known affine case we see
that $U^G\hookto U$ is closed. Since $U^G\simeq X^G\times_X U$,
it follows that $X^G\hookto X$ is closed. Note that this proof
may seem simpler than the one given above, but it uses the existence
of quotients.
\xrema

Finally we provide the proof of the existence of FPR-atlases.

\begin{proof}[Proof of Theorem~\ref{adapted-atlas-app}]
Proceeding as in the beginning of the proof of
Theorem~\ref{G_connected_implies_fixed_pts_closed},
we reduce to the case where $S$ is affine and $S,X$ are
of finite type over $\bbZ$. In particular, the set of closed
points $|X^{\cl}|\subset |X|$ is dense. It will thus be enough
to find an affine, $G$-equivariant \'etale neighbourhood
$U(x)\to X$ of each closed point $x\in |X^{\cl}|$ and to
eventually consider $U:=\amalg_{x\in |X^{\cl}|}U(x)$.

For such a closed point $x$ with image $s\in S$, the residue
field extension $\kappa(x)/\kappa(s)$ is finite.
By~\cite{AHR21}, Corollary~20.2 there exists an
affine pointed  scheme $(U_0=\Spec(A),u_0)$ and an \'etale,
$G$-equivariant morphism $(U_0,u_0)\to (X,x)$ which induces
an isomorphism of residue fields
$\kappa(x)\simeq \kappa(u_0)$ and an isomorphism of stabilizers
$G_{u_0}\simeq G_x$. Write
\[
H=D(M/Z)=T\times F
\]
as the product
of a split torus with a finite diagonalizable group scheme.
In case (ii) where $H$ has connected fibres (which means that
$F$ is infinitesimal), it follows from
Theorem~\ref{G_connected_implies_fixed_pts_closed} that
$U_0^H\to X^H\times_X U_0$ is an isomorphism, in other words
$U(x):=U_0$ is the desired $Z$-FPR atlas at $x$.
It remains to consider the case (i) where $X$ is separated.
In this case the stabilizer of the action of $F$ is finite
over $X$ hence there is a quotient $X\to X/F$ which is a
finite morphism (see~\cite[Th.~5.3 and Prop.~4.7]{Ry13}).
Of course, the same is true for the stabilizer of the action
of $F$ on $U_0$ so there exists also a quotient $U_0\to U_0/F$.
Let $U_1:=U_0/F\times_{X/F} X$ and $u_1:=(u_0,x)\in U_1$.
The following properties are seen to hold:
\begin{itemize}
\item[\tiny\textbullet] the scheme $U_1$ is a disjoint sum
of affine $S$-schemes (because this is true for $U_0/F$
and its pullback along the affine map $X\to X/F$),
\item[\tiny\textbullet] the group $G/F$ acts on $U_0/F$ and
on $X/F$, so $G$ acts on the fibred product $U_1$ diagonally
on the factors,
\item[\tiny\textbullet] the map of point stabilizers
$G_{u_1}\to G_x$ is an isomorphism,
\item[\tiny\textbullet] the map $U_1^F\to X^F\times_X U_1$
is an isomorphism.
\end{itemize}
Moreover, by Theorem~\ref{G_connected_implies_fixed_pts_closed}
the morphism $U_1^T\to X^T\times_X U_1$ is an isomorphism.
Since $U_1^H=U_1^T\cap U_1^F$ and $X^H=X^T\cap X^F$, we conclude
that $U(x):=U_1$ is a $Z$-FPR atlas at $x$.

The final claim that $U \to X$ can be chosen quasi-compact
if $X\to S$ is quasi-compact is obvious.
\end{proof}

\end{document}